\documentclass[11pt]{article}
 \usepackage{amsfonts,amssymb,enumerate,epsf}
 \usepackage{graphicx}
\date{}
\usepackage{hyperref}
\usepackage{xcolor}
\usepackage{amsmath}
\usepackage{enumitem}
\usepackage{verbatim}
\usepackage{authblk}

 \usepackage[english]{babel}
 \def\botcaption#1#2{\medskip\centerline{{\scshape #1.}\kern8pt
 {\rm #2}}\bigskip}

 \newtheorem{ittheorem}{Theorem}[section]
 \newtheorem{itlemma}[ittheorem]{Lemma}
 \newtheorem{itproposition}[ittheorem]{Proposition}
 \newtheorem{itdefinition}[ittheorem]{Definition}
 \newtheorem{itremark}[ittheorem]{Remark}
 \newtheorem{itclaim}[ittheorem]{Claim}
 \newtheorem{itcorollary}[ittheorem]{Corollary}
  \numberwithin{equation}{section}

 \newenvironment{theorem}{
 \begin{ittheorem}}{\end{ittheorem}}

 \newenvironment{lemma}{
 \begin{itlemma}}{\end{itlemma}}

 \newenvironment{proposition}{
 \begin{itproposition}}{\end{itproposition}}

 \newenvironment{definition}{
 \begin{itdefinition}}{\end{itdefinition}}

 \newenvironment{remark}{
 \begin{itremark}}{\end{itremark}}

 \newenvironment{claim}{
 \begin{itclaim}}{\end{itclaim}}

\newenvironment{corollary}{
	\begin{itcorollary}}{\end{itcorollary}}

 \newenvironment{proof}{\noindent {\bf Proof.\,}
 }{\hspace*{\fill}$\qed$\medskip}

 \newenvironment{proof*}{\noindent {\bf Proof\,}
}{\hspace*{\fill}$\qed$\medskip}

 \newcommand{\be}[1]{\begin{equation}\label{#1}}
 \newcommand{\ee}{\end{equation}}

 \newcommand{\bl}[1]{\begin{lemma}\label{#1}}
 \newcommand{\el}{\end{lemma}}

 \newcommand{\br}[1]{\begin{remark}\label{#1}}
 \newcommand{\er}{\end{remark}}

 \newcommand{\bt}[1]{\begin{theorem}\label{#1}}
 \newcommand{\et}{\end{theorem}}

 \newcommand{\bd}[1]{\begin{definition}\label{#1}}
 \newcommand{\ed}{\end{definition}}

 \newcommand{\bcl}[1]{\begin{claim}\label{#1}}
 \newcommand{\ecl}{\end{claim}}

 \newcommand{\bp}[1]{\begin{proposition}\label{#1}}
 \newcommand{\ep}{\end{proposition}}

 \newcommand{\bc}[1]{\begin{corollary}\label{#1}}
 \newcommand{\ec}{\end{corollary}}

 \newcommand{\bpr}{\begin{proof}}
 \newcommand{\epr}{\end{proof}}

 \newcommand{\bi}{\begin{itemize}}
 \newcommand{\ei}{\end{itemize}}

 \newcommand{\ben}{\begin{enumerate}}
 \newcommand{\een}{\end{enumerate}}

\def\ss{{\cX^s}}
\def\sm{{\cX^m}}

 \def \ba {\begin{array}}
 \def \ea {\end{array}}

 \def \qed {{\square\hfill}}
 \def \Z {{\mathbb Z}}
 \def \R {{\mathbb R}}
 
 \def \N {{\mathbb N}}
 \def \P {{\mathbb P}}
 \def \E {{\mathbb E}}

 \def \ra {\rightarrow}

 \def \cS {{\cal S}}
 
 \def \cF {{\cal F}}
 \def \cA {{\cal A}}
 \def \cR {{\cal R}}
 
\def \cI {{\cal I}}
\def \cL {{\cal L}}
 
 \def \cC {{\cal C}}
 
 \def \cX {{\cal X}}
 \def \cB {{\cal B}}
 
 \def \cW {{\cal W}}
 
 \def \cP {{\cal P}}
 \def \cL {{\cal L}}
 \def \cV {{\cal V}}

 \def \G {{\Gamma}}
 \def \L {{\Lambda}}

 \def \b {{\beta}}
 \def \e {{\varepsilon}}
 \def \D {{\Delta}}
 \def \r {{\rho}}
\def \m {{\mu}}

 \def \h {{\eta}}
 
 \def \z {{\zeta}}
 \def \g {{\gamma}}
 \def \t {{\tau}}
 \def \o {{\omega}}
 \def \d {{\delta}}
\def \p {{\pi}}

 \def \SES {{\hbox{\footnotesize\rm SES}}}

\def\pieno{{\underline 1}}
\def\vuoto{{\underline 0}}

\begin{document}
	\title{Metastability in a lattice gas \\ with strong anisotropic interactions under Kawasaki dynamics}
	
\author[
{}\hspace{0.5pt}\protect\hyperlink{hyp:email1}{1},\protect\hyperlink{hyp:affil1}{a}
]
{\protect\hypertarget{hyp:author1}{Simone Baldassarri}}

\author[
{}\hspace{0.5pt}\protect\hyperlink{hyp:email2}{2},\protect\hyperlink{hyp:affil1}{a,b}
]
{\protect\hypertarget{hyp:author2}{Francesca R.\ Nardi}}

\affil[ ]{\centering
	\small\parbox{365pt}{\centering
		\parbox{5pt}{\textsuperscript{\protect\hypertarget{hyp:affil2}{a}}}Dipartimento di Matematica e Informatica ``Ulisse Dini", Universit\`{a} degli Studi di Firenze, Firenze, Italy.
	}
}

\affil[ ]{\centering
	\small\parbox{365pt}{\centering
		\parbox{5pt}{\textsuperscript{\protect\hypertarget{hyp:affil2}{b}}}Department of Mathematics and Computer Science, Eindhoven University of Technology, Eindhoven, the Netherlands.
	}
}

\affil[ ]{\centering
	\small\parbox{365pt}{\centering
		\parbox{5pt}{\textsuperscript{\protect\hypertarget{hyp:email1}{1}}}\texttt{\footnotesize\href{mailto:simone.baldassarri@unifi.it}{simone.baldassarri@unifi.it}},
		\parbox{5pt}{\textsuperscript{\protect\hypertarget{hyp:email2}{2}}}\texttt{\footnotesize\href{mailto:francescaromana.nardi@unifi.it}{francescaromana.nardi@unifi.it}},
	}
}

\maketitle

\begin{center}
	{\it To our friend and colleague Carlo Casolo}
	
\end{center}

	\begin{abstract}
		
	In this paper we analyze metastability and nucleation in the context of a local version of the Kawasaki dynamics for the two-dimensional strongly anisotropic Ising lattice gas at very low temperature. Let $\Lambda=\{0,1,..,L\}^2\subset\mathbb{Z}^2$ be a finite box. Particles perform simple exclusion on $\Lambda$, but when they occupy neighboring sites they feel a binding energy $-U_1<0$ in the horizontal direction and $-U_2<0$ in the vertical one. Thus the Kawasaki dynamics is conservative inside the volume $\Lambda$. Along each bond touching the boundary of $\Lambda$ from the outside to the inside, particles are created with rate $\rho=e^{-\Delta\beta}$, while along each bond from the inside to the outside, particles are annihilated with rate $1$, where $\beta$ is the inverse temperature and $\Delta>0$ is an activity parameter. Thus, the boundary of $\Lambda$ plays the role of an infinite gas reservoir with density $\rho$. We consider the parameter regime $U_1>2U_2$ also known as the strongly anisotropic regime. We take $\Delta\in{(U_1,U_1+U_2)}$ and we prove that the empty (respectively full) configuration is a metastable (respectively stable) configuration. We consider the asymptotic regime corresponding to finite volume in the limit of large inverse temperature $\beta$. We investigate how the transition from empty to full takes place. In particular, we estimate in probability, expectation and distribution the asymptotic transition time from the metastable configuration to the stable configuration. Moreover, we identify the size of the \emph{critical droplets}, as well as some of their properties. For the weakly anisotropic model corresponding to the parameter regime $U_1<2U_2$, analogous results have already been obtained. We observe very different behavior in the weakly and strongly anisotropic regimes. We find that the \emph{Wulff shape}, i.e., the shape minimizing the energy of a droplet at fixed volume, is not relevant for the critical configurations.

		 \medskip
		 \noindent
		 
		 {\it AMS} 2010 {\it subject classifications.} 60J10; 60K35; 82C20; 82C22; 82C26\\
		 
		 {\it Key words and phrases.} Lattice gas, Kawasaki dynamics, metastability, critical droplet, large deviations.\\
		 
		 {\it Acknowledgment.} F.R.N. was partially supported by the Netherlands Organisation for Scientific Research (NWO) [Gravitation Grant number 024.002.003--NETWORKS].
		 
	 \end{abstract}

\newpage
\tableofcontents

\section{Introduction}
	\label{ch1}
Metastability is a dynamical phenomenon that occurs when a system is close to first order phase transition, i.e., a crossover that involves a jump in some intrinsic physical parameter such as the energy density or the magnetization. The phenomenon of metastability occurs when a system is trapped for a long time in a state (the metastable state) different from the equilibrium state (the stable state) for specific values of the thermodynamical parameters, and subsequently at some random time the system undergoes a sudden transition from the metastable to the stable state. So we call \emph{metastability} or \emph{metastable behavior} the transition from the metastable state to the equilibrium state. Investigating metastability, researches typically address three main question:
\begin{enumerate}
	\item What are the asymptotic properties of the first hitting time of the stable states for a process starting from a metastable state?
	\item What is the set of critical configurations that the process visits with high probability before reaching the set of stable states?
	\item What is the tube of typical trajectories that the process follows with high probability during the crossover from the metastable states to the stable states?
\end{enumerate}

In this paper we study the metastable behavior of the two-dimensional strongly anisotropic Ising lattice gas that evolves according to Kawasaki dynamics, i.e., a discrete time Markov chain defined by the Metropolis algorithm with transition probabilities given precisely later in (\ref{defkaw}) (see in subsections \ref{S1.1} and \ref{S1.2} for more precise definitions). We consider the local version of the model, i.e., particles live and evolve in a conservative way in a box $\Lambda\subset\mathbb{Z}^2$ and are created and annihilated at the boundary of the box $\Lambda$ in a way that reflects an infinite gas reservoir. More precisely, particles are created with rate $\rho=e^{-\Delta \beta}$ and are annhilated with rate $1$, where $\beta$ is the inverse temperature of the gas and $\Delta>0$ is an activity parameter. When two particles occupy horizontal (resp.\ vertical) neighbouring sites, each one feels a binding energy $-U_1<0$ (resp.\ $-U_2<0$). Without loss of generality we assume $U_1\geq U_2$ and we choose $\Delta\in{(U_1, U_1+U_2)}$, so that the system is in the metastable regime. For this value of the parameters the totally empty (resp.\ full) configuration can be naturally related to metastability (resp.\ stability). We consider the asymptotic regime corresponding to finite volume in the limit of large inverse temperature $\beta$. 

In this work we study the {\it strong anisotropic case}, i.e., the parameters $U_1$, $U_2$ and $\D$ are fixed and such that $U_1>2U_2$ and $\e:=U_1+U_2-\D>0$ sufficiently small. A special feature of Kawasaki dynamics is that in the metastable regime (see (\ref{regime})) particles move along the border of a droplet more rapidly than they arrive from the boundary of the box. More precisely, single particles attached to one side of a droplet tipycally detach before the arrival of the next particle (because $e^{U_1\b}\ll e^{\D\b}$ and $e^{U_2\b}\ll e^{\D\b}$), while bars of two or more particles tipycally do not detach (because $e^{\D\b}\ll e^{(U_1+U_2)\b}$). 

The goal of the paper is to investigate the critical configurations and the tunnelling time between $\underline{0}$ (empty box) and $\underline{1}$ (full box) for this model, answering questions $1$ and $2$ above. In subsection \ref{mainresults} we give four main results: Theorem \ref{metstate} states that the empty box is the metastable state and the full box is the stable state. Theorem \ref{t1} states that the random variable $X_{\b}:=\frac{1}{\b}\log(\t_{\pieno}^{\vuoto})$ converges in probability to $\G$ as $\b$ tends to infinity, where $\G>0$ is a suitable constant that is computed in (\ref{g}) and $\t_{\pieno}^{\vuoto}$ is the first hitting time to $\pieno$ starting from the metastable state $\vuoto$. Moreover, in the same theorem there is also its asymptotic behavior in $\cL^1$ and in law in the limit as $\b\ra\infty$. In particular, after a suitable rescaling, the tunnelling time from $\vuoto$ to $\pieno$ follows an exponential law. This is typical for models where a success occurs only after many unsuccessfull attempts. Theorem \ref{t1'} states that some set of configurations $\cP$, which we define precisely in (\ref{defP}) (see also Figure \ref{fig:P}), has a {\it domino shape} with $l_1\sim2l_2$ (see (\ref{ld0}), (\ref{ld1}) and (\ref{ld2}) for rigorous definitions) and it is a gate for the nucleation, a set with the property that has to be crossed during the transition (see subsection \ref{defind} for the definition of a gate). Theorem \ref{t2} states that with probability tending to 1 the configurations contained in $\cR(2l_2^*-3,l_2^*)$ or $\cR(2l_2^*-1,l_2^*-1)$ are subcritical, in the sense that they shrink to $\vuoto$ before growing to $\pieno$, and those containing $\cR(2l_2^*-2,l_2^*)$ are supercritical, in the sense that they grow to $\pieno$ before shrinking to $\vuoto$.

\begin{figure}[h]
	\begin{center}
		\includegraphics[width=16cm]{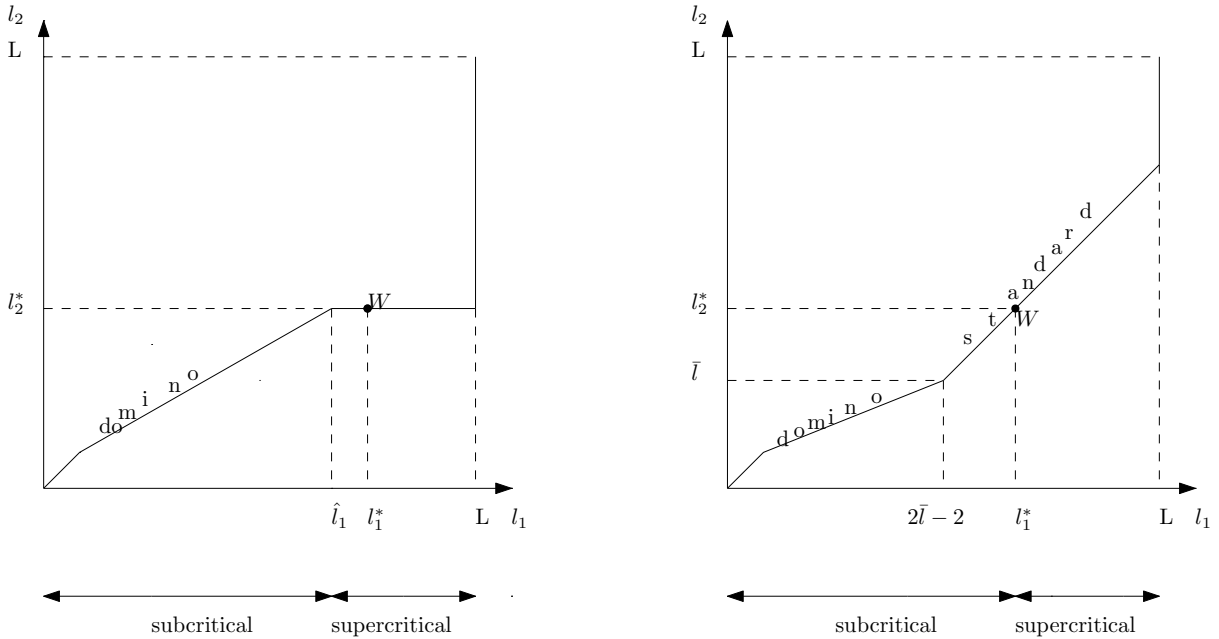}\\
		\caption{Typical path for strong anisotropy (on the left hand-side) and weak anisotropy (on the right hand-side).}
		\label{fig:strongweak}
	\end{center}
\end{figure}

In the regime with exponentially small transition probabilities, it is natural to call {\it Wulff shape} the one minimizing the energy of a droplet at fixed volume. Indeed at low temperature it is possible to show that only the energy is relevant, and the entropy is negligible. An interesting question is then how relevant is the Wulff shape in the nucleation pattern: is the shape of critical configurations Wulff? As mentioned above, if the evolution is according to Kawasaki dynamics, it turns out that particles can move along the border of a droplet more rapidly than they can arrive from the boundary of the container. For this reason particles will be rearranged before the growth of the droplet. Thus one could be tempted to conjecture that this displacement along the border of the growing droplet should establish the equilibrium shape at fixed volume, i.e., the Wulff shape. However, in Section \ref{euristicaa} we give an heuristic discussion based on a careful comparison between time scales of contraction, of growth and of different types of movements on the border of a droplet (see (\ref{sellestaccotrenini})), which indicates that the above conjecture is false. Indeed, as shown in Figure \ref{fig:strongweak} on the left hand-side, the critical droplet is not Wulff (see Theorem \ref{t1'}) and the Wulff shape is a supercritical configuration (see Theorem \ref{t2} and Remark \ref{remwulff}). When growing a nucleus of plus ones the system follows domino shapes up to the critical droplet $\cP$ and thus in this case the ratio between the side lengths is approximatively 2. In the Wulff shape the ratio between the side lengths is of order $\frac{U_1}{U_2}$, so we come to the conclusion that the critical droplets are not Wulff shape. In Section \ref{euristicaa}, we let $\hat{l}_1$, $l_2^*$ be the horizontal and vertical sides, respectively, of the critical droplet ($\hat{l}_1\in\{2l_2^*-2,2l_2^*-1\}$, see (\ref{defP1}) and (\ref{defP2})). In the strongly anisotropic case, the supercritical growth follows a sequence of rectangles with $l_2=l_2^*$ and $l_1=\hat{l}_1+m$, with $m=1,2,...$ up to $l_1=L$, the side of the container. During this epoch, the nucleation pattern crosses the Wulff shape with sides $(l_1^*,l_2^*)$ (see (\ref{defbarl}) and Remark \ref{remwulff}). Finally, after the formation of a strip $l_2^*\times L$ the system starts growing in the vertical direction up to the full configuration. Similarly, for any anisotropic Glauber dynamics the critical configurations are not Wulff-shaped and the tube of typical paths crosses the Wulff shape only during the supercritical growth (see \cite{KOd}). Indeed the tube of typical paths evolves along squared-shaped configurations in the subcritical part, along {\it horizontal-growing rectangles} until they wrap around the torus, and then stripes growing in vertical direction up to the configurations with all pluses.
	
On the one hand, if we change the definition of a gate not imposing that the energy of its configurations is $\G$, but requiring only that every optimal path must cross it, heuristically we have that the Wulff shape has a gate property. On the other hand, we want to underline that Theorem \ref{t2} can not be extended for the Wulff shape, since it is not true that all the configurations with rectangle smaller than Wulff shape are subritical (see Remark \ref{remwulff}). Moreover, Theorem \ref{t1} can not be adapted to the Wulff shape because in our paper we fix the values of the parameters $U_1$, $U_2$ and $\D$ such that (\ref{regime}) holds. The most interesting results are obtained in the cases when $\e$ is small, in which we have that the critical configurations are large. In our regime the energy of the Wulff shape is different from $\G$ and this holds also in the case in which $\D\ra U_1+U_2$ that corresponds to $\e\ra0$.

\subsection{Comparison with models subject to Kawasaki dynamics}

In this subsection we make a comparison between the model considered in this paper and other models that also consider Kawasaki dynamics and were already studied in literature.

The bidimensional isotropic case $U_1=U_2$ has already been studied using the {\it pathwise approach} in \cite{HOS} with results concerning question 1 giving estimates in probability, law and distribution and, concerning question 2, giving intrinsically the critical configurations without their geometrical description. In \cite{GOS} the authors investigated question 3 identifying the tube of typical trajectories, again using the pathwise approach. For the three-dimensional case, in \cite{HNOS} there are results concerning questions 1 and 2. It is interesting that, concerning the asymptotic expectation of the tunnelling time as in (\ref{theorem3}), using the pathwise approach, it is not possible to distinguish the presence of a certain function $f(\b)$ such that $\log f(\b)/\b\ra0$ in the limit as $\b\ra+\infty$ and $\E_{\vuoto}\t_{\pieno}=f(\b)e^{\G\b}(1+o(1))$, or the presence of a constant factor. To this end, a more detailed study of the so-called {\it pre-factor} $f(\b)$ is given in \cite{BHN} for two and three dimensions, using the {\it potential theoretic approach}. In \cite{BHN} the authors estimated the constant pre-factor  and found that it does not depend on the parameter $\b$, but on the size of the box and the cardinality of the set of critical droplets with size $l_c$. These estimates of the pre-factor are possible once the geometrical description of the critical configurations and of its neighborhood are found. See also \cite{GoL}, where \cite{BL1} is applied to derive again results for this model. Furthermore, in the three-dimensional case similar results are obtained but with less control over the geometry and the constant. Since in the isotropic models the Wulff-shape concides with the critical shape, it is not possible to distinguish among them. This motivates together with applications the study of anisotropic models.

The weak anisotropic case $U_2<U_1<2U_2-\e$ has already been studied in \cite{NOS}. We observe very different behavior with respect to (w.r.t.) strong anisotropy. Indeed for weak anisotropy, as we can see in Figure \ref{fig:strongweak} on the right hand-side, after the growth along \emph{domino} shapes with $l_1\sim2l_2$ (see (\ref{ld0}),(\ref{ld1}) and (\ref{ld2})) during the early stage of nucleation, the nucleation pattern consists of a growing sequence of a certain class of rectangles, called \emph{standard rectangles} (see (\ref{ls0}) and (\ref{ls1}) for a rigorous definition). In this regime we have that $l_1^*$ and $l_2^*$ are the critical sizes (see (\ref{defbarl}) and (\ref{critinteri}) for the definition). In \cite{NOS} it is proved that the critical droplet is close to the Wulff shape, whereas during the other stages of nucleation the shape of the growing droplet is not Wulff, but standard (see \cite[Theorems 2,3]{NOS}). This argument leads to say that in both strong and weak anisotropy the Wulff shape is not relevant in the nucleation pattern as for anisotropic Glauber dynamics (see \cite{KOd}). Our choice to work with Kawasaki dynamics rather than Glauber in this very low temperature regime is a first step in showing the robustness of the argument rooted in the dynamical nature of metastable systems. The locally conservative dynamics and the movement of particles along the border of the droplet give a regularization effect. Surprisingly, as mentioned above, this effect does not drive the nucleation process along Wulff-shaped configurations especially in the assumption of strong anisotropy. More precisely, for weak anisotropy the critical configuration is Wulff-shaped but the tube of typical paths is evolving via domino and standard configurations while for strong anisotropy the critical configurations are not Wulff-shaped and the tube of typical paths is evolving via domino rectangles in the subcritical part and via {\it horizontal-growing rectangles} in the supercritical part.

Results similar to the ones obtained in this paper were given in \cite{HNT2,HNT1,HNT3} where two types of particles are present in $\L$. In particular, the authors analyzed the two-dimensional lattice gas subject to Kawasaki dynamics, where neighboring particles have negative binding energy if and only if their types are different. The authors obtained results regarding the identification of the critical droplets and their geometrical properties, i.e., question 2, in \cite{HNT3}. With this knowledge they studied the transition time from the metastable state to the stable state in \cite{HNT2} in law and in distribution (question 1), using the potential-theoretic approach. In particular, they were able to identify the pre-factor.

It turns out that a complete description of the tube of typical trajectories (question 3), as given in \cite{KOd} for the anisotropic Ising model evolving under Glauber dynamics, is much more complicated when we consider Kawasaki dynamics. Using Kawasaki dynamics, the tube of typical trajectories is analyzed only in \cite{GOS} for the two-dimensional isotropic case. There are no known results for three dimensions, either for the anisotropic model or for the two-particle-types model. We remark that in many previous papers (\cite{{AC},{CN},{CO},{HNOS},{KOd},{KOs},{NO}}) the asymptotic of the tunnelling time and the tube of typical trajectories realizing the transition were treated simultaneously by exploiting a detailed control of the energy landscape in connection with the paths allowed by the dynamics.

\subsection{State of the art}
\label{art}
A mathematically description was first attempted in \cite{LPr,KUH} inspired on Gibbsian equilibrium Statistical Mechanics. A more faithful approach, known as \emph{pathwise approach}, was initiated in 1984~\cite{CGOV} and was developed in \cite{OS,OS2,OV}. This approach focuses on the {\it dynamics} of the transition from metastable to stable state. Independently, a graphical approach was introduced in \cite{CC} and later used for Ising-like models \cite{CaTr}. With the pathwise approach they obtained a detailed description of metastable behavior of the system and it made possible to answer the three questions of metastability. By identifying the most likely path between metastable states, the time of the transition and the tube of typical trajectories can be determined. A modern version of the pathwise approach containing the information about time and critical droplets disentangled w.r.t.\ the tube of typical trajectories can be found in \cite{{MNOS},{CNbc},{CNS2},{BNZ}}. This approach developed over the years has been extensively applied to study metastability in Statistical Mechanics lattice models. In this context, this approach and the one that follows (\cite{BEGK,MNOS,OV}) have been developed with the aim of finding answers valid with maximal generality and to reduce as much as possible the number of model dependent inputs necessary to describe the metastable behavior of any given system.

Another approach is the \emph{potential-theoretic approach}, initiated in \cite{BEGK}. We refer to \cite{BH} for an extensive discussion and applications to different models. In this approach, the metastability phenomenon is interpreted as a sequence of visits of the path to different metastable sets. This method focuses on a precise analysis of hitting times of these sets with the help of \emph{potential theory}. In the potential-theoretic approach the mean transition time is given in terms of the so-called \emph{capacities} between two sets. Crucially capacities can be estimated by exploiting powerful variational principles. This means that the estimates of the average crossover time that can be derived are much sharper than those obtained via the pathwise approach. The quantitative success of the potential-theoretic approach is however limited to the case of reversible Markov processes. 

These mathematical approaches, however, are not completely equivalent as they rely on different definitions of metastable states (see \cite[Section 3]{CNbc} for a comparison) and thus involve different properties of hitting and transition times. The situation is particularly delicate for evolutions of infinite-volume systems, for irreversible systems, and degenerate systems, i.e., systems where the energy landscape has configurations with the same energy (as discussed in \cite{CNbc,CNS2,CNS2017}). More recent approaches are developed in \cite{BL1,BL2,BiGa}.

Statistical mechanical models for magnets deal with dynamics that do not conserve the total number of particles or the total magnetization. They include single spin-flip Glauber dynamics and many probabilistic cellular automata (PCA), that is parallel dynamics. The pathwise approach was applied in finite volume at low temperature in \cite{CGOV,NevSchbehavdrop,CaTr,KOd,KOs,CO,NO,CL} for single-spin-flip Glauber dynamics and in \cite{CN,CNSp,CNS22,CNS2016} for parallel dynamics. The potential theoretic approach was applied to models at finite volume and at low temperature in \cite{BM,BHN,HNT1,HNT2,NS1,HNTA2018}. The more involved infinite volume limit at low temperature or vanishing magnetic field was studied in \cite{DS1,DS2,S2,SS,MO1,MO2,HOS,GHNOS,GN,BHS,CeMa,GMV} for Ising-like models under single-spin-flip Glauber and Kawasaki dynamics.

\subsection{Outline of the paper}

The outline of the paper is as follows. In Section \ref{S1} we define the model, give some definitions in order to state our main theorems (see Theorems \ref{metstate}, \ref{t1}, \ref{t1'} and \ref{t2}), give a comparison between strong and weak anisotropy and a heuristic discussion of the dynamics. In Section \ref{depres} we obtain our main model-dependent results regarding the metastable and stable states (Theorem \ref{metstate}), tunnelling time (Theorem \ref{t1}) and the gate of the transition (Theorem \ref{t1'}). In Section \ref{proofs} we give the proof of Theorem \ref{3.1}, which consists in a careful analysis of the minimal energy along all the possible communicating configurations from a particular set $\cB$ to $\cB^c$ (see (\ref{defB}) for the precise definition). In Section \ref{recurrence} we prove an important result that allows us to deduce Theorem \ref{t2}. In the Appendix we give additional explicit computations.

\section{Model and results}
\label{S1}

\subsection{Definition of the model}
\label{S1.1}

Let $\Lambda=\{0,..,L\}^2\subset \Z^2$ be a finite box centered at the
origin. The side length $L$ is fixed, but arbitrary, and later we will require $L$ to be sufficiently large. Let

\be{inbd}
\partial^- \L:= \{x\in\L\colon ~\exists\; y \notin\L\colon ~|y-x|=1\},
\ee

\noindent be   the interior  boundary of $ \Lambda$ and let
$ \Lambda_0:= \Lambda\setminus\partial^- \Lambda$ be the interior of $\L$.
With each $x\in \Lambda$ we associate an occupation variable
$\eta(x)$, assuming values 0 or 1. A lattice configuration is
denoted by $\eta\in {\cal X} =\{ 0,1\} ^{ \Lambda }$.

\noindent
Each configuration $\h\in \cX$ has an energy given by the following Hamiltonian:

\be{hamilt} H(\eta):= -U_1 \sum_{(x,y)\in   \Lambda_{0,h}^{*}}
\eta(x)\eta(y) -U_2\sum_{(x,y)\in   \Lambda_{0,v}^{*}} \eta(x)
\eta(y)+ \D \sum _{x\in \L} \eta (x), \ee

\noindent
where $ \Lambda_{0,h}^{*}$ (resp.\ $ \Lambda_{0,v}^{*}$) is the set of
the horizontal (resp.\ vertical) unoriented  bonds joining nearest-neighbors points in
$ \Lambda_0$. Thus the interaction is acting only inside $
\Lambda_0$; the binding energy associated to a horizontal
(resp.\ vertical) bond is $-U_1<0$ (resp.\ $-U_2<0$). We may assume without
loss of generality that $U_1\ge U_2$. (Note that $H-\D\sum _{x\in
	\partial^-\L} \eta (x) $  can be viewed as the
Hamiltonian, in lattice gas variables, for an Ising  system
enclosed in $\L_0$, with $0$ boundary conditions.)

The grand-canonical Gibbs measure associated with $H$ is

\be{misura} \m(\eta):= {  e^{- \b H(\eta) }\over Z} \qquad \h\in
\cX, \ee

\noindent
where

\be{partfunc} Z:=\sum_{\eta\in {\cal X}}e^{-\b H(\eta)} \ee

is the so-called {\it partition function}.

\subsection{Local Kawasaki dynamics}
\label{S1.2}

Next we define Kawasaki dynamics on $\L$ with boundary
conditions that mimic the effect of an infinite gas reservoir
outside $\L$ with density $ \r = e^{-\D\b}.$ Let $b=(x \to y)$ be
an oriented bond, i.e., an {\it ordered} pair of nearest neighbour
sites, and define

\be{Loutindef}
\ba{lll}
\partial^* \L^{out} &:=& \{b=(x \to y)\colon \; x\in\partial^- \L,
y\not\in\L\},\\
\partial^* \L^{in}  &:=& \{b=(x \to y)\colon \; x\not\in
\L, y\in\partial^-\L\},\\
\L^{*, orie} &:=& \{b=(x \to y)\colon \; x,y\in\L\},
\ea
\ee

\noindent and put $ \bar\L^{*, orie}:=\partial^* \L ^{out}\cup
\partial^* \L ^{in}\cup\L^{*,\;orie}$.
Two configurations $  \eta,
\eta'\in {\cal X}$ with $ \eta\ne \eta'$ are said to be {\it
	communicating states} if there exists a bond
$b\in  \bar\L^{*,orie}$ such that $ \eta' = T_b \eta$, where $T_b   \eta$ is the
configuration obtained from $ \eta$ in any of these ways:

\begin{itemize}
	
	\item
	for $b=(x \to y)\in\L^{*,\;orie}$, $T_b \eta$ denotes the
	configuration obtained from $ \eta$ by interchanging particles
	along $b$:
	
	\be{Tint}
	T_b \h(z) =
	\left\{\ba{ll}
	\h(z) &\mbox{if } z \ne x,y,\\
	\h(x) &\mbox{if } z = y,\\
	\h(y) &\mbox{if } z = x.
	\ea
	\right.
	\ee
	
	\item
	For  $b=(x \to y)\in\partial^*\L^{out}$ we set:
	
	\be{Texit}
	T_b \h(z) =
	\left\{\ba{ll}
	\h(z) &\mbox{if } z \ne x,\\
	0     &\mbox{if } z = x.
	\ea
	\right.
	\ee
	
	\noindent
	This describes the annihilation of particles along the border;
	
	\item
	for  $b=(x  \to y)\in\partial^*\L^{in}$ we set:
	
	\be{Tenter}
	T_b \h(z) =
	\left\{\ba{ll}
	\h(z) &\mbox{if } z \ne y,\\
	1     &\mbox{if } z=y.
	\ea
	\right.
	\ee
	
	\noindent
	This describes the creation of particles along the border.
	
\end{itemize}

The Kawasaki dynamics is  the discrete time Markov chain
$(\eta_t)_{t\in \mathbb{N}}$ on state space $ {\cal X} $ given by
the following transition  probabilities: for  $  \eta\not= \eta'$:

\be{defkaw}
\P( \eta,  \eta'):=\left\{\ba{ll}
{ |\bar\L^{*,\;orie}|}^{-1} e^{-\b[H( \eta') - H( \eta)]_+}
&\mbox{if }  \exists b\in \bar\L ^{*, orie}: \eta' =T_b \eta   \\
0   &\mbox{ otherwise }  \ea \right.
\ee

\noindent
and $\P(\h,\h):=1-\sum_{\h'\neq\h}\P(\h,\h')$, where $[a]_+ =\max\{a,0\}$.
This is a standard Metropolis dynamics with an open boundary: along each
bond touching $\partial^-\L$ from the outside, particles are created with
rate $\rho=e^{-\D\b}$ and are annihilated with rate 1, while inside $\L_0$ particles
are conserved.
Note that an exchange of occupation
numbers $\h(x)$ for any $x$ inside the ring $ \L\setminus  \L_0$
does not involve any change in energy.\par

\br{p1}
The stochastic dynamics defined by
(\ref{defkaw}) is reversible w.r.t. Gibbs measure (\ref{misura}) corresponding
to $ H$.
\er

\subsection{Definitions and notations}
We will use italic capital letters for subsets of $\L$, script
capital letters for subsets of $\cX$, and boldface capital letters for events
under the Kawasaki dynamics. We use this convention in order to keep the
various notations apart. We will denote by $\P_{\h_0}$ the probability law of the Markov
process  $(\eta_t)_{t\geq 0}$ starting at $\h_0$ and by
$\E_{\h_0}$ the corresponding expectation.

In order to formulate our main results in Theorem \ref{metstate}, Theorem \ref{t1}, Theorem \ref{t1'} and Theorem \ref{t2},
we first need some definitions.

\subsubsection{Model-independent definitions and notations}
\label{defind}
\medskip\noindent
{\bf 1. Paths, boundaries and hitting times.}
\bi

\item
A {\it path\/} $\o$ is a sequence $\o=\o_1,\dots,\o_k$, with
$k\in\N$, $\o_i\in\cX$ and $P(\o_i,\o_{i+1})>0$ for $i=1,\dots,k-1$.
We write $\o\colon\;\h\to\h'$ to denote a path from $\h$ to $\h'$,
namely with $\o_1=\h,$ $\o_k=\h'$. A set
$\cA\subset\cX$ with $|\cA|>1$ is {\it connected\/} if and only if for all
$\h,\h'\in\cA$ there exists a path $\o\colon\;\h\to\h'$ such that $\o_i\in\cA$
for all $i$.

\item[$\bullet$]
Given a non-empty
set $\cA\subset\cX$, define  the {\it first-hitting time of} $\cA$
as

\be{tempo}
\t_{\cA}:=\min \{t\geq 0\colon\, \eta_t \in \cA \}.
\ee
\ei
\noindent

\medskip
\noindent
{\bf 2. Min-max and gates}
\bi
\item 
The {\it bottom} $\cF(\cA)$ of a  non-empty
set $\cA\subset\cX$ is the
set of {\it global minima} of the Hamiltonian $H$ in  $\cA$:

\be{Fdef}
\cF(\cA):= \{\h\in\cA\colon\; H(\h)=\min_{\z\in\cA} H(\z)\}.
\ee
With an abuse of notation, for a set $\cA$ whose points have the same
energy we denote this energy by
$H(\cA)$.

\item 
The {\it communication height\/} between a pair $\h,\h'\in\cX$ is

\be{comh}
\Phi(\h,\h'):= \min_{\o\colon\;\h\to\h'}
\max_{\z\in\o} H(\z).
\ee

\item
We call
{\it stability level} of
a state $\z \in \cX$
the energy barrier

\be{stab}
V_{\z} :=
\Phi(\z,\cI_{\z}) - H(\z),
\ee

\noindent where $\cI_{\z}$ is the set of states with
energy below $H(\z)$:

\be{iz} \cI_{\z}:=\{\eta \in \cX:\; H(\eta)<H(\z)\}. \ee

\noindent We set $V_\z:=\infty$ if $\cI_\z$ is empty.

\noindent

\item
We call set of {\it $V$-irreducible states}
the set of all states with stability level larger than  $V$:

\be{xv} \cX_V:=\{\h\in\cX:\quad V_{\h}>V\}. \ee

\item
The set of {\it stable states} is the set of the global minima of
the Hamiltonian:

\be{st.st.} \cX^s:=\cF(\cX). \ee

\item
The set of {\it metastable states} is given by

\be{st.metast.} \sm:=\{\h\in\cX:\quad
V_{\h}=\max_{\z\in\cX\backslash \ss}V_{\z}\}. \ee
\noindent
We denote by $\G_m$ the stability level of the states in $\cX^m$.

\item 
We denote by $(\h\to\h')_{opt} $ the {\it set of optimal paths\/}, i.e., the set of all
paths from $\h$ to $\h'$ realizing the min-max in $\cX$, i.e.,

\be{optpath}
(\h\to\h')_{opt}:=\{\o:\h\to\h'\; \hbox{such that} \; \max_{\xi\in\o} H(\xi)=  \Phi(\h,\h') \}.
\ee

\item
The set of {\it minimal saddles\/} between
$\h,\h'\in\cX$
is defined as

\be{minsad}
\cS(\h,\h'):= \{\z\in\cX\colon\;\; \exists\o\in (\h\to\h')_{opt},
\ \o\ni\z \hbox{ such that} \max_{\xi\in\o} H(\xi)= H(\z)\}.
\ee

Given two non-empty sets $\cA,\cB\subseteq\cX$, put

\be{SAB}
\cS(\cA,\cB):= \bigcup_{{\h\in\cA,\,\h'\in\cB\colon \atop
		\Phi(\h,\h') = \Phi(\cA,\cB)}} \cS(\h,\h').
\ee

\item
Given a pair $\h,\h'\in\cX$,
we say that $\cW\equiv\cW(\h,\h')$ is a {\it gate\/}
for the transition $\h\to\h'$ if $\cW(\h,\h')\subseteq\cS(\h,\h')$
and $\o\cap\cW\neq\emptyset$ for all $\o\in (\h\to\h')_{opt}$.

\item
We say that $\cW$ is a {\it minimal gate\/} for the transition
$\h\to\h'$ if it is a gate and for any $\cW'\subset \cW$ there
exists $\o'\in (\h\to \h')_{opt}$ such that $\o'\cap\cW'
=\emptyset$. In words, a minimal gate is a minimal subset of $\cS(\h,\h')$ by inclusion that is visited by all optimal paths.

\ei

\subsubsection{Model-dependent definitions and notations}
\label{moddef}
We briefly give some model-dependent definitions and notations in order to state our main theorems. For more details see subsection \ref{extmoddef}.

\bi

\item[$\bullet$]
For $x\in\L_0$, let $ \hbox{nn}(x):=\{ y\in \L_0\colon\;|y-x|=1\}$ be the set of
nearest-neighbor sites of $x$ in $\L_0$.
\item[$\bullet$]
A {\it free particle\/} in $\h\in\cX$ is a site $x$ such that either $x\in\eta\cap\partial^-\L$
or $x\in \eta\cap\L_0$, and $\sum_{y\in nn(x)\cap\L_0}\eta(y)$ $=0$.

We denote by $\h_{fp}$ the union  of free particles in $\partial^-
\L$ and  free particles in $\L_0$ and by $\h_{cl}$ the clusterized
part of $\h$

\be{hcl}\h_{cl} :=\h\cap \L_0 \setminus \h_{fp}.\ee

We denote by $|\h_{fp}|$ the number of free particles in $\h$ and by $|\h_{cl}|$ the cardinality of the clusterized part of $\h$.

\item[$\bullet$]
Given a configuration $\h \in\cX$, consider the set $C(
\h_{cl}) \subset \R^2$ defined as the union of the $1\times 1$
closed squares centered at the occupied sites of $ \h_{cl}$ in $
\L_0$.

\item[$\bullet$]
For $\h\in\cX$, let $|\eta|$ be the number of particles in $\eta$,
$\g(\h)$ the Euclidean boundary of $C(\h_{cl})$, $\g(\h)=\partial
C(\h_{cl})$; we denote by $g_1(\eta)$
(resp. $g_2(\eta)$) one half of the horizontal (resp. vertical) length
of $\g(\h)$.

\item[$\bullet$]
Let $p_1(\h)$
and $p_2(\h)$ be the total lengths of  horizontal and vertical
projections of $ C(\h_{cl})$ respectively.

\item[$\bullet$]
We define $g'_i(\h):= g_i(\h)- p_i(\h)\ge 0$; we call {\it monotone} a
configuration such that $g_i(\h)= p_i(\h)$ for $i=1,2$.

\item[$\bullet$]
We write
\be{defdep1}
\ba{lll}
s(\h)&:=& p_1(\h)+p_2(\h),\\
v(\h)&:=& p_1(\h)p_2(\h)- |\h_{cl}|,\\
n(\h)&:=& |\h_{fp}|.
\ea
\ee

\item[$\bullet$]
We denote by $\cR (l_1,l_2)$ the set of configurations whose
single contour is a rectangle $R(l_1,l_2)$, with $l_1,l_2\in\N$.
For any $\h,\h'\in\cR (l_1,l_2)$ we have immediately:

\be{enerettan} H(\h)=H(\h')=H(\cR(l_1,l_2))={U_1}l_2 +{U_2} l_1
-\varepsilon l_1 l_2, \ee
where
\be{defepsilon}\e:= U_1+U_2 -\D. \ee

\item[$\bullet$]
Let $l_2\geq2$. A rectangle $R(l_1,l_2)$ with $l_1=2l_2$ is called \emph{$0$-domino rectangle} and in this case we have $[s]_3=[0]_3$ with $s=l_1+l_2$, where for $x\in \Z$, $n\in \N$ we denote $[x]_n:=x \hbox { mod } n$. Thus we define the set of $0$-domino rectangles as $\mathcal{R}^{0-dom}(s)=\mathcal{R}(l_1(s),l_2(s))$ with
\be{ld0}
l_1(s):=\frac{2s}{3}, \  \ \ \ \ 
l_2(s):=\frac{s}{3} \ \ \ \hbox{for} \ [s]_3=[0]_3.
\ee
If $l_1=2l_2-2$, we have $[s]_3=[1]_3$, so we define the set of \emph{$1$-domino rectangles} as $\cR^{1-dom}(l_1(s),l_2(s))$ with
\be{ld1}
l_1(s):=\frac{2s-2}{3}, \  \ \ \ \ 
l_2(s):=\frac{s+2}{3} \ \ \ \hbox{for} \ [s]_3=[1]_3.
\ee
If $l_1=2l_2-1$, we have $[s]_3=[2]_3$, so we define the set of \emph{$2$-domino rectangles} as $\mathcal{R}^{2-dom}(l_1(s),l_2(s))$ with
\be{ld2}
l_1(s):=\frac{2s-1}{3}, \  \ \ \ \ 
l_2(s):=\frac{s+1}{3} \ \ \ \hbox{for} \ [s]_3=[2]_3.
\ee

\ei
\subsection{Main results}
\label{mainresults}
Let

\be{squaredefempty}
\vuoto:= \{\h\in\cX\colon\; \h(x)=0 \; \forall x \in \L\}
\ee

\noindent
be the  configuration empty. By (\ref{hamilt}) and (\ref{squaredefempty}) we have that $H(\vuoto)=0$. Let
\be{squaredeffull}
\pieno:= \{\h\in\cX\colon\; \h(x)=1 \quad\forall x \in
\L_0,\ \h(x)=0 \quad\forall x \in\L\backslash \L_0 \}
\ee

\noindent
be the  configuration that is full in $\L_0$ and empty in  $\L\backslash\L_0$.

\bt{metstate} If the side $L$ of the box $\Lambda$ is sufficiently large, then $\G_m=V_{\underline{0}}=\Phi(\underline{0}, \underline{1})=\G$ and $\underline{1}=\mathcal{X}^s$ and $\underline{0}=\mathcal{X}^m$.
\et

The proof of Theorem \ref{metstate} is analogue to the one in \cite[Proposition 15]{NOS} for a different value of $L$ (say $L>\lceil\frac{U_1+U_2}{\e}\rceil$) once we have given a specific upper bound $V^*$ for the stability levels of the configurations different from $\vuoto$ and $\pieno$ (see Proposition \ref{V^*}) and proved $\Phi(\vuoto,\pieno)=\G$ (see Corollary \ref{phi}).

\medskip
\noindent
We define the critical vertical length

\be{critinteri}  l_2^*:=\left\lceil {U_2 \over U_1+U_2 -\D}\right\rceil={U_2 \over U_1+U_2 -\D}+\d,
\ee

\noindent
where $\lceil \;\rceil$ denotes the integer part plus 1, and $0<\d<1$ is fixed. Furthermore we set the critical value of $s$ and the critical configurations $\cP$ as
\be{defs^*} s^*:=3l_2^*-1, \ee
\be{defP}
\cP:=\cP_{1}\cup\cP_{2}, \ee

\noindent
where
\be{defP1}
\ba{ll}
\cP_{1}:= \{\h&:\, n(\h)= 0,\, v(\h)=\ell_1(s^*)-1,\,
\h_{cl} \hbox{ is connected},\, g_1'(\h)=0,\,g_2'(\h)=1,\\
&\hbox{ with circumscribed rectangle in }
\cR(\ell_1(s^*),\ell_2(s^*))\},
\ea
\ee
\be{defP2}
\ba{ll}
\cP_{2}:= \{\h&:\, n(\h)= 1,\, v(\h)=\ell_2(s^*-1)-1,\,
\h_{cl} \hbox{ is connected},\, \hbox{monotone},\\
&\hbox{ with circumscribed rectangle in }
\cR(\ell_1(s^*-1),\ell_2(s^*-1))\},
\ea
\ee
\noindent 
with $l_i(s),\, i=1,2$ defined as in (\ref{ld0}), (\ref{ld1}) and
(\ref{ld2}) (recall (\ref{defdep1})). See Figure \ref{fig:P} for an example of configurations in $\cP$.

From (\ref{energeta}) below, it follows that $H(\h)$ is constant on $\cP$. We write

$$\G := H(\cP)-H(\vuoto)= H(\cP)=H(\cR(2l_2^*-2,l_2^*))+\e(l_2^*-1)+\D=$$
\be{g}
=H(\cR(2l_2^*-1,l_2^*-1))+\D-U_2+U_1. \ee

\setlength{\unitlength}{1.1pt}
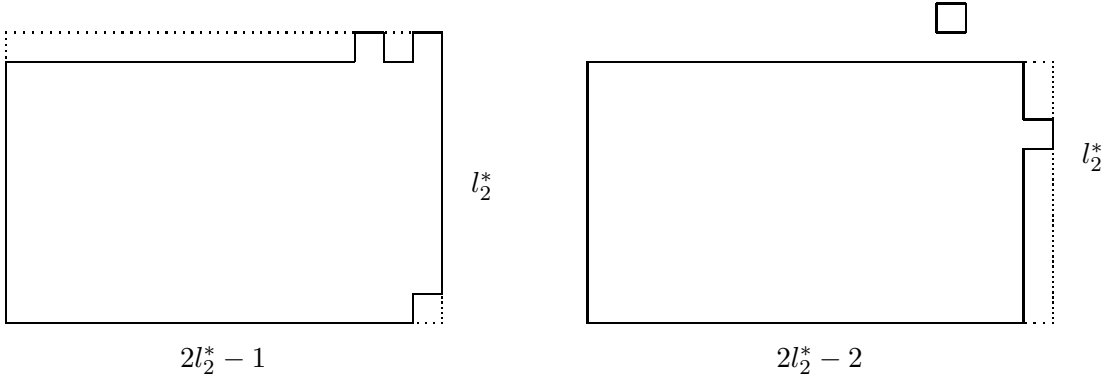
\begin{figure}
	\begin{picture}(400,90)(0,30)
	\thinlines
	\qbezier[51](20,0)(100,0)(170,0)
	\qbezier[51](20,100)(100,100)(170,100)
	\qbezier[51](20,0)(20,45)(20,100)
	\qbezier[51](170,0)(170,45)(170,90)
	\thinlines
	\put(20,0){\line(1,0){140}}
	\put(20,90){\line(1,0){120}}
	\put(20,0){\line(0,1){90}}
	\put(150,90){\line(1,0){10}}
	\put(160,90){\line(0,1){10}}
	\put(160,100){\line(1,0){10}}
	\put(170,90){\line(0,1){10}}
	\put(170,10){\line(0,1){90}}
	\put(160,10){\line(1,0){10}}
	\put(160,0){\line(0,1){10}}
	\put(140,90){\line(0,1){10}}
	\put(150,90){\line(0,1){10}}
	\put(140,100){\line(1,0){10}}
	\thinlines  \put(180,45){$l_2^*$}
	\put(80,-15){$2l_2^*-1$}
	\thinlines
	\thinlines
	\qbezier[51](220,0)(300,0)(380,0)
	\qbezier[51](220,90)(300,90)(380,90)
	\qbezier[51](220,0)(220,45)(220,90)
	\qbezier[51](380,0)(380,45)(380,90)
	\thinlines
	\put(220,0){\line(1,0){150}}
	\put(220,90){\line(1,0){150}}
	\put(220,0){\line(0,1){90}}
	\put(370,0){\line(0,1){60}}
	\put(370,70){\line(0,1){20}}
	\put(370,60){\line(1,0){10}}
	\put(380,60){\line(0,1){10}}
	\put(370,70){\line(1,0){10}}
	\put(340,100){\line(1,0){10}}
	\put(350,100){\line(0,1){10}}
	\put(340,100){\line(0,1){10}}
	\put(340,110){\line(1,0){10}}
	\thinlines  
	\put(390,55){$l_2^*$}
	\put(285,-15){$2l_2^*-2$}
	\end{picture}
	\vskip 1.5 cm
	\caption{Configurations in $\cP$: on the left hand-side is represented a configuration in $\cP_1$ and on the right hand-side a configuration in $\cP_2$}
	\label{fig:P}
\end{figure}

\medskip
\noindent
The behavior of the model strongly depends on the different values of the parameters. We will not consider all the possible regimes and we will not be interested in characterizing the broadest parameter regime for which our results hold. We will assume

\be{regime} 0< \e\ll U_2\qquad \hbox{ and  }\qquad U_1> 
2U_2, \ee

\noindent where $\ll$ means sufficiently smaller; for instance
$\e\le{U_2\over 100}$ is enough.

The main results about the asymptotics of the tunnelling time is contained in the following:

\bt{t1} Let $U_1, U_2, \D$ be  such that $U_2/(U_1+U_2-\D)$ is not
integer and (\ref{regime}) holds. Let $\L$ be a box with side
$L+2$. For $L$ sufficiently large and for any $\d>0$,

\be{lim} \lim_{\b\to\infty}P_{\vuoto}\big( e^{\b(\G-\d)} \leq
\t_{\pieno} \leq e^{\b(\G+\d)}\big)\, =\, 1, \ee

\be{theorem3}
\lim_{\b\to\infty}{1\over \b}\log  \E_{\vuoto}\t_{\pieno}=\G.
\ee

\noindent
Moreover, letting
$T_\b:=
\inf
\{n \ge 1:\, \P_{\vuoto}({\t_{\pieno} \le n})\ge
1-e^{-1}\},$
we have

\be{th2}
\lim_{\b \to \infty} \P_{\vuoto}(\t_{\pieno} > t T_\b) = e^{-t},
\ee

\noindent and

\be{th2'}
\lim_{\b \to \infty} {\E_{\vuoto} (\t_{\pieno})\over T_\b} =1.
\ee

\et

\noindent
In words, Theorem \ref{t1} says that

\bi

\item[--]
(\ref{lim}): For $\b\to\infty$  the transition
time from $\vuoto$ to $\pieno$ behaves asymptotically, in
probability, as  $e^{\G\b+o(\b)}$.

\item[--]
(\ref{theorem3}) and (\ref{th2}):
The mean value of the transition time from $\vuoto$ to $\pieno$ is asymptotic to $e^{\G\b}$ as $\b\ra\infty$. Moreover, the rescaled transition time converges to an exponential distribution.

\ei

\noindent
We refer to subsection \ref{timegate} for the proof of Theorem \ref{t1}. The main results about the gate to stability are contained in the following:

\bt{t1'} Let $U_1, U_2, \D$ be  such that $U_2/(U_1+U_2-\D)$ is
not integer and (\ref{regime}) holds. Let $\L$ be a box with side
$L+2$. For $L$ sufficiently large, $\cP$ is a gate and there exists $c>0$ such that, for sufficiently large $\b$,
\be{gate} \P_{\vuoto} \big(\t_{\cP}>\t_{\pieno}\big)\leq e^{-\b c}.\ee \et

In words, Theorem \ref{t1'} says that the set $\cP$ is a gate for the nucleation; all paths from the metastable state $\vuoto$ to the stable state $\pieno$ go through this set with probability close to 1 as $\b\to\infty$.
Note that in this theorem we do not establish the minimality of
the gate $\cP$ (see (\ref{defP})), which would involve a much
more detailed analysis.

We refer to subsection \ref{timegate} for the proof of Theorem \ref{t1'}.

 \bt{t2}
Let $U_1, U_2, \D$ be  such that $U_2/(U_1+U_2-\D)$ is not integer
and (\ref{regime}) holds. Let $\L$ be a box with side $L+2$ and
let $\cR^{\le(l_1,l_2)}$ (resp. $\cR^{\ge(l_1,l_2)}$) be the set of
configurations whose single contour is a rectangle contained in
(resp. containing) a rectangle with sides $l_1$ and $l_2$. For $L$
sufficiently large,

\be{shrgrw}
\ba{lll}\hbox{ if }\quad
\h\in \cR^{\le(2l_2^*-3,l_2^*)} \hbox{ or } \h\in \cR^{\le(2l_2^*-1,l_2^*-1)}&\Longrightarrow&
\lim\limits_{\b\to\infty}\P_\h (\t_\vuoto < \t_\pieno) = 1,\\
\hbox{ if }\quad \h\in \cR^{\ge(2l_2^*-2,l_2^*)}&\Longrightarrow&
\lim\limits_{\b\to\infty}\P_\h (\t_\pieno < \t_\vuoto) = 1.
\ea
\ee \et

In other words, $2l^*_2-2$ and $l^*_2$ are the critical sizes, i.e., subcritical  rectangles shrink to $\vuoto$, supercritical rectangles grow to $\pieno$. We refer to subsection \ref{dimt2} for the proof of Theorem \ref{t2}.

\br{remwulff}
Theorem \ref{t2} implies that the Wulff shape is supercritical, indeed its circumscribed rectangle is $\cR(l_1^*,l_2^*)$ (see (\ref{defbarl})), that is in $\cR^{\geq(2l_2^*-2,l_2^*)}$. For example, consider $U_1=10U_2$, the circumscribed rectangle of a Wulff-shaped configuration is $\cR(10l_2^*,l_2^*)$. We consider a rectangle strictly smaller $\h=\cR(10l_2^*-1,l_2^*)\in\cR^{\geq(2l_2^*-2,l_2^*)}$, then by Theorem \ref{t2} it follows that $\h$ is supercritical (has tendency to grow) and not subcritical (has not tendency to shrink).
\er

\subsection{Comparison with weak anisotropy and dynamical heuristic discussion}
\label{euristicaa}
In this subsection we provide a detailed comparison between the strongly and the weakly anisotropic case. As we have already said, the behavior in the two regimes is very different and now we elaborate on why. For weak anisotropy we need the following additional definitions. Let

\be{defbarl} l_1^*:=\left\lceil{U_1 \over U_1+U_2
	-\D}\right\rceil, \qquad
\quad \bar l:=\Bigg\lceil{U_1-U_2\over U_1+U_2-\D}\Bigg\rceil.
\ee

For any $s> \bar l+2 $, if $s$ has the same parity as
$\bar l$ i.e.,
$[s-\bar l]_2=[0]_2$, then we define the set of {\it $0$-standard
	rectangles} as
$\cR^{0-st}(s):=\cR(\ell_1(s),\ell_2(s))$
with side lengths 

\be{ls0} \ell_1(s):= {s+\bar l\over 2} \qquad \ell_2(s):= {s-\bar
	l\over 2},\qquad \hbox{
	for } [s-\bar l]_2=[0]_2.
\ee
If $s$ has the same parity as $\bar l-1$ i.e.,  $[s-\bar
l]_2=[1]_2$, we define the set of {\it $1$-standard rectangles} to be
$\cR^{1-st}(s):=\cR(\ell_1(s),\ell_2(s))$
with side lengths

\be{ls1}
\ell_1(s):={s+\bar l-1\over 2}, \qquad \ell_2(s):={s-\bar l +1\over 2} \qquad
\hbox{ for } [s-\bar l]_2=[1]_2.
\ee

\noindent
For this value of $s$ we define the set of {\it quasi-standard}
rectangles as $\cR^{q-st}(s):=\cR(\ell_1(s)+1,\ell_2(s)-1)$. Finally, we set

\be{standardgen} \cR^{st}(s):=\left\{\ba{ll}
\cR^{0-st}(s)  &\mbox{if } [s-\bar l]_2= [0]_2\\
\cR^{1-st}(s)  &\mbox{if } [s-\bar l]_2= [1]_2. \ea\right. \ee

\bigskip
What happens for weak anisotropy is that, after an initial stage of the nucleation which consists of a growth along domino shapes (independently on the parameters of the interaction), the nucleation pattern consists of a growing sequence of standard rectangles up to configurations that have horizontal length equal to the side of the box (see Figure \ref{fig:strongweak} on the right hand-side). The critical droplet belongs to this sequence. This is a crucial difference with the strongly anisotropic case, for which the nucleation pattern follows the domino shape up to the critical droplet and then increases only adding columns up to the configurations with horizontal length equal to the side of the box (see Figure \ref{fig:lat1} in the middle), without involving the standard shape. In both cases, when the horizontal side is the same as the side of the box, the nucleation pattern grows in the vertical direction via the mechanisms ``change one column in row" and ``add column" (see Figure \ref{fig:lat1} and \ref{fig:lat0} on the right hand-side).

In \cite{NOS} and in the present paper we use the strategies suggested in \cite [Section 4.2]{MNOS} points I) and II). The point II) is the so-called {\it recurrence property} and is similar for the weak and strong anisotropic case (see \cite[Section 3.5]{NOS} and Proposition \ref{V^*} proved in Section \ref{recurrence}). Concerning point I), in both cases it is difficult to find a detailed description of the energy landscape, so the authors follow a general criterion to find a set $\cB$ satisfying properties (a) and (b). One of the ideas that were used to carry out this preliminary analysis in some cases consists of finding a suitable {\it foliation} of the state space $\cX$ into manifolds according to a certain parameter (for instance  the value of the semiperimeter $s$ for our model). In \cite{NOS} the authors introduced the foliation $\cV_s = \{\h\in\cX\colon\; p_1(\h)+p_2(\h)=s\}$ and they characterized the main property of standard rectangles $\cR^{st}(s)$: configurations that minimize the energy in $\cV_s$ for $s$ fixed. For the detailed result we refer to \cite[Proposition 8]{NOS}.

For the strongly anisotropic case we can not use this foliation and this result, because standard rectangles remain those minimizing the energy in the $s$-manifolds, but the nucleation pattern does not involve the standard shape. Indeed the energy needed to change a column into a row (see (\ref{sellestaccotrenini}) for the explicit formula) becomes much bigger in the case of strong anisotropy. To establish properties (a) and (b) see Corollary \ref{omegauppbou}, Theorem \ref{3.1} and Corollary \ref{phi}. Thus, summarizing, the idea behind the definition of the set $\cB$ in the weak and strong regime is similar, even though the strategies used in the proofs are quite different. Indeed without the tool of foliations and the identification of configurations with minimal energy on them, we need to carefully subdivide the proof in different cases.

In \cite{NOS} the authors defined

\be{defPweak}
\ba{ll}
\tilde\cP:= \{\h&:\, n(\h)= 1,\, v(\h)=\ell_2(l_1^*+l_2^*-1)-1,\,
\h_{cl} \hbox{ is connected, monotone,}\\
&\hbox{ with circumscribed rectangle in }
\cR(\ell_1(l_1^*+l_2^*-1)+1,\ell_2(l_1^*+l_2^*-1))\},
\ea
\ee

and proved that this set is a gate for the transition between $\vuoto$ and $\pieno$ (see Figure \ref{fig:Pweak} for an example of configurations in $\tilde\cP$ for weak anisotropy).

\setlength{\unitlength}{1.1pt}
\begin{figure}
	\begin{picture}(400,90)(0,30)
	\thinlines
	\qbezier[51](20,0)(100,0)(180,0)
	\qbezier[51](20,90)(100,90)(180,90)
	\qbezier[51](20,0)(20,45)(20,90)
	\qbezier[51](180,0)(180,45)(180,90)
	\thinlines
	\put(20,0){\line(1,0){150}}
	\put(20,90){\line(1,0){150}}
	\put(20,0){\line(0,1){90}}
	\put(170,0){\line(0,1){60}}
	\put(170,70){\line(0,1){20}}
	\put(170,60){\line(1,0){10}}
	\put(180,60){\line(0,1){10}}
	\put(170,70){\line(1,0){10}}
	\put(140,100){\line(1,0){10}}
	\put(150,100){\line(0,1){10}}
	\put(140,100){\line(0,1){10}}
	\put(140,110){\line(1,0){10}}
	\thinlines  \put(190,55){$l_2^*$}
	\put(100,-15){$l_1^*$}
	\thinlines
	\qbezier[51](220,0)(300,0)(380,0)
	\qbezier[51](220,90)(300,90)(380,90)
	\qbezier[51](220,0)(220,45)(220,90)
	\qbezier[51](380,0)(380,45)(380,90)
	\thinlines
	\put(230,10){\line(1,0){10}}
	\put(240,10){\line(0,-1){10}}
	\put(240,0){\line(1,0){130}}
	\put(370,0){\line(0,1){20}}
	\put(370,20){\line(1,0){10}}
	\put(380,20){\line(0,1){60}}
	\put(370,80){\line(1,0){10}}
	\put(370,80){\line(0,1){10}}
	\put(230,10){\line(0,1){10}}
	\put(230,20){\line(-1,0){10}}
	\put(220,20){\line(0,1){60}}
	\put(220,80){\line(1,0){20}}
	\put(240,80){\line(0,1){10}}
	\put(240,90){\line(1,0){130}}
	\put(300,100){\line(1,0){10}}
	\put(310,100){\line(0,1){10}}
	\put(300,100){\line(0,1){10}}
	\put(300,110){\line(1,0){10}}
	\put(390,55){$l_2^*$}
	\put(300 ,-15){$l_1^*$}
	\end{picture}
	\vskip 2. cm
	\caption{Configurations in $\tilde\cP$ for weak anisotropy}
	\label{fig:Pweak}
\end{figure}
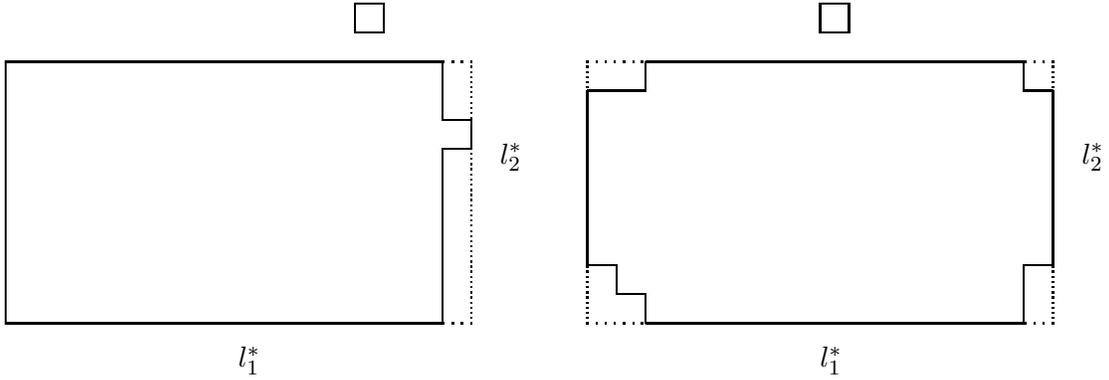

Note that, except for the length of the vertical and horizontal sides, $\tilde\cP$ coincides with $\cP_2$. This feature is due to the definition of the set $\cB$. The other set of saddles $\cP_1$ for the strong anisotropic case is a peculiar feature of this model, because in this case we have two minimal saddles at the same height from an energetical point of view. 

Roughly speaking, the difference between the two parameter regimes depends on the fact that in order to go up in energy by a factor $U_1$ for weak anistropy it is sufficient to go up by $2U_2$, since $U_1<2U_2$. For strong anisotropy this is not possible, since $U_1>2U_2$ and we have not an upper bound of $U_1$ in terms of $U_2$, so that more effort is needed: the key strategy of our proof is to analyze in a very detailed way all the possible exit moves from the set $\cB$, see (\ref{defB}). 

For the anisotropic case in \cite[Sections 2.1,2.2]{NOS}, the static and dynamic heuristics are discussed as was done in \cite[Sections 1.3,1.4]{HNOS}. Here we will summarize the key ideas of the dynamical heuristic description for the strong anisotropy.

\medskip
\noindent 
{\bf Key transitions}

\noindent
We start with a coarse-grained description: we will
restrict ourselves to the determination of the sequence of rectangles visited
by typical trajectories. By the continuity properties of the dynamics it is reasonable to
expect that only transitions between neighboring rectangles have
to be taken into consideration. More precisely, starting from a configuration  $\h\in \cR (l_1,l_2)$, with $l_1,l_2\geq 2$, the possible successive rectangles in the tube have to belong to one of the following classes: $\cR(l_1+1 , l_2)$, $\cR (l_1 , l_2+1)$, $\cR (l_1-1 , l_2)$, $\cR(l_1 , l_2-1)$, $\cR (l_1-1 , l_2+1)$ and $\cR (l_1+1 , l_2-1)$. So we shall consider the following transitions:
\bi

\item
from  $\cR (l_1, l_2)$ to $\cR (l_1, l_2+1)$,   corresponding to
vertical growth, that will be denominated  {\it add row } and
symbolically denoted by the arrow $\uparrow$ pointing north
direction;

\item
from  $\cR (l_1, l_2)$ to $\cR (l_1+1, l_2)$,   corresponding to
horizontal growth, that will be denominated  {\it add column } and
denoted by the arrow $\rightarrow$ pointing east;

\item
from  $\cR (l_1, l_2)$ to $\cR (l_1, l_2-1)$,   corresponding to
vertical contraction, that will be denominated  {\it remove row}
and denoted by the arrow $\downarrow$ pointing south;

\item
from  $\cR (l_1, l_2)$ to $\cR (l_1-1, l_2)$,   corresponding to
horizontal  contraction, that will be denominated  {\it remove
	column } and denoted by the arrow $\leftarrow$ pointing west;

\item
from  $\cR (l_1, l_2)$ to $\cR (l_1-1, l_2+1)$,   corresponding to
a readjustment of the  edges, making higher and narrower the
rectangle by removing a column and simultaneously adding a row. It
will be denominated  {\it column to row } and denoted by the arrow
$\nwarrow$ pointing northwest;

\item
from  $\cR (l_1, l_2)$ to $\cR (l_1+1, l_2-1)$,  corresponding to
a readjustment opposite to the previous one. It will be
denominated  {\it  row to column } and denoted by the arrow
$\searrow$ pointing southeast.
\ei

\medskip
The transition  from $\cR (l_1, l_2)$ to $\cR (l_1-1, l_2-1)$ and
$\cR (l_1+1, l_2+1)$ are not considered as elementary since, as
it can be easily seen, a suitable combination of two of the above
transitions takes place with larger probability.

\begin{figure}
	\begin{center}
		\includegraphics[width=14cm]{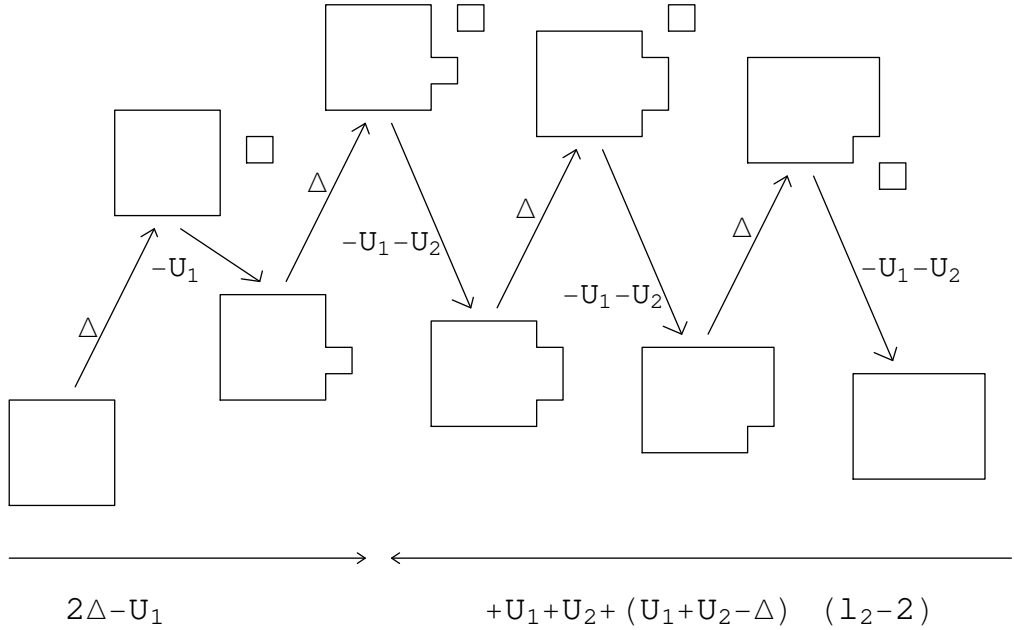}\\
		\caption{The procedure to grow a column}
		\label{fig:growc}
	\end{center}
\end{figure}

\begin{figure}[h]
	\begin{center}
		\includegraphics[width=14cm]{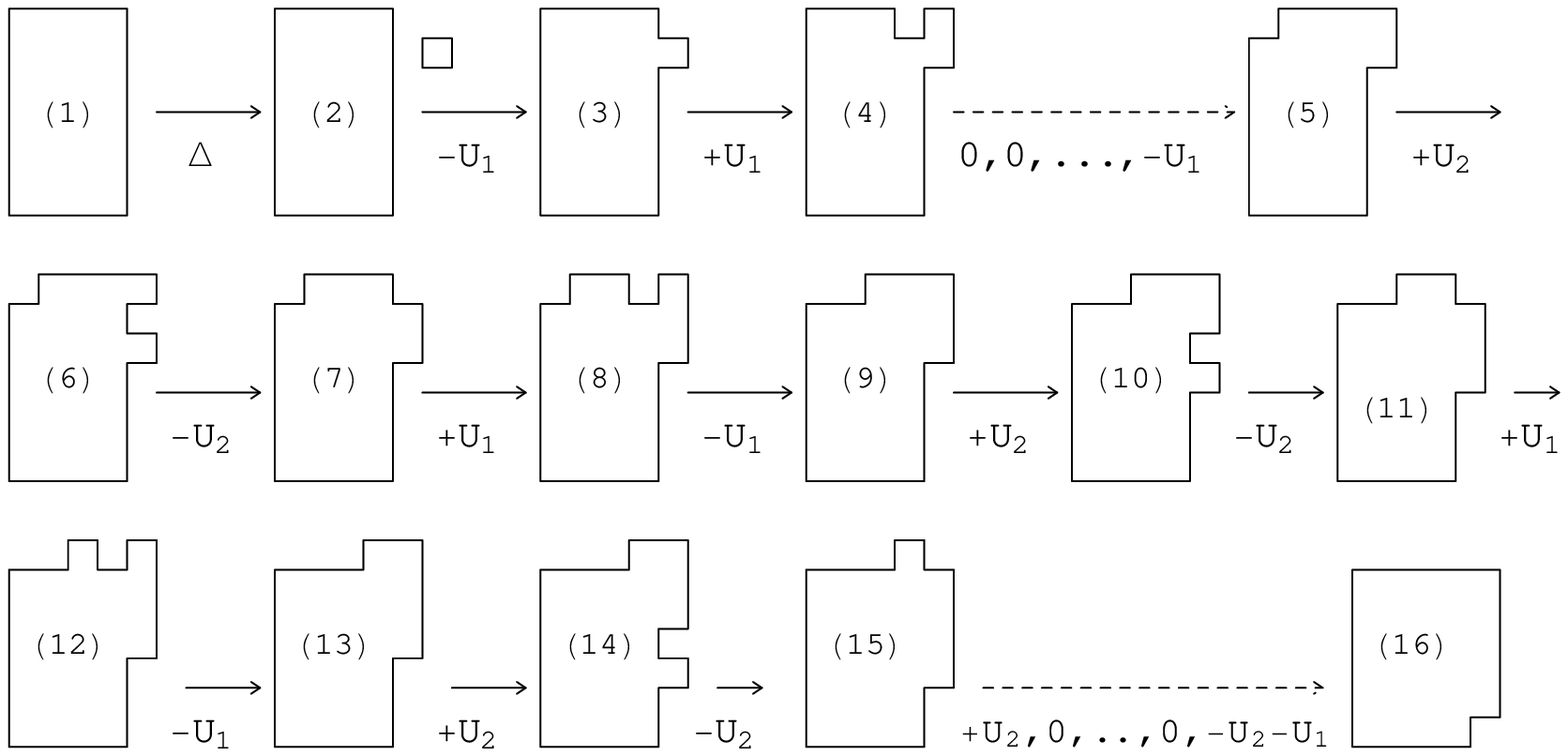}\\
		\caption{A path in $\Omega^{(5)}$}
		\label{fig:O5}
	\end{center}
\end{figure}

\begin{figure}[h]
	\begin{center}
		\includegraphics[width=14cm]{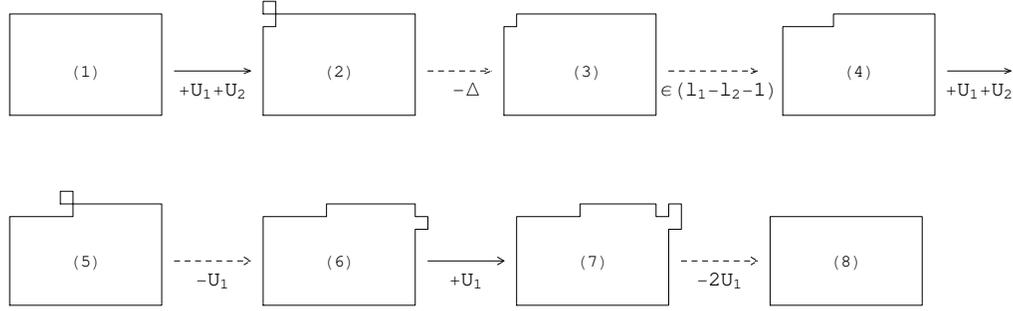}\\
		\caption{A path in $\bar{\Omega}^{(5)}$}
		\label{fig:barO5}
	\end{center}
\end{figure}

At first sight the optimal interpolation paths realizing the above
transitions between contiguous rectangles are the ones depicted in
Figures \ref{fig:growc}, \ref{fig:O5} and \ref{fig:barO5}. Let us call $ \Omega^{(1)} $ the set of
paths as the one depicted in Figure \ref{fig:growc}. They are the natural
candidates to realize, in an optimal way, the transition
$\rightarrow$. For the transition $\uparrow$ we have an analogous
set of  paths that we call $\Omega^{(2)}$.

Let us call $B$ the time-reversal operator acting on finite
paths; we have for $\o = \o_1, \ldots ,\o_T$ \be{timereversal}
B\o= \o' \hbox{  with  } \o_i' = \o_{T+1-i} \quad i=1, \ldots,
T.
\ee

\noindent 
For the transition $\downarrow$ we choose the set of
paths $\Omega^{(3)}$ obtained  by time-reversal from the paths,
analogous to the ones in $ \Omega^{(1)} $, that realize the
transition $\cR(l_1-1, l_2)$ to $\cR(l_1, l_2)$.

Similarly, for the transition $\leftarrow$ we use the set of
paths $ \Omega^{(4)} $ obtained by time-reversal from the paths,
analogous to the ones in $ \Omega^{(2)} $, that realize the
transition $\cR(l_1, l_2-1)$ to $\cR(l_1, l_2)$.

The set of paths that we consider as the optimal interpolation for the transition from  $\cR(l_1, l_2)$ to $\cR(l_1-1, l_2+1)$ in the two cases $l_1<l_2$ and $l_1\geq l_2$ are called $\Omega^{(5)}$ and $\bar \Omega^{(5)}$ respectively. A path in $\Omega^{(5)}$ is represented in figure \ref{fig:O5} where each arrow corresponds to a move and the quantities under the arrows represent the corresponding energy barriers $\Delta H$. Dotted arrows indicate sequences of moves. The maximal energy along the path is reached in the configuration (2). A path in $\bar \Omega^{(5)}$ is represented in figure \ref{fig:barO5} where to simplify we indicate under the dotted arrows the sum of the corresponding $\Delta H$. Along this path the maximal energy is reached in configuration (5). In a similar way we define the optimal interpolation paths $\Omega^{(6)}$ and $\bar\Omega^{(6)}$ for the transition from $\cR(l_1, l_2)$ to $\cR(l_1+1, l_2-1)$. We call {\it canonical} the paths in the above sets.

\medskip
Given $(l_1, l_2)$, in order to determine the most probable transition
between $\cR (l_1,l_2)$  and one of the previous six contiguous
rectangles, we will use the criterion of the smallest energy
barrier, defined as the difference between the communication
height and $H(\cR(l_1, l_2))$. We  call {\it energy barrier  from
	$\h$ to $\h'$ along the path $\o=
	(\o_1=\h, \ldots , \o_n=\h') $} the difference between the maximal
height reached along this path and $H(\h)$.
We compute the energy barriers  along the canonical paths and we use
them to estimate the true energy
barriers. We denote by  $\D H \mbox{(add  row)}$ the energy
barrier along  the paths in $\Omega^{(1)}$; similarly for the
other transitions.

From Figures \ref{fig:growc}, \ref{fig:O5} and \ref{fig:barO5} via  easy computations, we get:
\be{sellestaccotrenini}
\left.\ba{llllll}
\D H\mbox{(add  row)} &=  2\D-U_2 &\qquad\\
\D H\mbox{(add column)} &=  2\D-U_1&\qquad\\
\D H\mbox{(remove row)}
&=\varepsilon(l_1 -2) +U_1+U_2&\qquad\\
\D H\mbox{(remove column)}
&=\varepsilon(l_2 -2) +U_1+U_2&\qquad\\
\D H\mbox{(row} \hbox{ to }\mbox{column)  }
&=\D\;  &\mbox{if } l_1<l_2 \\
\D H\mbox{(row} \hbox{ to }\mbox{column)  }
&=U_1+U_2+\e(l_1-l_2)\; &\mbox{if } \; l_1\geq l_2 \\
\D H\mbox{(column} \hbox{ to }\mbox{row)  }
&=\D-U_2+U_1\; &\mbox{if }  \; l_1> l_2\\
\D H\mbox{(column} \hbox{ to }\mbox{row)  }
&=\D-U_2+U_1+\varepsilon(l_2-l_1+1)\;
&\mbox{if }  \; l_1\leq l_2\\
\ea \right.\ee

These  estimated energy barriers  are, of course, larger than or
equal to the true ones; the equality does not hold in general,
since the above canonical paths sometimes happen to  be
non-optimal. For example a deeper analysis leads to the conclusion
that to add a row,  instead of using  a path in $\Omega^{(1)}$, it
is more convenient to  compose $\Omega^{(2)}$ and $\Omega^{(5)}$,
resp. $\bar \Omega^{(5)}$, when $l_1<l_2$, resp. $l_1\ge l_2$.

Let us now make a comparison between the estimated energy barriers
appearing in  equation (\ref{sellestaccotrenini}). For $l_1\leq l_2$, we can easily check that $\D H\mbox{(row } \hbox{ to }\mbox{column) }\leq U_1+U_2 =\D+\e$
is the smallest estimated energy barrier. So in the sequel we will consider only the case $l_1>l_2$. For $l_1>l_2$, since in the strongly anisotropic case $U_1>2U_2$ and 
	$$2\D-U_1< 2 \D -U_2, \;  2\D-U_1<\D-U_2+U_1  \hbox{ and } $$
	\be{comp5'}
	U_1+U_2+\e (l_2-2)<U_1+U_2+\e (l_1-2), \ee
	by (\ref{sellestaccotrenini}), we deduce that  we have only to
	compare $\D H\mbox{(remove  column)}$, 
	
	$\D H\mbox{(add  column)}$ and $\D
	H\mbox{(row to column)}$.
	We get
	
	\be{comp1}\D H\mbox{(remove  column)}<\D H\mbox{(add  column)}\iff
	l_2<  l_2^*,\ee

	\be{comp3} \D H\mbox{(row}\hbox{ to }\mbox{column)}<\D H\mbox{(add
		column)}
	\iff  l_1< l_2+  l_2^*-2 .\ee
	
	\be{comp8}\D H\mbox{(remove  column)}\leq \D H\mbox{(row}\hbox{ to
	}\mbox{column)}
	\iff 2l_2-2 \leq l_1\ee
	
	\begin{figure}[h!]
		\begin{center}
			\includegraphics[width=14cm]{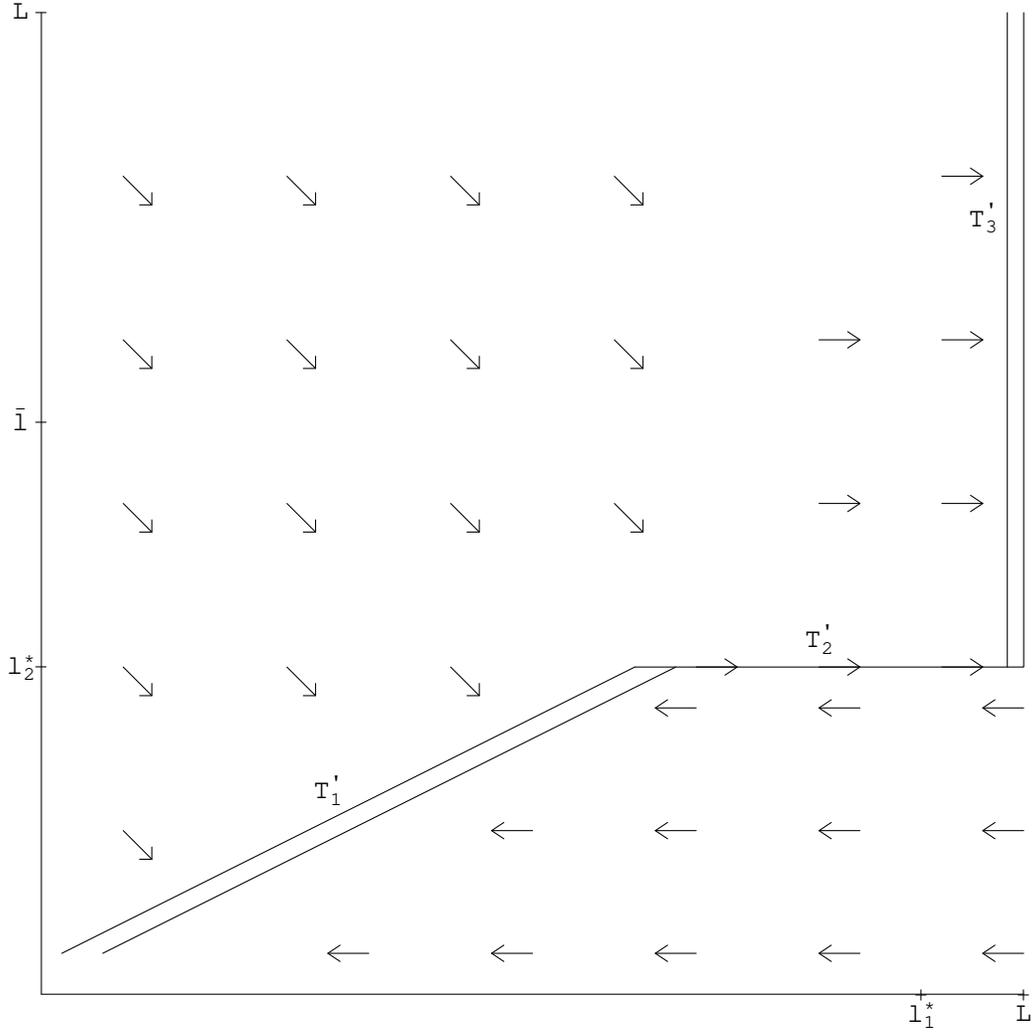}\\
			\caption{Strong anisotropy: minimal transitions and tube of typical trajectories. We choose to indicate a few arrows to have a comprehensible picture.}
			\label{fig:strongtube}
		\end{center}
	\end{figure}

Summarizing we have that:

$\bullet$ in the set $A'=\{l_2\leq l_2^*-1, \; l_1> 2l_2-2\}$ the
minimal estimated energy barrier is $\D H\mbox{(remove}$ $\mbox{
	column)}$;

$\bullet$ in the set $B'=\{l_1< l_2+ l_2^*-2,\; l_1< 2l_2-2\}$
the minimal estimated energy barrier is
$\D H\mbox{(row} \hbox{ to }\mbox{column)  }$;

$\bullet$ in the set $C'=\{l_2 \geq  l_2^*,\; l_1\geq l_2+ l_2^*
-2\}$ the minimal estimated energy barrier is $\D H\mbox{(add
	column)}$.

$\bullet$ in the set $D'=\{l_2\leq l^*_2-1,\; l_1=2 l_2 -2\}$ we
have degeneracy of the minimal estimated estimated energy barrier:
$\D H\mbox{(remove  column)}=$ $\D H(\mbox{row} \hbox{ to
}\mbox{column})$

\medskip
Note  that  $B'=\{ l_2\leq l_2^*-1,\;  l_1< 2l_2-2 \}\, \cup \,\{
l_2\geq l_2^* ,\; l_1 < l_2+ l_2^* -2\}$, so that $A'\cup
B'\cup C'\cup D'= \{l_1>l_2\}$.

In Figure \ref{fig:strongtube} we represent $\cR(l_1, l_2)$ as points in $\Z^2$ of coordinates $l_1, l_2$ (representing, respectively, the horizontal and
vertical edges). Emerging from any representative point, we draw
the arrows  corresponding to  transitions with minimal $\D H$
between $\cR(l_1, l_2)$ and contiguous rectangles.

\medskip
\noindent
{\bf Strongly anisotropic case}

\noindent
In the strongly anisotropic case, from Figure \ref{fig:strongtube}, it is evident that in the plane $(l_1, l_2)$
there is a connected region $T'$, which is attractive in the sense that if we follow the oriented paths given by the sequences of arrows emerging from every points outside $T'$ we end up inside $T'$. It consists of three
parts $T'_1=\{(l_1, l_2): l_2< l_2^* \hbox{ and }  2 l_2-3 \leq
l_1 \leq 2l_2 -1 \}\cup \cR(2l^*_2-3,l^*2)$  containing domino
shape rectangles,  $T'_2=\{(l_1, l_2): l_2 = l^*_2\hbox{ and }
l_2+ l_2^* -2 \leq l_1 <L \}$,   and  $T'_3=\{(l_1, l_2):  l^*_2
\leq l_2 \hbox{ and } L-1 \leq l_1 \leq L \}$.

Let us now consider the arrows inside the region $T'$. From each $\h\in{T_1'}$, with $l_1=2l_2-2$, as a consequence of the deneracy $\D H\mbox{(remove  column)}=\D H\mbox{(row to  column)}$, we have two arrows, one pointing to $\h'\in\cR(l_1-1,l_2)$ and the other pointing to $\h''\in\cR(l_1+1,l_2-1)$. Subsequently, starting from $\eta'$ the minimal estimated $\D H$ is unique and it corresponds to an arrow pointing to $\cR(l_1,l_2-1)$; analogously starting from $\h''$ the minimal $\D H$ is unique and it corresponds to an arrow also pointing to $\cR(l_1,l_2-1)$ (see Figure \ref{fig:lat1} on the left hand-side). 

\setlength{\unitlength}{1.15 pt}
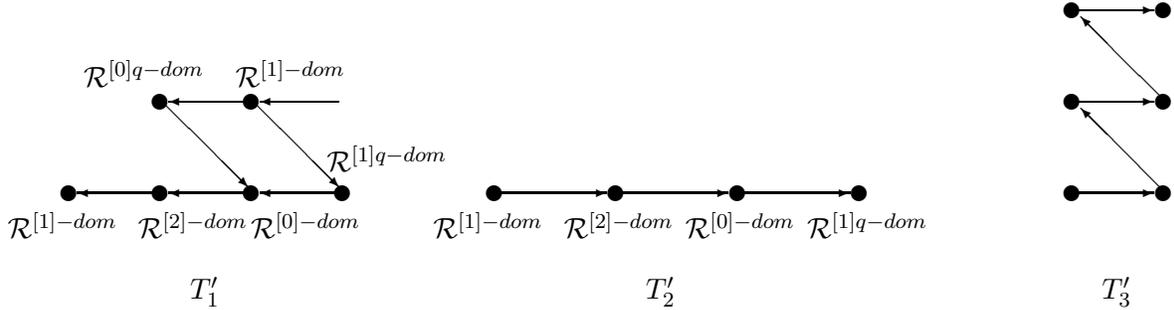
\begin{figure}
	\begin{picture}(400,80)(-20,20)
		\put(0,10){\circle*{5}}\put(-20,-5){${\cal R}^{[1]-dom}$}
		\put(30,10){\circle*{5}}\put(23,-5){${\cal R}^{[2]-dom}$}
		\put(60,10){\circle*{5}} \put(60,-5){${\cal R}^{[0]-dom}$}
		\put(90,10){\circle*{5}} \put(85,17){${\cal R}^{[1]q-dom}$}
		\put(30,40){\circle*{5}}\put(5,45){${\cal R}^{[0]q-dom}$}
		\put(60,40){\circle*{5}}\put(55,45){${\cal R}^{[1]-dom}$}
		\put(89,40){\vector(-1,0){26}}\put(59,40){\vector(-1,0){26}}
		\put(89,10){\vector(-1,0){26}}
		\put(59,10){\vector(-1,0){26}} \put(29,10){\vector(-1,0){26}}
		\put(31,40){\vector(1,-1){28.2}} \put(61,40){\vector(1,-1){28.2}}
		\put(40,-25){$ T_1'$}
		
		\put(140,10){\circle*{5}}\put(120,-5){${\cal R}^{[1]-dom}$}
		\put(180,10){\circle*{5}}\put(163,-5){${\cal R}^{[2]-dom}$}
		\put(220,10){\circle*{5}}\put(203,-5){${\cal R}^{[0]-dom}$}
		\put(260,10){\circle*{5}}\put(243,-5){${\cal R}^{[1]q-dom}$}
		\put(139,10){\vector(1,0){38}} 
		\put(179,10){\vector(1,0){38}} 
		\put(219,10){\vector(1,0){38}} 

		\put(190,-25){$ T_2'$}

		
		\put(330,10){\circle*{5}}
		\put(360,10){\circle*{5}}
		\put(330,40){\circle*{5}}
		\put(360,40){\circle*{5}}
		\put(331,10){\vector(1,0){26}} 
		\put(331,40){\vector(1,0){26}} 
		\put(361,10){\vector(-1,1){28}} 
		\put(330,70){\circle*{5}}
		\put(360,70){\circle*{5}}
		\put(331,70){\vector(1,0){26}} 
		\put(361,40){\vector(-1,1){28}} 
		
		\put(340,-25){$ T_3'$}

	\end{picture}\vskip 1.8 cm
	\caption{Minimal transition inside $T_1'$, $T_2'$ and $T_3'$.}
	\label{fig:lat1}
\end{figure}

For every configuration in $T'_2$, the
minimal estimated energy barrier is $\D H\mbox{(add  column)}$,
which implies
that the rectangles in $T'_2$ will grow  in the horizontal direction until they  become
a complete horizontal strip with  length $L$ (see Figure \ref{fig:lat1} in the middle). In $T'_3$ the
minimal estimated energy barrier is $\D H\mbox{(add row)}$, which
implies that every horizontal strip with $l_1=L$  will
grow in the vertical direction until it covers the
whole box (see Figure \ref{fig:lat1} on the right hand-side). 

It is natural at this point to distinguish two parts in the set $T'$: the subcritical part $T'_{sub}$ corresponding to $T_1'$ and the supercritical part $T'_{sup}$, corresponding to the configurations in $T'_2$ and $T'_3$.

Let us now summarize our heuristic discussion in the strongly anisotropic case. We expect that every rectangle outside $T'$ is attracted by $T'$; the configurations in $T'_{sub}$ are subcritical in the sense that they tend to shrink along $T'$ following domino shapes; configurations in $T'_{sup}$ are supercritical in the sense that they tend to grow following domino shapes in $T'_2$ and a sequences of rectangles with bases $L-1$ or $L$ in $T'_3$. The nucleation pattern in the strongly anisotropic case contains a sequence of increasing domino shaped rectangles up to $\cR(2l_2^*,l_2^*)$; then a sequence of rectangles with $l_2=l_2^*$ and $l_1$ going up to $L$ (the size of the container); finally a sequences of horizontal strips whose width grows from $l_2^*$ to $L$. We note that the nucleation pattern in the strongly anisotropic case is very similar to the one that we would have for non-conservative Glauber dynamics for \emph{any} anisotropy.

This heuristic discussion provides a description of the tube of the typical nucleating path. Suppose first to consider the typical paths going from the maximal subcritical rectangle to $\vuoto$. From the discussion on the subcritical part, we have that the sequence of cycles follows the arrows as in Figure \ref{fig:lat1}. Looking at $T'_1$ we see that there are no loops there, so we can associate to each rectangular configuration $\h$ in $T'_1$ the maximal cycle containing $\h$ and not containing other rectangular configurations: by using the arrows of the figure we obtain, in this way, a coarse-grained cycle path corresponding to the first domino part of the tube. The \textcolor{red}{coarse-grained} of these cycle paths can be resolved by introducing a suitable interpolation between rectangular configurations corresponding to each arrow in the picture, obtaining, in this way, a family of true cycle paths $T'_{sub}$, describing the tube of typical paths going from the maximal subcritical rectangle to $\vuoto$.

A similar discussion can be applied to the study of the tube of typical paths going from the minimal supercritical rectangle to $\pieno$ obtaining in the same way the family of cycles $T'_{sup}$.

We expect that the first supercritical rectangle configuration is contained in $\cR(2l_2^*,l_2^*)$: we see in the next that it is $\cR(2l_2^*-2,l_2^*)$.

The tube of the typical nucleating paths describing the first excursion from $\vuoto$ to $\pieno$, can be obtained by applying general arguments based on reversibility and by providing a suitable interpolation between the maximal subcritical rectangle and the minimal supercritical one. More precisely, to obtain the typical tube from $\vuoto$ to $\pieno$ we apply the time reversal operator $B$ (see \ref{timereversal}), to the tube $T'_{sub}$ and we join it to $T'_{sup}$ by means of this interpolation. These interpolations between rectangular configurations can be obtained by using the reference path $\o^*$ described in subsection \ref{refpathsec}; $\o^*$ can be considered as a representative of a typical nucleation path.

Some aspects of the behavior that we have heuristically described are rigorously discussed in this work; in particular we determine a gate $\cP$ for the transition between $\vuoto$ and $\pieno$ (see Theorem \ref{t1'}), even though we say nothing about its minimality, and we give a sufficient condition to discriminate subcritical and supercritical domino rectangles (see Theorem \ref{t2}).

\setlength{\unitlength}{1.15 pt}
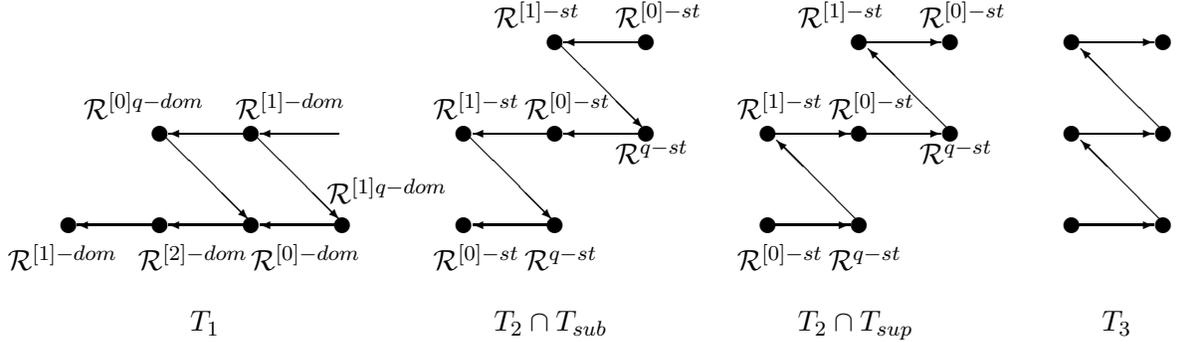
\begin{figure}
	\begin{picture}(400,80)(-20,20)
		\put(0,10){\circle*{5}}\put(-20,-5){${\cal R}^{[1]-dom}$}
		\put(30,10){\circle*{5}}\put(23,-5){${\cal R}^{[2]-dom}$}
		\put(60,10){\circle*{5}} \put(60,-5){${\cal R}^{[0]-dom}$}
		\put(90,10){\circle*{5}} \put(85,17){${\cal R}^{[1]q-dom}$}
		\put(30,40){\circle*{5}}\put(5,45){${\cal R}^{[0]q-dom}$}
		\put(60,40){\circle*{5}}\put(55,45){${\cal R}^{[1]-dom}$}
		\put(89,40){\vector(-1,0){26}}\put(59,40){\vector(-1,0){26}}
		\put(89,10){\vector(-1,0){26}}
		\put(59,10){\vector(-1,0){26}} \put(29,10){\vector(-1,0){26}}
		\put(31,40){\vector(1,-1){28.2}} \put(61,40){\vector(1,-1){28.2}}
		\put(40,-25){$ T_1$}
		
		\put(130,10){\circle*{5}}\put(120,-5){${\cal R}^{[0]-st}$}
		\put(160,10){\circle*{5}}\put(150,-5){${\cal R}^{q-st}$}
		\put(130,40){\circle*{5}}\put(120,45){${\cal R}^{[1]-st}$}
		\put(160,40){\circle*{5}}\put(150,45){${\cal R}^{[0]-st}$}
		\put(159,40){\vector(-1,0){26}} 
		\put(159,10){\vector(-1,0){26}} 
		\put(131,40){\vector(1,-1){28}} 
		\put(190,40){\circle*{5}}\put(180,30){${\cal R}^{q-st}$}
		\put(160,70){\circle*{5}}\put(140,75){${\cal R}^{[1]-st}$}
		\put(190,70){\circle*{5}}\put(180,75){${\cal R}^{[0]-st}$}
		\put(189,40){\vector(-1,0){26}} 
		\put(189,70){\vector(-1,0){26}} 
		\put(161,70){\vector(1,-1){28}} 
		
		\put(140,-25){$ T_2\cap T_{sub}$}
		
		
		\put(230,10){\circle*{5}}\put(220,-5){${\cal R}^{[0]-st}$}
		\put(260,10){\circle*{5}}\put(250,-5){${\cal R}^{q-st}$}
		\put(230,40){\circle*{5}}\put(220,45){${\cal R}^{[1]-st}$}
		\put(260,40){\circle*{5}}\put(250,45){${\cal R}^{[0]-st}$}
		\put(231,40){\vector(1,0){26}} 
		\put(231,10){\vector(1,0){26}} 
		\put(261,10){\vector(-1,+1){28}} 
		\put(290,70){\circle*{5}}\put(280,30){${\cal R}^{q-st}$}
		\put(260,70){\circle*{5}}\put(240,75){${\cal R}^{[1]-st}$}
		\put(290,40){\circle*{5}}\put(280,75){${\cal R}^{[0]-st}$}
		\put(261,40){\vector(1,0){26}} 
		\put(261,70){\vector(1,0){26}} 
		\put(291,40){\vector(-1,1){28}} 
		
		\put(240,-25){$ T_2\cap T_{sup}$}

		
		\put(330,10){\circle*{5}}
		\put(360,10){\circle*{5}}
		\put(330,40){\circle*{5}}
		\put(360,40){\circle*{5}}
		\put(331,10){\vector(1,0){26}} 
		\put(331,40){\vector(1,0){26}} 
		\put(361,10){\vector(-1,1){28}} 
		\put(330,70){\circle*{5}}
		\put(360,70){\circle*{5}}
		\put(331,70){\vector(1,0){26}} 
		\put(361,40){\vector(-1,1){28}} 
		
		\put(340,-25){$ T_3$}

	\end{picture}\vskip 1.8 cm
	\caption{Minimal transition inside $T_1$, $T_2$ with $l_2<l_2^*$, $T_2$ with $l_2\geq l_2^*$ and $T_3$.}
	\label{fig:lat0}
\end{figure}

\medskip
\noindent
{\bf Weakly anisotropic case}

\noindent
In the weakly anisotropic case the main difference is that, since $U_1<2U_2$, $\D H(\hbox{column to row})$ is smaller than in the strong anisotropic case, thus it plays an important role. In the region above $A'$ this saddle is the minimal up to the standard shape where $\D H(\hbox{row to column})=\D H(\hbox{column to row})$. We refer to \cite[Section 2.2]{NOS} for a more detailed discussion of the weak anisotropic case. The region $T$ consists of three parts: $T_1=\{(l_1, l_2):
l_2\leq \bar l \hbox{ and }  2 l_2-3 \leq l_1 \leq 2l_2 -1 \}$
containing domino shape rectangles and $T_2=\{(l_1, l_2): l_2 >
\bar l \hbox{ and } l_2+\bar l -1\leq l_1\leq l_2+\bar l    \}$
containing standard rectangles (see (\ref{ls0}) and (\ref{ls1})) and
$T_3=\{(l_1, l_2): l_1 =L \hbox{ and } l_2\geq L-\bar l \}$.

The properties of $T_1$ can be discussed in analogy with $T_1'$ (see Figure \ref{fig:lat0} on the left hand-side).

In $T_2$, for each value of the semi-perimeter $s$, there are
pairs of configurations $(\h, \h')$ such that the minimal  among
the estimated energy barriers starting from $\h$ corresponds to
the transition from $\h$ to $\h'$ and conversely the minimal
estimated energy barrier from $\h'$ corresponds to the transition
from $\h'$ to $\h$. So inside $T_2$ there are pairs of arrows
forming
two-states loops that we represent as $\searrow$ \hskip-15pt$\nwarrow$.
This suggests that in $T_2$  a more detailed study  is necessary,
based  on the analysis of suitable  cycles  containing  the above
described loops. These cycles represent a sort of generalized
basin of attraction of the standard rectangles contained in the
loops: they are the maximal cycles containing a unique standard
rectangle. These cycles contain, among others,  rectangular configurations and in each cycle all the rectangular configurations  have  the same semiperimeter $s$, i.e., belong to the same manifold $\cV_s$.

We draw in our picture the arrows between rectangular configurations
corresponding to these most probable exits. It turns out that
these arrows are
horizontal
pointing east if $l_2\ge l^*_2$ and pointing west if $l_2<l^*_2$ (see Figure \ref{fig:lat0} in the middle). In both cases these horizontal arrows point to
configurations which are again in the set $T$, so that we can
iterate the argument to analyze all the arrows in $T$. Thus we
associate to the loops $\searrow$ \hskip-15pt$\nwarrow$
in the
picture cycles containing rectangles in $\cV_s$ and  transitions
given by the horizontal arrows. In $T_3$ we can argue like in
$T_2$ (see Fig. \ref{fig:lat0} on the right hand-side).

It is natural at this point to distinguish two parts in the set
$T$: the subcritical part $T_{sub}$ corresponding to $T_1$ plus
the part of $T_2$ with horizontal arrows pointing  west, i.e.,
with $l_2<l^*_2$ and the supercritical part $T_{sup}$,
corresponding to the configurations  in $T_2$ with horizontal
arrows pointing east, i.e., with $l_2\ge l^*_2$ and $T_3$.

Let us  now summarize our heuristic discussion in the weakly
anisotropic case. We expect that every rectangle outside $T$ is
attracted by $T$; the configurations in $T_{sub}$ tend to shrink along $T$ following either the standard or the domino shape, depending on $l_2$; configurations  in $T_{sup}$ tend to grow following standard shapes in $T_2$ and a sequences of rectangles with bases $L-1$ or $L$ in $T_3$.

\section{Model-dependent results}
\label{depres}

\subsection{Extensive model-dependent definitions and notations}
\label{extmoddef}
In this subsection we extend the model-dependent definitions given in Subsection \ref{moddef} following \cite{NOS}, that will be useful for characterizing configurations from a geometrical point of view.

	\medskip\noindent
	{\bf Clusters and projections.}
	
	Next we introduce  a geometric description of the configurations
	in terms of contours.
	
	\begin{itemize}
	\item[$\bullet$]
	Given a configuration $\h \in\cX$, consider the set $C(
	\h_{cl}) \subset \R^2$ defined as the union of the $1\times 1$
	closed squares centered at the occupied sites of $ \h_{cl}$ in $
	\L_0$. The maximal connected components $C_1,\dots, C_m$
	($m\in\N$) of $C( \h_{cl})$ are called {\it clusters} of $ \h$.
	There is a one-to-one
	correspondence between configurations $\h_{cl}\subset\L_0$ and sets $C(\h_{cl})$.
	A configuration $\h\in\cX$ is characterized by a set $C(\h_{cl})$, depending
	only on $\h\cap\L_0$, plus possibly a set of free particles in $\partial^-\L$
	and in $\L_0$. We are actually identifying three different
	objects: $\h\in\cX$, its support $\hbox{supp}(\h)\subset\L$, and the pair $(C(\h_{cl}),
	\h_{fp})$; we  write $x\in\h$ to indicate that $\h$ has
	a particle at $x \in \L$.
	
	\item[$\bullet$]
	For $\h\in\cX$, let $|\eta|$ be the number of particles in $\eta$,
	$\g(\h)$ the Euclidean boundary of $C(\h_{cl})$, $\g(\h)=\partial
	C(\h_{cl})$;   we denote by $g_1(\eta)$
	(resp. $g_2(\eta)$) one half of the horizontal (resp. vertical) length
	of $\g(\h)$, i.e., one half of the number of horizontal (vertical) broken bonds in $\h_{cl}$. Then
	the
	energy associated with $\h$ is given by
	
	\be{Hcont}
	H(\h) = - (U_1+U_2-\D)|\eta_{cl}| + {U_1}
	g_2(\eta)+
	{U_2} g_1(\eta)
	+\D |\h_{fp}|.
	\ee
	The maximal connected components of $\partial C(\h_{cl})$ are
	called {\it contours} of $\h$.
	
	\item[$\bullet$]
	Let $p_1(\h)$
	and $p_2(\h)$ be the total lengths of  horizontal and vertical
	projections of $ C(\h_{cl})$ respectively. More precisely let
	$r_{j,1}=\{x \in \Z^2:(x)_1=j\}$ be the $j$-th column and
	$r_{j,2}=\{x \in \Z^2:(x)_2=j\}$ be the $j$-th row, where
	$(x)_1$ or $(x)_2$ denote the first or second component of $x$. We
	say that a line $r_{j,1}$ ($r_{j,2}$) is {\it active} if
	$r_{j,1}\cap C(\h_{cl})\not = \emptyset$ ($r_{j,2}\cap
	C(\h_{cl})\not = \emptyset$).
	
	Let
	
	\be{proie1} \p_1(\h):=\{j \in \Z:\, r_{j,1}\cap
	C(\h_{cl})\not=\emptyset\} \ee
	
	\noindent and $p_1(\h):=|\p_1(\h)|$. In a similar way we define
	the vertical projection $\p_2(\h)$ and $p_2(\h)$. We also call
	$\p_1(\h)$ and $\p_2(\h)$ the horizontal and vertical {\it shadows
	} of $\h_{cl}$, respectively.
	
	Note that  $ g_1,g_2, \p_1, \p_2, p_1, p_2$ are actually depending
	on $\h$ only through $\h_{cl}$, even though, for notational
	convenience, we omit the subscript $_{cl}$ in their functional
	dependence.
	
	Note that $\h_{cl}$ is not necessarily a connected set and thus both
	the horizontal and vertical projections $\p_1(\h),\,\p_2(\h)$ are
	not in general connected. We have obviously:
	
	\be{g'} g'_i(\h):= g_i(\h)- p_i(\h)\ge 0. \ee
	
	\item[$\bullet$]
	A single cluster $C$ is called {\it monotone} if $g_i(C)= p_i(C)$
	for $i=1,2$, i.e.,
	$g_1$ and $g_2$ equal respectively  the horizontal and
	vertical side lengths of the rectangle $\cR(C)$ circumscribed to the
	unique cluster $C$. More generally,  we call {\it monotone} a
	configuration such that $g_i(\h)= p_i(\h)$ for $i=1,2$.

	Note that $s(\h)$ coincides with the semi-perimeter if $\h$ is
	a configuration with a single monotone cluster. It is immediate
	to show that $v(\h)$ is a non negative integer and  that it is
	equal to zero if $\h_{cl}$ has a  unique rectangular cluster with
	semi-perimeter $s(\h)$; it represents the number of vacancies in
	$\h$. Define:
	
	\be{r12} P_1(\h):=\bigcup_{j\in \p_1(\h)}r_{j,1}\qquad
	P_2(\h):=\bigcup_{j\in \p_2(\h)}r_{j,2} \ee
	
	\noindent the minimal unions of columns and rows, respectively, in
	$\Z^2$ containing $\h_{cl}$. By definition we have
	
	\be{vacancies1} P_1(\h)\cap P_2(\h)\supseteq \h_{cl}, \ee
	
	\noindent where $P_1(\h)\cap P_2(\h)$ is, in general, the union of
	rectangles such that $|P_1(\h)\cap P_2(\h) |=p_1(\h)p_2(\h)$. The
	{\it vacancies} of $\h$ are the sites in $P_1(\h)\cap
	P_2(\h)\backslash \h_{cl}.$

	\item[$\bullet$]
	Given a non-empty set $\cA\subset\cX$, define its
	{\it external and internal boundary\/} as, respectively
	
	\be{Abd}
	\partial^+ \cA:= \{\z\notin\cA\colon\; P(\z,\h)>0 \mbox{ for some } \h\in\cA\},
	\ee
	
	\be{Abd-}
	\partial^- \cA:= \{\z\in\cA\colon\; P(\z,\h)>0 \mbox{ for some } \h\notin\cA\}.
	\ee
	
	Moreover, let
	
	\be{bordoB}
	\partial\cA:=\{(\bar\h,\h):\bar\h\in\partial^{-}\cA,\,\h\in\partial^+\cA\,
	\hbox{ with } P(\bar\h,\h)>0\}, \ee
	
	\noindent
	be the set of moves exiting from $\cA$.

	We define
	
	\be{defHmin} H_{min}(\partial\cA):=
	\min_{(\bar\h,\h)\in\partial\cA}\{\max{\{H(\bar\h),H(\h)\}}\} \ee
	and we
	denote by $(\partial\cA)_{min}$ the subset of $\partial\cA$ where
	this minimum is realized:
	
	\be{defminbd} (\partial\cA)_{min}:= \{(\bar\h,\h)\in\partial\cA:\;
	\max{\{H(\bar\h),H(\h)\}}=H_{min}(\partial\cA)\}. \ee
	
\end{itemize}

\subsection{Reference path}
\label{refpathsec}
We recall (\ref{ld0}), (\ref{ld1}) and (\ref{ld2}) for the definitions of domino rectangles.

We will construct a particular set $\Omega^*$ whose elements are reference paths $\omega^*:\vuoto\ra \pieno$. Each path will be given by a particular sequence of growing domino rectangles, followed by a sequence of rectangles growing in the horizontal direction, followed by a sequence of rectangles growing in the vertical direction. The maximum of the energy along $\o^*$, $\{{arg \ max}_{\o^*}H\}$, is reached on particular configurations given by circumscribed rectangle $\cR(2l_2^*-1,l_2^*-1)$ as in Figure \ref{fig:P} on the left hand-side. 

We will prove in Corollary \ref{omegauppbou} that $\o^*\in{(\vuoto\ra\pieno)_{opt}}$ so that $\{{arg \ max}_{\o^*}H\}\in{\cS(\vuoto,\pieno)}$.

We want to recall that in this work we get only a partial solution to the problem of determination of the tube of typical paths, i.e., the set of paths followed by the process with high probability during the transition from $\vuoto$ to $\pieno$. Note that this set is much larger than $\Omega^*$; in the construction of the paths $\o^*$ we have a lot of freedom, so we choose this particular path from $\vuoto$ to $\pieno$ that suggests the structure of the tube of typical paths. 
\noindent
The idea behind the construction of the reference path is the following: we first construct a skeleton path $\{\bar{\omega}_{s}\}_{s=0}^{2L}$ given by a sequence of domino rectangles of semi-perimeter $s$. We point out that the transition from $\bar{\omega}_s$ to $\bar{\omega}_{s+1}$ can not be given in a single step, since  $\bar{\omega}_s$ and  $\bar{\omega}_{s+1}$ are rectangles and so this is not a path in that sense. Thus we have to interpolate each transition of the skeleton path in order to obtain a path. This is done in two different steps. First we introduce a sequence $\tilde{\omega}_{s,0}$, $\ldots$ , $\tilde{\omega}_{s,i_{s}}$ between $\bar{\omega}_s$ and $\bar{\omega}_{s+1}$, given by $\bar{\omega}_s$ plus a growing column. There are some cases (for $0$-domino rectangles) in which growing a column is equivalent to the operation of column to row from an energetic point of view. Since we consider a specific path, we choose one of the two operations arbitrarily. What we are doing is considering the time-reversal dynamic in the subcritical region and the usual dynamics in the supercritical region. The last interpolation consists of inserting between every pair of consecutive configurations in $\tilde{\omega}$, for which the cluster is increased by one particle, a sequence of configurations with one new particle created at the boundary of the box and brought to the correct site with consecutive moves of this free particle. In this way, from the sequence of configurations $\tilde{\omega}_{s,i}$, we obtain a path $\omega^{*}$, i.e., such that $\mathbb{P}(\omega^{*}_j, \ \omega^{*}_{j+1})>0$. \\ \\
$\mathbf{Skeleton: \ \bar{\omega}}$ \\
Let us construct a sequence of rectangular configurations $\bar{\omega}=\{\bar{\omega}_s\}$, with $s=0,\ \ldots, \ L$, such that $\bar{\omega}_1=\underline{0}$, $\bar{\omega}_2=\{x_0\}$, $\ldots$, $\bar{\omega}_{2L}=\mathcal{F}(\cX)\in{\underline{1}}$, where $x_0$ is a given site in $\Lambda_0$ and for every $s$, $\bar{\omega_s}\subset \bar{\omega}_{s+1}$. 
\begin{description}
	\item[Step a.] For any $s<3l_2^{*}-2$, $\{\bar{\omega}_s\}$ is a growing sequence of domino rectangles, depending on the value of $s$. If $[s]_3=[0]_3$, we have that $\bar{\omega}_s\in{\mathcal{R}}(2l_2,\ l_2)$ is a $0$-domino rectangle. If $[s]_3=[1]_3$ we have that $\bar{\omega}_s\in{\mathcal{R}}(2l_2-2,\ l_2)$ is a $1$-domino rectangle. If $[s]_3=[2]_3$ we have that $\bar{\omega}_s\in{\mathcal{R}}(2l_2-1,\ l_2)$ is a $2$-domino rectangle. 
	\item [Step b.] For any $3l_2^{*}-2\leq s\leq l_2^{*}+L-1$, $\{\bar{\omega}_s\}\in{\mathcal{R}}(s-l_2^{*}, \ l_2^{*})$.
	\item [Step c.] For any $s\geq l_2^{*}+L-1$, if $l_1=L-1$ we have $\bar{\omega}_s\in{\mathcal{R}}(L-1, \ s-L+1)$, and if $l_1=L$ we have $\bar{\omega}_s\in{\mathcal{R}}(L, \ s-L)$
\end{description}
$\mathbf{First \ interpolation: \ \tilde{\omega}}$ \\ 
Given a choice for $\bar{\omega}_s$, we can construct the path $\tilde{\omega}_{s, i}$ such that $\tilde{\omega}_{s, 0}=\bar{\omega}_s$ and insert between each pair $(\bar{\omega}_s, \ \bar{\omega}_{s+1})$ for any $s$ a sequence composed by configurations $\tilde{\omega}_{s, i}$ for $i=0, \ 1, \ \ldots, \ i_s$.
\begin{description}
	\item[Step a.1.] If $s<3l_2^{*}-2$ and $[s]_3=[1]_3$ add a column as in Figure \ref{fig:growc}, passing from $\tilde{\omega}_{s, 0}\in{\mathcal{R}}(2l_2, l_2+1)$ to $2$-domino rectangle $\tilde{\omega}_{s, i_s}\in{\mathcal{R}}(2l_2+1, l_2+1)$.
	\item[Step a.2.] If $s<3l_2^{*}-2$ and $[s]_3=[2]_3$ add a column as in Figure \ref{fig:growc}, passing from $\tilde{\omega}_{s, 0}\in{\mathcal{R}}(2l_2-1, l_2)$ to $0$-domino rectangle $\tilde{\omega}_{s, i_s}\in{\mathcal{R}}(2l_2, l_2)$.
	\item[Step a.3.] If $s<3l_2^{*}-2$ and $[s]_3=[0]_3$ add a column as in Figure \ref{fig:growc}, passing from $\tilde{\omega}_{s, 0}\in{\mathcal{R}}(2l_2, l_2)$ to $q$-domino rectangle $\tilde{\omega}_{s, i_s}\in{\mathcal{R}}(2l_2+1, l_2)$. Then use the path described in Figure \ref{fig:barO5} to define the path between $\tilde{\omega}_{s, l_2}\in{\mathcal{R}}(2l_2+1, l_2)$ to $1$-domino rectangle $\tilde{\omega}_{s, i_s}\in{\mathcal{R}}(2l_2, l_2+1)$.
	\item[Step b.1.] If $3l_2^{*}\leq s \leq l_2^{*}+L-1$ add a column as in Figure \ref{fig:growc}, passing from $\tilde{\omega}_{s, 0}\in{\mathcal{R}}(s-l_2^{*},  l_2^{*})$ to $\tilde{\omega}_{s, i_s}\in{\mathcal{R}}(s-l_2^{*}+1, l_2^{*})$.
	\item[Step c.1.] If $s\geq l_2^{*}+L-1$ and $l_1=L-1$, add a column as in Figure \ref{fig:growc}, passing from  $\tilde{\omega}_{s, 0}\in{\mathcal{R}}(L-1,  s-L+1)$ to $\mathcal{R}(L-1,  s-L+1)$, so we have $l_1=L$. Then use the path described in Figure \ref{fig:barO5} to obtain $\mathcal{R}(L-1, s-L+2)$.
\end{description}
$\mathbf{Second \ interpolation: \ \omega^{*}}$ \\ 
For any pair of configurations $(\tilde{\omega}_{s, i},  \tilde{\omega}_{s, i+1})$ such that $|\tilde{\omega}_{s, i}|<|\tilde{\omega}_{s, i+1}|$, by construction of the path $\tilde{\omega}_{s, i}$ the particles are created along the external boundary of clusters. Thus there exists $x_1,  \ldots, x_{j_{i}}$ connected chain of nearest-neighbor empty sites of $\tilde{\omega}_{s, i}$ such that $x_1\in{\partial ^{-}\Lambda}$ and $x_{j_{i}}$ is the site where the additional particle in $\tilde{\omega}_{s, i+1}$ is located. Hence we define $\omega_{s,  i,  0}^{*}=\tilde{\omega}_{s, i}$ and $\omega_{s, i,  j_{i}}^{*}=\tilde{\omega}_{s, i+1}$ for $s=0,  \ldots,  2(L+2)$. Otherwise, if $|\tilde{\omega}_{s, i}|=|\tilde{\omega}_{s, i+1}|$, we define $\omega_{s, i, 0}^{*}=\tilde{\omega}_{s, i}$ and  $\omega_{s,  i+1, 0}^{*}=\tilde{\omega}_{s, i+1}$.

We recall a useful lemma to compute the energy of different configurations:

\bl{lemma2.1}\cite[Lemma 7]{NOS}
For any configuration $\eta$: 
\be{energeta}
H(\eta)=H(\mathcal{R}(p_1(\eta),p_2(\eta)))+\varepsilon v(\eta)+U_1g_2'(\eta)+U_2g_1'(\eta)+n(\eta)\Delta,
\ee
with $\e$ as in (\ref{defepsilon}) and $g_i'$ as in (\ref{g'}).

\el

\begin{corollary}\label{energiapieno}
	For $L$ sufficiently large (say $L>\Big\lfloor\frac{U_1+U_2}{\e}\Big\rfloor$) we have $H(\pieno)<0=H(\vuoto)$.
\end{corollary}

The main property of the path $\o^*$ is the following.

\bp{refpath} If $U_1>2U_2$ and $L$ is large enough, we have that
\be{}
\{{\hbox{arg max}}_{\o^*}H\}\subseteq\o^*\cap\cP_1.
\ee
\ep
\bpr
Let us consider the skeleton path $\{\bar{\omega}_s\}_{s=0,..,2(L+2)}$ and let $\omega^*(\bar{\omega}_s,\bar{\omega}_{s+1})$ be the part of $\omega^*$ between $\bar{\omega}_s$ and $\bar{\omega}_{s+1}$. Letting
\be{}
g(s):=\max_{\h\in{\o^*}(\bar{\o}_s,\bar{\o}_{s+1})}H(\h),
\ee
we have
\be{}
\max_{\eta\in{\omega^*}}H(\h)=\max_{s=0,..,2(L+2)}g(s).
\ee
For the values of $s$ corresponding to steps $a.1$ ($s\leq3l_2^*-2$ and $[s]_3=[1]_3$), $a.2$ ($s\leq3l_2^*-2$ and $[s]_3=[2]_3$) and $b.1$ ($s>3l_2^*-2$), we can verify directly that $g(s)=H(\bar{\omega}_s)+2\Delta-U_1$, indeed $2\Delta-U_1$ is the energy barrier for adding a column. Now let us consider the case $s\leq3l_2^*-2$ and $[s]_3=[0]_3$, then the path described in step $a.3$ has a first part going from $\bar{\omega}_s$ to $\cR^{q-dom}(s+1)$, reaching its maximal value of energy in $H(\bar{\omega}_s)+2\Delta-U_1$. The second part of the path in the step $a.3$ goes from $\cR^{q-dom}(s+1)$ to $\bar{\omega}_{s+1}\in{\cR^{1-dom}(s+1)}$ with the operation column to row, so it reaches its maximal value of energy in
\be{} 
\ba{lll}
H(\cR^{q-dom}(s+1))+\D-U_2+U_1=\\
=H(\bar{\omega}_s)+2\Delta-U_1+\Delta-U_2+U_1-(U_1+U_2+\e(l_2-2))=\\
=H(\bar{\omega}_s)-\e\Big(\frac{s}{3}\Big)+\Delta+U_1, 
\ea\ee
where $U_1+U_2+\e(l_2-2)$ is the energy barrier for removing a column and it consists in the difference between $H(\bar{\omega}_s)+2\D-U_1$ and $H(\cR^{q-dom}(s+1)$ (see Figure \ref{fig:growc}). For the last equality we use the fact that $l_2(s)=\frac{s}{3}$ for a $0$-domino rectangle.

Explicit computations show that
\be{}
\max\{H(\bar{\omega}_s)+2\Delta-U_1, \ H(\bar{\omega}_s)-\e\frac{s}{3}+\Delta+U_1\}=H(\bar{\omega}_s)-\e\frac{s}{3}+\Delta+U_1.
\ee
Indeed, since $U_1-2U_2>0$ and we can write $s\leq3l_2^*-3$, because $[s]_3=[0]_3$ and thus $s\neq3l_2^*-2$, we have
\be{stepa.3}
\ba{lll}
H(\bar{\omega}_s)-\e\frac{s}{3}+\Delta+U_1\geq H(\bar{\omega}_s)-\e\Big(\frac{3l_2^*-3}{3}\Big)+\D+U_1=\\

=H(\bar{\omega}_s)-\e\Big(\Big\lceil{\frac{U_2}{\e}}\Big\rceil\Big)+\e+\D+U_1\geq H(\bar{\omega}_s)-U_2+\D+U_1>H(\bar{\omega}_s)+2\D-U_1 \\
\Leftrightarrow \quad \D<2U_1-U_2 \Leftrightarrow \quad U_1-2U_2>-\e. 
\ea
\ee

\noindent
We obtain
\be{}
g(s)=
\left \{ \ba{ll}
H(\bar{\o}_s)+2\D-U_1 & \mbox{if } s\leq3l_2^*-2 \mbox{ } \mbox{ and } [s]_3\neq[0]_3, \mbox{  or } s>3l_2^*-2, \\
H(\bar{\o}_s)-\e\frac{s}{3}+\D+U_1 & \mbox{if } s\leq3l_2^*-2 \mbox{ and } [s]_3=[0]_3.
\ea
\right.
\ee
We want now to evaluate the maximal value of $g(s)$ for $s\leq3l_2^*-2$. Let us consider the energy of domino configurations: 

\be{hdom} \ba{lll}
h^{(0-dom)}(n):=H(\cR^{0-dom}(3n))=U_1n+2U_2n-2\e n^2, \quad n=0,..,l_2^*-1,\\
h^{(1-dom)}(n):=H(\cR^{1-dom}(3n+1))=U_1(n+1)+2U_2n \\
\qquad \qquad \qquad \ \ -2\e n(n+1), \quad n=0,..,l_2^*-1,\\
h^{(2-dom)}(n):=H(\cR^{2-dom}(3n+2))=U_1(n+1)+U_2(2n+1) \\ \qquad \qquad \qquad \ \ -\e(n+1)(2n+1), \quad n=0,..,l_2^*-2 . \ea \ee
\bigskip
We observe that $h^{(0-dom)}(n)$ is an increasing function of $n$, indeed 
\be{hdomn}
\frac{dh}{dn}^{(0-dom)}(n)=U_1+2U_2-4\e n>0 \ \Leftrightarrow \ n<\frac{U_1+2U_2}{4\e}. 
\ee
In this case $n\leq l_2^*-1$ and, if $n=l_2^*-1$, we have that
\begin{center}
	$l_2^*-1<\frac{U_1+2U_2}{4\e} \ \Leftrightarrow \ \frac{2U_2-U_1}{\e}<4(1-\delta).$
\end{center} 

\noindent
Since $U_1<2U_2$ and $\d<1$, we have proved (\ref{hdomn}).
In a similar way we obtain that $h^{(1-dom)}(n)$ and $h^{(2-dom)}(n)$ are also increasing function of $n$. This implies that
\be{}
\max_{s\leq3l_2^*-2}g(s)=\max\{H(\bar{\o}_{3l_2^*-3})-\e(l_2^*-1)+\D+U_1, \ H(\bar{\o}_{3l_2^*-2})+2\Delta-U_1, \ H(\bar{\o}_{3l_2^*-4})+2\D-U_1\}.
\ee
By a direct comparison we obtain immediately
\be{maxg(s)}
\max_{s\leq3l_2^*-2}g(s)=H(\bar{\o}_{3l_2^*-3})-\e(l_2^*-1)+\D+U_1.
\ee
Indeed, since $\bar{\o}_{3l_2^*-2}\in\cR(2l_2^*-2,l_2^*)$, $\bar{\o}_{3l_2^*-3}\in\cR(2l_2^*-2,l_2^*-1)$ and $\bar{\o}_{3l_2^*-4}\in\cR(2l_2^*-3,l_2^*-1)$, we can write
\bi
\item[]  $H(\bar{\o}_{3l_2^*-2})=U_1l_2^*+2U_2(l_2^*-1)-2\e l_2^*(l_2^*-1)=U_1l_2^*+2U_2l_2^*-2U_2-2\e((l_2^*)^2-l_2^*)$,
\item[] $H(\bar{\o}_{3l_2^*-3})=U_1(l_2^*-1)+2U_2(l_2^*-1)-2\e(l_2^*-1)^2=U_1l_2^*-U_1+2U_2l_2^*-2U_2-2\e((l_2^*)^2-2l_2^*+1)$,
\item[]  $H(\bar{\o}_{3l_2^*-4})=U_1(l_2^*-1)+U_2(2l_2^*-3)-2\e(l_2^*-1)(2l_2^*-3)=U_1l_2^*-U_1+2U_2l_2^*-3U_2-\e(2(l_2^*)^2-5l_2^*+3)$.

\ei
\medskip
Since $U_1>2U_2$ and $\d<1$, we get $H(\bar{\o}_{3l_2^*-2})+2\D-U_1>H(\bar{\o}_{3l_2^*-4})+2\D-U_1$. Indeed
\be{}
\ba{lll}
H(\bar{\o}_{3l_2^*-2})+2\D-U_1<H(\bar{\o}_{3l_2^*-4})+2\D-U_1\\
\Leftrightarrow \ -2U_2+2\e l_2^*<-U_1-3U_2+5\e l_2^*-3\e \\
\Leftrightarrow \ l_2^*>\frac{U_1+U_2}{3\e}+1 \Leftrightarrow \ \frac{U_2}{\e}+\d>\frac{U_1+U_2}{3\e}+1 \Leftrightarrow \ \frac{U_1-2U_2}{\e}<3(\d-1).
\ea
\ee

\noindent
Since $\d>0$, we get $H(\bar{\o}_{3l_2^*-3})+\D+U_1-\e(l_2^*-1)>H(\bar{\o}_{3l_2^*-2})+2\D-U_1$, which concludes the proof of (\ref{maxg(s)}). Indeed
\be{}
\ba{lll}
H(\bar{\o}_{3l_2^*-3})+\D+U_1-\e(l_2^*-1)>H(\bar{\o}_{3l_2^*-2})+2\D-U_1 \\
\Leftrightarrow \ 4\e l_2^*-2\e+\D+U_1-\e l_2^*+\e>2\e l_2^*+2\D \\

\Leftrightarrow \ \e l_2^*-\e>\D-U_1 \Leftrightarrow \ \e\Bigg(\frac{U_2}{\e}+\d\Bigg)-\e>\D-U_1 \Leftrightarrow \ \e\d>0.

\ea\ee

\medskip
If $s>3l_2^*-2$ we have that $g(s)=H(\bar{\o}_s)+2\Delta-U_1$: for which value of $s$ do we obtain the maximum of the function $g(s)$? 

Since $\bar{\o}_s=\cR(s-l_2^*,l_2^*)$, we have that $H(\bar{\o}_s)=H(\cR(s-l_2^*,l_2^*))=U_1l_2^*+U_2(s-l_2^*)-\e l_2^*(s-l_2^*)$. By a direct computation we observe that $H(\bar{\o}_s)$ is a decreasing function of $s$, so the maximum value is reached for the minimum possible value of $s$, i.e., $s_0=3l_2^*-1=s^*$. Indeed we have
\be{}
\frac{dH}{ds}(\bar{\o}_s)=U_2-\e l_2^*<0.
\ee

Since the energy of the configurations in $T_3'$ can be made arbitrary small by choosing $L$ large enough, it remains only to compare the maximum values of $g(s)$ for $s\leq3l_2^*-2$ and $s>3l_2^*-2$. 
\be{argmax}
\ba{lll}
\max_{s=0,..,2(L+2)}g(s)=\max\Big\{\max_{s\leq3l_2^*-2}g(s),\max_{s=3l_2^*-1,..,2(L+2)}g(s)\Big\}=\\

=\max\{H(\bar{\o}_{3l_2^*-3})-\e(l_2^*-1)+\D+U_1, \ H(\bar{\o}_{s^*})+2\D-U_1\}=\\

=H(\bar{\o}_{3l_2^*-3})-\e(l_2^*-1)+\D+U_1. 
\ea
\ee

Since $l_2^*=\left\lceil{\frac{U_2}{\e}}\right\rceil$, set $\bar{\omega}_{s^*}=\cR(2l_2^*-1,l_2^*)$ and $\bar{\omega}_{3l_2^*-3}=\cR(2l_2^*-2,l_2^*-1)$, by a direct computation we obtain that $H(\bar{\o}_{3l_2^*-3})-\e(l_2^*-1)+\D+U_1>H(\bar{\o}_{s^*})+2\D-U_1$. Indeed
\be{}
\ba{lll}
H(\bar{\o}_{3l_2^*-3})-\e(l_2^*-1)+\D+U_1>H(\bar{\o}_{s^*})+2\D-U_1 \\
\Leftrightarrow \ 2U_1+2U_2-2\e-U_2+\e l_2^*<-2U_2+4\e l_2^*-2\e-\e l_2^*+U_1+U_2+U_1 \\
\Leftrightarrow \ 2U_2<2\e l_2^* \Leftrightarrow \ l_2^*>\frac{U_2}{\e}.
\ea
\ee

By the definition of $\o^*\in\Omega^*$, it is immediate to show that the configurations where the maximum value of the energy is reached are the configurations in $\o^*\cap\cP_1$.

\epr

\bc{omegauppbou}
We have 

$$\Phi(\vuoto,\pieno)\leq\G.$$
\ec

\bpr
By definition of communication height $\Phi(\vuoto,\pieno)$, given a path $\o^*\in\Omega^*$ by
Proposition \ref{refpath}, (\ref{squaredefempty}) and the sentence below, we have immediately

\be{e3.1} \Phi(\vuoto,\pieno)\leq \max_i H(\o^*_i) = H(\cP_1)=\G, \ee

where $\cP_1$ and $\G$ are defined in (\ref{defP1}) and (\ref{g}) respectively. The first inequality follows because $\Phi(\vuoto,\pieno)$ is the minimum over all paths between $\vuoto$ and $\pieno$, and we bound it by taking the reference path $\o^*$.

\epr
\subsection{Definition of the set $\cB$ and exit from $\cB$}
\label{exit}
We give the definition of the set $\cB$ that is a basin of attraction of $\vuoto$ and satisfies the strategy following \cite[Section 4.2]{MNOS}.
\bd{B}
We define the set $\cB$ as follows:

\be{defB}
\cB:=\left\{ \eta:
\ba{lll}
s(\h)&\leq& s^*-2,  \hbox{ or }\\
s(\h)&\geq&s^*-1 \hbox{ and }
p_2(\h)\leq l_2^*-1,\hbox{ or }\\
s(\h)&=&s^*-1, \hbox{ } p_2(\h)\geq l_2^* \hbox{ and } v(\h)\geq p_{min}(\h)-1,  \hbox{ or} \\
s(\h)&\geq&s^*, \hbox{ } p_2(\h)=l_2^* \hbox{ and } v(\h)\geq p_{max}(\h)-1
\ea
\qquad\right \}, \ee

\noindent where $s^*$ is defined in (\ref{defs^*}), $p_{min}(\h)= \min{\{p_1(\h),p_2(\h)\}}$ and $p_{max}(\h)= \max{\{p_1(\h),p_2(\h)\}}$.
\ed

In \cite{LPr} there is the construction of the probability measure using the restricted ensemble $\m_{\cR}=\m(\cdot|\cR)$ while in \cite{CGOV} and \cite[Sections 4.2, 6.3]{OV} the authors consider the restricted dynamics associated to the exit time $T_{\cX\setminus\cR}$ for some $\cR$ that can be thought as basin of attraction of the metastable state. The set $\cB$ defined in (\ref{defB}) satisfies the hypotheses of this set $\cR$ used in \cite{BGM}.

\bl{lemmaB}
For $\cB$ as in (\ref{defB}) we get $\vuoto\in\cB$ and $\pieno\notin\cB$.

\el

\bpr
Since $s(\vuoto)=0$ it immediately follows that $\vuoto\in\cV_{\leq s^*-2}$ and thus we get $\vuoto\in\cB$. Similarly, since $\pieno\in\cV_{\geq s^*}$ and $p_2(\pieno)>l_2^*$, we get $\pieno\notin\cB$.
\epr
\bigskip

The main result of this subsection is given by the following:

\bt{3.1} For  $H_{min}(\partial\cB)$ as in (\ref{defHmin}),
$(\partial\cB)_{min}$ as in (\ref{defminbd}),
$\G$ as in (\ref{g}) and $\cP$ as in (\ref{defP}), we have

\be{H^*} H_{min}(\partial\cB)=\G. \ee

\noindent
In addition
\bi
\item[(i)]  if $(\bar\h,\h)\in(\partial\cB)_{min}$, then $H(\bar\h)\ge H(\h)$ and  $\bar\h\in\cP_2$;
 
\item[(ii)] if $\bar\h\in\cP_1$, there exists at least a path $\bar\o\in(\vuoto\ra\pieno)_{opt}$ such that $\bar\o=(\vuoto,\bar\o_1,..,\bar\o_{j-1},\bar\h,\bar\o_{j+1},..,\pieno)$. Moreover, we have that $\bar\o_{j+1}\in\cB$.
\ei
\et

\bc{corollario1}
For any path $\o^*\in\Omega^*$ we have

$$\{{\arg \max}_{\o^*}H\}\cap(\partial\cB)_{min}\neq\emptyset.$$
\ec

\bpr
Combining Proposition \ref{refpath} and Theorem \ref{3.1}, we directly get the conclusion.

\epr

\bc{phi} We have

$$\Phi(\vuoto,\pieno)=\G.$$

\ec
\bpr
Since every path going from $\vuoto$ to $\pieno$ has to leave $\cB$, we have by Theorem \ref{3.1} that

\be{e3.2}
\Phi(\vuoto,\pieno):=\min_{\o:\vuoto\to\pieno}\max_{\zeta\in\o}H(\zeta)\ge
H_{min}(\partial\cB)=\G. \ee

Combining (\ref{e3.2}) and Corollary \ref{omegauppbou}, we get $\Phi(\vuoto,\pieno)=\G$.
\epr

Note that to prove Corollary \ref{phi} we have proposed a suitable set $\cB$ (see (\ref{defB})) and we have applied the argument developed in \cite[Section 4.2]{MNOS} with some small variations. In \cite{MNOS} the set $\partial^+\cB$ (external boundary of $\cB$, see (\ref{Abd})) was considered, while here we use the set $\partial\cB$ (see (\ref{bordoB})), so in the present case $H_{min}(\partial\cB)$ substitutes $H(\cF(\partial^+\cB))$.

\bp{V^*} There exists $V^*<\Gamma$ such that $\mathcal{X}_{V^{*}}\subseteq \{\underline{0}, \ \underline{1}\}$, i.e., $\forall \ \eta\neq \underline{0}$, $\underline{1}$ there exists $\eta'\in{\mathcal{X}}$ and a path $\omega:\eta\rightarrow \eta'$ such that $H(\eta')<H(\eta)$ and $V_{\h}\leq\max_{\xi\in{\omega}}H(\xi)-H(\h)\leq V^*$.
\ep

In this proposition we prove that that each configuration $\eta\notin \{\underline{0}, \underline{1}\}$ is $V^*$-reducible, namely we can find a configuration $\eta'\in \mathcal{I}_{\eta}$ with smaller energy and $\Phi(\eta, \ \eta')\leq H(\eta)+V^*$. Thus $V^*$ consists in an upper bound for the quantity $V_{\h}$ for any $\h\notin\{\vuoto,\pieno\}$ and so we can write $V^*=\max\{V_\h \ |\ \h\in\cX\setminus\{\vuoto,\pieno\}\}$. In other words, we assert that there are no too deep wells in the energy landscape, i.e., no deeper than the well with bottom $\vuoto$. Moreover, we observe that being a metastable state is equivalent to the absence of too deep energy.
For the proof of Proposition \ref{V^*} we refer to subsection \ref{reduction}.

\subsection{Proof of the main theorems \lowercase{{\ref{t1}}} and \lowercase{\ref{t1'}}}
\label{timegate}

The proof of theorems \ref{t1} and \ref{t1'} consist of an application of theorems \cite[Theorem 4.1]{MNOS}, \cite[Theorem 4.9]{MNOS}, \cite[Theorem 4.15]{MNOS} for the asymptotics of the tunnelling time and Theorem \cite[Theorem 5.4]{MNOS} for the gates in the general setup of reversible Markov chains. 

\bigskip
\begin{proof*}{\bf of Theorem \ref{t1}}

	Combining theorems \ref{metstate}, \cite[Theorem 4.1]{MNOS}, \cite[Theorem 4.9]{MNOS}, \cite[Theorem 4.15]{MNOS} and Corollary \ref{phi}, we get the conclusion.
	
\end{proof*}

\bigskip
\begin{proof*}{\bf of Theorem \ref{t1'}}
	
	Given any optimal path $\o\in(\vuoto\ra\pieno)_{opt}$, since $\vuoto\in\cB$ and $\pieno\notin\cB$ by Lemma \ref{lemmaB}, $\o$ has to leave the set $\cB$ and this has to be done with a pair of configurations $(\bar\h,\h)\in(\partial\cB)_{min}$. If this is not satisfied, by Theorem \ref{3.1} we would have $\max_{i}H(\o_i)>\G$ and, since $\Phi(\vuoto,\pieno)=\G$, the path $\o$ would not be optimal, which is a contradiction. By Theorem \ref{3.1} we obtain that $\o$ must intersect the set $\cP$ and, thus the set $\cP$ is a gate for the transition $\vuoto\ra\pieno$; see subsection \ref{defind} for the definition of a gate. Thus (\ref{gate}) is proven using \cite[Theorem 5.4]{MNOS}.
	
\end{proof*}

\section{Proof of the main Theorem \lowercase{\ref{3.1}}}
\label{proofs}
In this Section we give the proof of the main Theorem \ref{3.1}.
In particular, in order to
analyze the exiting move from $\cB$, we prove preliminary results on single moves.

\medskip
{\bf Definitions and notations.} We first give important definitions and notations that will be useful in the section.

\bd{move}
We say that $(\bar\h,\h)$ is a {\it move} if $(\bar\h,\h)\in\cX\times\cX$ and $P(\bar\h,\h)>0$.

We say that a move $(\bar\h,\h)$ is an {\it exiting move from $\cB$} if $\bar\h\in\cB$, $\h\notin\cB$ and $\P(\bar\h,\h)>0$.

We define

\be{deltas} \D s:=s(\h)-s(\bar\h) \ee

the {\it variation of $s$} in the move.

\ed

We recall the definition of active line that is in the bullet before (\ref{proie1}). We say that a line $r$ (column or row in $\Z^2$) {\it becomes active} in the move from $\bar\h$ to $\h$ if it was not active in $\bar\h$ ($ r\cap\bar\h_{cl}=\emptyset$ see (\ref{hcl}))
and it is active in $\h$ ($r\cap\h_{cl}\neq \emptyset$). Furthermore, we say that a line $r$ {\it becomes inactive} in the move from $\bar\h$ to $\h$ if it was active in $\bar\h$ and it is not active in $\h$. If a line does not become active nor inactive we say that it does not change its behavior.
We will call $x_1$ the site containing the moving particle, $x_2$ the site containing the particle after the move and $x_3$, $x_4$ and $x_5$ the nearest neighbor sites of $x_2$ together with $x_1$ (see Figure \ref{fig:rappresentazione}). The sites $y_1$, $y_2$ and $y_3$ are the nearest neighbors of the particle in $x_1$ and $z_1$, $z_2$ and $z_3$ are the nearest neighbors of the site $y_3$ together with the site $x_1$. Furthermore we set $t$ the site above $y_1$ and $s$ the site under $y_2$.

We will call $r_1$ the line of the move, $r_2$ the line orthogonal to it passing trough site $x_2$, and $r_3$, $r_4$ and $r_5$ the lines passing through the three nearest neighbor sites of $x_2$, i.e., $x_3$, $x_4$ and $x_5$ respectively. Furthermore we will call $r_6$ the line orthogonal to $r_1$ passing trough the site $x_1$ and $r_7$ the line orthogonal to $r_1$ passing trough the site $y_3$.

\bigskip
\setlength{\unitlength}{0.8pt}
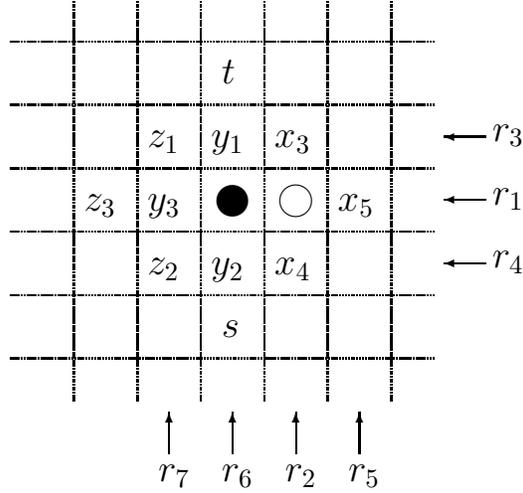
\begin{figure}
	\begin{center}
		\begin{picture}(380,60)(-30,40)
		
		\thinlines 
		\qbezier[150](90,0)(190,0)(290,0)
		\qbezier[150](90,30)(190,30)(290,30) \qbezier[150](90,60)(190,60)(290,60)
		\qbezier[150](90,90)(190,90)(290,90)
		\qbezier[150](90,-30)(190,-30)(290,-30)
		\qbezier[150](90,-60)(190,-60)(290,-60)
		\qbezier[150](240,-80)(240,15)(240,110)
		\qbezier[150](270,-80)(270,15)(270,110)
		\qbezier[150](120,-80)(120,15)(120,110)
		\qbezier[150](150,-80)(150,15)(150,110)
		\qbezier[150](180,-80)(180,15)(180,110)
		\qbezier[150](210,-80)(210,15)(210,110)
		\put(195,15){\circle*{15}}
		\put(225,15){\circle{15}}
		\begin{Large}
		\put(245,10){$x_5$}
		\put(215,40){$x_3$}
		\put(215,-20){$x_4$}
		\put(185,40){$y_1$}
		\put(185,-20){$y_2$}
		\put(190,70){$t$}
		\put(190,-50){$s$}
		\put(155,10){$y_3$}
		\put(155,40){$z_1$}
		\put(155,-20){$z_2$}
		\put(125,10){$z_3$}
		
		\thinlines
		\put(225,-105){\vector(0,1){20}} \put(220,-118){$r_2$}
		\put(165,-105){\vector(0,1){20}} \put(160,-118){$r_7$}
		
		\put(250,-118){$r_5$} \put(255,-105){\vector(0,1){20}}
		\put(195,-105){\vector(0,1){20}} \put(190,-118){$r_6$}
		\put(315,-15){\vector(-1,0){20}} \put(318,-17){$r_4$}
		\put(315,15){\vector(-1,0){20}} \put(318,13){$r_1$}
		\put(315,45){\vector(-1,0){20}} \put(318,43){$r_3$}
		\end{Large}
		\end{picture}
		\vskip 4.3 cm
		\caption{In order to help the reader we depict part of the configuration of $\bar\h$.}
		\label{fig:rappresentazione}
	\end{center}
\end{figure}

In Figure \ref{fig:rappresentazione} and consequently in the rest of the paper, we assume that the move is horizontal, i.e., $r_1$ is an horizontal line. Since the binding energies are not equal in the vertical and horizontal directions, one could be tempted to consider also the case $r_1$ vertical, but in the latter we can conclude with analogue arguments, except in the Proposition \ref{lemma6}, where we distinguish $r_1$ vertical or horizontal (case (h) and case (v)). We observe that the horizontal case is the worst case scenario.

\medskip
\noindent
We recall a useful lemma.

\bl{lemma3.2}\cite[Lemma 12]{NOS} Let $p_{min}(\h)\ge 4$, we have

\item{(i)}
$|\D s|\le 5$

\item{(ii)}
if $\D s = 1$ then $v(\h)\ge p_{min}(\h)-3$

if $\D s = 2$ then $n(\bar\h)\ge 1$ and $v(\h)\ge 2p_{min}(\h)-5$

if $\D s = 3$ then $n(\bar\h)\ge 2$ and $v(\h)\ge 3p_{min}(\h)-6$

if $\D s = 4$ then $n(\bar\h)\ge 3$ and $v(\h)\ge 4p_{min}(\h)-7$

if $\D s = 5$ then $n(\bar\h)\ge 4$ and $v(\h)\ge 5p_{min}(\h)-8$

\item{(iii)}
if $\D s=1$ and $v(\h)<p_{min}(\h)-1$ then $n(\bar\h)\ge 2$ \el

\br {remark1} We note that if $n(\bar\h)=0$ the unique line that
can become active is $r_2$ and in this case in $\bar\h$ sites
$x_3$ and $x_4$ are empty and $x_5\in\bar\h_{cl}$, where $x_5\in
r_1$ (see Figure \ref{fig:rappresentazione}), so that $\bar\h$ is not monotone, i.e., $g'_1(\bar \h)+g'_2(\bar \h)\ge 1$.
\er

\br {remarkulteriore} We note that if $n(\bar\h)=0$, $g_1'(\bar\h)=1$ and $g_2'(\bar\h)=0$ then the site $x_5$ must be empty, otherwise $g_2'(\bar\h)\geq1$. Thus in this case we obtain that no lines can become active.
\er

\noindent
Now we give an estimate of $s(\bar\h)$ for each possible value of $\D s\in{\{-5,-4,..,4,5\}}$. This will be useful for the entire section.
\br{remark2} Let $(\bar\h,\h)\in\partial\cB$ be an exiting move from $\cB$ (see definition \ref{move}): 

\bi
\item If $\D s=-5$ the exiting move is admissible only in the case $s(\bar\h)\geq s^*+4$. Indeed, if $s(\bar\h)\leq s^*+3$, then $s(\h)=s(\bar\h)-5\leq s^*-2$, which implies $\h\in\cB$.
\item If $\D s=-4$ the exiting move is admissible only in the case $s(\bar\h)\geq s^*+3$. Indeed, if $s(\bar\h)\leq s^*+2$, then $s(\h)=s(\bar\h)-4\leq s^*-2$, which implies $\h\in\cB$.
\item If $\D s=-3$ the exiting move is admissible only in the case $s(\bar\h)\geq s^*+2$. Indeed, if $s(\bar\h)\leq s^*+1$, then $s(\h)=s(\bar\h)-3\leq s^*-2$, which implies $\h\in\cB$.
\item If $\D s=-2$ the exiting move is admissible only in the case $s(\bar\h)\geq s^*+1$. Indeed, if $s(\bar\h)\leq s^*$, then $s(\h)=s(\bar\h)-2\leq s^*-2$, which implies $\h\in\cB$.
\item If $\D s=-1$ the exiting move is admissible only in the case $s(\bar\h)\geq s^*$. Indeed, if $s(\bar\h)\leq s^*-1$, then $s(\h)=s(\bar\h)-1\leq s^*-2$, which implies $\h\in\cB$.
\item If $\D s=0$ the exiting move is admissible only in the case $s(\bar\h)\geq s^*-1$. Indeed, if $s(\bar\h)\leq s^*-2$, then $s(\h)=s(\bar\h)\leq s^*-2$, which implies $\h\in\cB$.
\item If $\D s=1$ the exiting move is admissible only in the case $s(\bar\h)\geq s^*-2$. Indeed, if $s(\bar\h)\leq s^*-3$, then $s(\h)=s(\bar\h)+1\leq s^*-2$, which implies $\h\in\cB$.
\item If $\D s=2$ the exiting move is admissible only in the case $s(\bar\h)\geq s^*-3$. Indeed, if $s(\bar\h)\leq s^*-4$, then $s(\h)=s(\bar\h)+2\leq s^*-2$, which implies $\h\in\cB$.
\item If $\D s=3$ the exiting move is admissible only in the case $s(\bar\h)\geq s^*-4$. Indeed, if $s(\bar\h)\leq s^*-5$, then $s(\h)=s(\bar\h)+3\leq s^*-2$, which implies $\h\in\cB$.
\item If $\D s=4$ the exiting move is admissible only in the case $s(\bar\h)\geq s^*-5$. Indeed, if $s(\bar\h)\leq s^*-6$, then $s(\h)=s(\bar\h)+4\leq s^*-2$, which implies $\h\in\cB$.
\item If $\D s=5$ the exiting move is admissible only in the case $s(\bar\h)\geq s^*-6$. Indeed, if $s(\bar\h)\leq s^*-7$, then $s(\h)=s(\bar\h)+5\leq s^*-2$, which implies $\h\in\cB$.
\ei
\er
We can justify the previous remark looking at the definition of the set $\cB$.

\subsection{Main propositions}
\label{mainprop}
In this subsection we give the proof of the main Theorem \ref{3.1}, subdivided in the several following propositions, where we emphasize the cases in which $\h$ is in $\cB$, or $\h$ is not in $\cB$ but $H(\bar\h)>\G$, or the cases in which $\h$ is not in $\cB$ but $H(\bar\h)=\G$ and $\bar\h\in\cP$. We recall (\ref{defP}), (\ref{defB}) and (\ref{deltas}) for the definitions of $\cP$, $\cB$ and $\D s$ respectively.

\bp{lemma2}
Let $(\bar\h,\h)$ be a move with $\bar\h\in\cB$. If $\D s\leq-2$, then either $\h\in\cB$ or $\h\notin\cB$ and $H(\bar\h)>\G$.
\ep

We refer to subsection \ref{dim2} for the proof of the Proposition \ref{lemma2}.

\bp{lemma3}
Let $(\bar\h,\h)$ be a move with $\bar\h\in\cB$. If $\D s\geq-1$ and $p_2(\bar\h)\leq l_2^*-1$, then either $\h\in\cB$ or $\h\notin\cB$ and $H(\bar\h)>\G$.
\ep

We refer to subsection \ref{dim3} for the proof of the Proposition \ref{lemma3}.

\noindent
The following results consider the case $p_2(\bar\h)\geq l_2^*$, which we subdivide in different propositions because of the lengthy proof. Note that in some of these propositions we identify the set $\cP$.

\bp{lemma4}
Let $(\bar\h,\h)$ be a move with $\bar\h\in\cB$. If $\D s\geq3$, $p_{min}(\h)\geq4$ and $p_2(\bar\h)\geq l_2^*$, then either $\h\in\cB$, or $\h\notin\cB$ and $H(\bar\h)>\G$.
\ep

We refer to subsection \ref{dim4} for the proof of the Proposition \ref{lemma4}.

\bp{lemma5}
Let $(\bar\h,\h)$ be a move with $\bar\h\in\cB$. If $\D s=-1$ and $p_2(\bar\h)\geq l_2^*$, then we have one of the following:

\bi
\item[(i)] either $\h\in\cB$;
\item[(ii)] or $\bar\h\in\cP_1$, $H(\bar\h)=\G$ and $\h\in\cB$;
\item[(iii)] or $\h\notin\cB$ and $\max\{H(\bar\h),H(\h)\}>\G$, with $\bar\h\notin\cP$.
\ei
\ep

We refer to subsection \ref{dim5} for the proof of the Proposition \ref{lemma5}.

\bp{lemma6}
Let $(\bar\h,\h)$ be a move with $\bar\h\in\cB$. If $\D s=0$ and $p_2(\bar\h)\geq l_2^*$, then we have one of the following:

\bi
\item[(i)] either $\h\in\cB$;
\item[(ii)] or $\h\notin\cB$ and $H(\bar\h)=\G$, with $\bar\h\in\cP_2$;
\item[(iii)] or $\h\notin\cB$ and $H(\bar\h)>\G$, with $\bar\h\notin\cP$.
\ei
\ep

We refer to subsection \ref{dim6} for the proof of the Proposition \ref{lemma6}.

\bp{lemma7}
Let $(\bar\h,\h)$ be a move with $\bar\h\in\cB$. If $\D s=1$ and $p_2(\bar\h)\geq l_2^*$, then we have one of the following:

\bi
\item[(i)] either $\h\in\cB$;
\item[(ii)] or $\h\notin\cB$ and $H(\bar\h)=\G$, with $\bar\h\in\cP_2$;
\item[(iii)] or $\h\notin\cB$ and $H(\bar\h)>\G$, with $\bar\h\notin\cP$.
\ei
\ep

We refer to subsection \ref{dim7} for the proof of the Proposition \ref{lemma7}.

\bp{lemma8}
Let $(\bar\h,\h)$ be a move with $\bar\h\in\cB$; if $\D s=2$, $p_{min}(\h)\geq4$ and $p_2(\bar\h)\geq l_2^*$ then either $\h\in\cB$ or $\h\notin\cB$ and $\max\{H(\bar\h),H(\h)\}>\G$.
\ep

We refer to subsection \ref{dim8} for the proof of the Proposition \ref{lemma8}.

\bigskip
\begin{proof*}{\bf of the main Theorem \ref{3.1}.} 

Let $(\bar\h,\h)\in\partial\cB$ be the exiting move from $\cB$ and $\D s$ be its corresponding variation of $s$. If $p_{min}(\h)\leq 3$, for $\e\ll U_2$, from explicit computations follows that

\be{} H(\h)>\G. \ee

\noindent
From now on we assume that $p_{min}(\h)\geq 4$.

Combining Lemma \ref{lemma3.2} and Propositions \ref{lemma2}, \ref{lemma3}, \ref{lemma4}, \ref{lemma5}, \ref{lemma6}, \ref{lemma7} and \ref{lemma8} we complete the proof. In particular, using the proof of Proposition \ref{lemma5}{\it (ii)} below, Remark \ref{-1} with the proof and Proposition  \ref{refpath}, we can get point {\it (i)} and {\it (ii)}.

\end{proof*}

\subsection{Useful Lemmas}
\label{lemmi}
In this subsection we give some lemmas that will be useful for the proof of the propositions reported in subsection \ref{mainprop}. The proof of the lemmas is postponed to subsection \ref{dimlemmi}. Recall def. \ref{move} for the definition of a move.

\bl{decrescita}
Let $(\bar\h,\h)\in\cX\times\cX$ be a move, then each horizontal (resp. vertical) line becoming inactive decreases the vertical (resp. horizontal) projection $p_2(\bar\h)$ (resp. $p_1(\bar\h)$).
\el

\br{crescita}
With a similar argument as in the proof of Lemma \ref{decrescita}, we note that each horizontal (resp. vertical) line becoming active increases the vertical (resp. horizontal) projection $p_2(\bar\h)$ (resp. $p_1(\bar\h)$).
\er

\bl{lemma1}
In a single move the lines $r_2$ and $r_5$ can not become inactive and the lines $r_6$ and $r_7$ can not become active. 
\el

\bl{lemma0}
Let $(\bar\h,\h)$ be a move. We have 

\bi
\item [(i)] if the line $r_1$ becomes inactive with the move, then no line can become active;
\item [(ii)] if the line $r_1$ becomes active with the move, then the lines that can become inactive are $r_3$ and $r_4$.
\ei
\el

\bl{ulteriore}
Let $(\bar\h,\h)$ be a move with $\bar\h\in\cB$ and $\D s$ be its corresponding variation of $s$, with $\D s\geq-2$ and $s(\h)\geq s^*-1$. If $p_2(\h)\geq l_2^*$ and either the line $r_3$ (respectively $r_4$) becomes active in the move or $r_1$ becomes active in the move and $x_5$ is empty, we have 

\bi

\item [(i)] if $s(\h)=s^*-1$, then $\h\in\cB$;
\item [(ii)] if $s(\h)\geq s^*$ and $p_2(\h)=l_2^*$, then $\h\in\cB$.

\ei
\el

\bl{ult2}
Let $(\bar\h,\h)$ be a move with $\bar\h\in\cB$ and $\D s$ be its corresponding variation of $s$, with $\D s\geq-2$. If $p_2(\bar\h)\leq l_2^*-1$, then

\bi
\item [(i)] if one line among $r_3$ and $r_4$ becomes active and the line $r_1$ does not become active in the move, then $\h\in\cB$;
\item [(ii)] if only two horizontal lines become active in the move and either $p_2(\bar\h)\leq l_2^*-2$ or $p_2(\bar\h)=l_2^*-1$ and $s(\h)=s^*-1$, then $\h\in\cB$.
\ei
\el

\subsection{Proof of Proposition \lowercase{\ref{lemma2}}}
\label{dim2}

\bpr
By Lemma \ref{lemma3.2}{\it (i)} we can distinguish four different cases corresponding to $\D s=-5,-4,-3,-2$. We analyze separately each case.

\medskip
\noindent
{\bf Case $\D s=-5$.} Let $\D s=-5$ and $(\bar\h,\h)\in\partial\cB$ is the exiting move from $\cB$, i.e., $\bar\h\in\cB$, $\h\notin\cB$ and $\P(\bar\h,\h)>0$. By Remark \ref{remark2} we may consider only the case $s(\bar\h)\geq s^*+4$. Since $\bar\h\in\cB$, by (\ref{defB}) and $s(\bar\h)\geq s^*+4$ we have only the following two cases:

\bi
\item [(a)] $p_2(\bar\h)=l_2^*$ and $v(\bar\h)\geq p_{max}(\bar\h)-1$;
\item [(b)] $p_2(\bar\h)\leq l_2^*-1$.
\ei

Since by Lemma \ref{lemma1} the lines $r_2$ and $r_5$ cannot become inactive, in order to obtain $\D s=-5$ necessarily five lines become inactive and these lines are $r_1$, $r_3$, $r_4$, $r_6$ and $r_7$. Among them, three are horizontal and two are vertical. By Lemma \ref{decrescita}, we get
\be{5.1}
\left\{ \ba{lll}
p_1(\h)=p_1(\bar\h)-2, \\
p_2(\h)=p_2(\bar\h)-3.
\ea\right.
\ee
\noindent
In both cases (a) and (b) we consider $p_2(\bar\h)\leq l_2^*$, thus by (\ref{5.1}) we obtain that $p_2(\h)\leq l_2^*-3\leq l_2^*-1$. Thus we can conclude that $\h\in\cB$, i.e., in this case it is impossible to leave the set $\cB$.

\medskip
\noindent
{\bf Case $\D s=-4$.} Let $\D s=-4$ and $(\bar\h,\h)\in\partial\cB$ is the exiting move from $\cB$. By Remark \ref{remark2} we can consider only the case $s(\bar\h)\geq s^*+3$. Again, since $\bar\h\in\cB$, by (\ref{defB}) and $s(\bar\h)\geq s^*+3$, we only have the following two cases:

\bi
\item [(a)] $p_2(\bar\h)=l_2^*$ and $v(\bar\h)\geq p_{max}(\bar\h)-1$;
\item [(b)] $p_2(\bar\h)\leq l_2^*-1$.
\ei

We observe that in order to obtain $\D s=-4$ we have the two following possibilities:

\begin{description}
	\item [Case I.] four lines becoming inactive and no line become active;
	\item [Case II.] five lines becoming inactive and one line becoming active.
\end{description}

\noindent
{\bf Case I.} By Lemma \ref{lemma1} we have that at least one horizontal line becomes inactive, so by Lemma \ref{decrescita} in both cases (a) and (b) we get $p_2(\h)\leq p_2(\bar\h)-1\leq l_2^*-1$ and thus it is impossible to leave $\cB$.

\medskip
\noindent
{\bf Case II.} By Lemma \ref{lemma1} the lines becoming inactive are $r_1$, $r_3$, $r_4$, $r_6$ and $r_7$, so the line becoming active is either $r_2$ or $r_5$. If the line that becomes active is $r_5$, then the site $x_5$ must contain a free particle in $\bar\h$. Thus the line $r_1$ can not become inactive, so this case is not admissible (see Figure \ref{fig:rappresentazione}). If the line becoming active is $r_2$, in at least one site among $x_3$, $x_4$ and $x_5$ there must be a free particle in $\bar\h$. If the free particle is in $x_3$ or in $x_4$, then it cannot happen that both lines $r_3$ and $r_4$ become inactive. Thus we consider the case in which the free particle is in $x_5$. As in the case in which $r_5$ becomes active, the line $r_1$ can not become inactive. Thus this case is not admissible.

\medskip
\noindent
{\bf Case $\D s=-3$.} Let $\D s=-3$ and $(\bar\h,\h)\in\partial\cB$ is the exiting move from $\cB$. By Remark \ref{remark2} we may consider only the case $s(\bar\h)\geq s^*+2$. Again, since $\bar\h\in\cB$, by (\ref{defB}) and $s(\bar\h)\geq s^*+2$, we only have the following two cases:

\bi
\item [(a)] $p_2(\bar\h)=l_2^*$ and $v(\bar\h)\geq p_{max}(\bar\h)-1$;
\item [(b)] $p_2(\bar\h)\leq l_2^*-1$.
\ei

We note that in order to obtain $\D s=-3$ we have the three following possibilities:

\begin{description}
	\item [Case I.] five lines becoming inactive and two lines becoming active; 
	\item [Case II.] four lines becoming inactive and one line becoming active; 
	\item [Case III.] three lines becoming inactive and no line becoming active.
\end{description}

\noindent
{\bf Case I.} By Lemma \ref{lemma1} the lines that must become inactive are $r_1$, $r_3$, $r_4$, $r_6$ and $r_7$: thus the two lines becoming active are $r_2$ and $r_5$ (both vertical lines). Thus by Lemma \ref{decrescita} in both cases (a) and (b) we have $p_2(\h)=p_2(\bar\h)-3\leq l_2^*-3\leq l_2^*-1$, so it is impossible to leave $\cB$.

\medskip
\noindent
{\bf Case II.} Again by Lemma \ref{lemma1} the lines becoming inactive are four among $r_1$, $r_3$, $r_4$, $r_6$ and $r_7$. Every choice of four of these lines includes at least two horizontal lines, so, even though one horizontal line becomes active, by Lemma \ref{decrescita} we obtain $p_2(\h)\leq p_2(\bar\h)-1\leq l_2^*-1$ in both cases (a) and (b). Thus it is impossible to leave $\cB$.

\medskip
\noindent
{\bf Case III.} Three lines must become inactive. By Lemma \ref{lemma1} we know that at least one horizontal line must become inactive. Since in both cases (a) and (b) $p_2(\bar\h)\leq l_2^*$, by Lemma \ref{decrescita} we get $p_2(\h)\leq p_2(\bar\h)-1\leq l_2^*-1$. Thus it is impossible to leave $\cB$.

\medskip
\noindent
{\bf Case $\D s=-2$.} Let $\D s=-2$ and $(\bar\h,\h)\in\partial\cB$ is the exiting move from $\cB$. By Remark \ref{remark2} we can only consider the case $s(\bar\h)\geq s^*+1$. Again, since $\bar\h\in\cB$, by (\ref{defB}) and $s(\bar\h)\geq s^*+1$, we only have the following two cases:

\bi
\item [(a)] $p_2(\bar\h)=l_2^*$ and $v(\bar\h)\geq p_{max}(\bar\h)-1$;
\item [(b)] $p_2(\bar\h)\leq l_2^*-1$.
\ei

We note that in order to obtain $\D s=-2$ we have the three following possibilities:

\begin{description}
	\item [Case I.] four lines becoming inactive and two lines becoming active; 
	\item [Case II.] three lines becoming inactive and one line becoming active, 
	\item [Case III.] two lines becoming inactive and no line becoming active. 
\end{description}

\noindent
{\bf Case I.} By Lemma \ref{lemma1} the lines becoming inactive are four among $r_1$, $r_3$, $r_4$, $r_6$ and $r_7$. Every choice of four of these lines includes at least two horizontal lines and thus the two lines becoming active include at most one horizontal line. Then in both cases (a) and (b) we get $p_2(\h)\leq p_2(\bar\h)-1\leq l_2^*-1$. Thus it is impossible to leave $\cB$.

\medskip
\noindent
{\bf Case II.} Three lines become inactive and one line becomes active. We distinguish the following cases:

\bi
\item[$\bullet$] If the line $r_1$ becomes inactive, by Lemma \ref{lemma0}{\it (i)} we know that no line can become active, so this case is not admissible.
\item[$\bullet$] If the line $r_1$ becomes active, by Lemma \ref{lemma0}{\it (ii)} we know that the only lines that can become inactive are $r_3$ and $r_4$. Since we require three deactivating lines, we deduce that this case is not admissible.
\item[$\bullet$] If the line $r_1$ does not become active nor inactive, by Lemma \ref{lemma1} the lines becoming inactive are three among $r_3$, $r_4$, $r_6$ and $r_7$ and the line becoming active is one among $r_2$, $r_3$, $r_4$ and $r_5$. 
\bi
\item [-)] If the line becoming active is $r_3$, the three lines becoming inactive are necessarily $r_4$, $r_6$ and $r_7$. Since $r_4$ becomes inactive the site $x_4$ is empty and $y_2$ is occupied; since $r_3$ becomes active in the site $x_3$ there is a free particle in $\bar\h$ (see Figure \ref{fig:caso-2}). Since the line $r_7$ becomes active, the site $y_3$ must be occupied. Furthermore, the site $s$ is empty, otherwise lines $r_4$ and $r_6$ can not become inactive. In this case we have one horizontal line becoming active and one horizontal becoming active, so we obtain
\be{5.2}
\left\{ \ba{lll}
p_1(\h)=p_1(\bar\h)-2, \\
p_2(\h)=p_2(\bar\h).
\ea\right.
\ee 

\setlength{\unitlength}{0.8pt}
\begin{figure}[h!]
	\begin{center}
		\begin{picture}(380,60)(-30,40)
		
		\thinlines 
		\qbezier[150](90,0)(190,0)(290,0)
		\qbezier[150](90,30)(190,30)(290,30) \qbezier[150](90,60)(190,60)(290,60)
		\qbezier[150](90,90)(190,90)(290,90)
		\qbezier[150](90,-30)(190,-30)(290,-30)
		\qbezier[150](90,-60)(190,-60)(290,-60)
		\qbezier[150](240,-80)(240,15)(240,110)
		\qbezier[150](270,-80)(270,15)(270,110)
		\qbezier[150](120,-80)(120,15)(120,110)
		\qbezier[150](150,-80)(150,15)(150,110)
		\qbezier[150](180,-80)(180,15)(180,110)
		\qbezier[150](210,-80)(210,15)(210,110)
		\put(95,15){\vector(1,0){25}}
		\put(75,12){$z_3$}
		\put(95,45){\vector(1,0){55}}
		\put(75,42){$z_1$}
		\put(135,-45){\vector(1,1){15}}
		\put(125,-52){$z_2$}
		\put(135,-15){\vector(1,1){15}}
		\put(125,-22){$y_3$}
		\put(165,75){\vector(1,-1){15}}
		\put(155,78){$y_1$}
		\put(165,-45){\vector(1,1){15}}
		\put(155,-52){$y_2$}
		\put(225,-45){\vector(-1,0){20}}
		\put(225,-48){$s$}
		\put(255,75){\vector(-1,-1){15}}
		\put(250,78){$x_3$}
		\put(255,-45){\vector(-1,1){15}}
		\put(250,-55){$x_4$}
		\put(195,15){\circle*{15}}
		\put(225,15){\circle{15}}
		\begin{Large}
		\put(245,10){$x_5$}
		\put(225,45){\circle*{15}}
		\put(225,75){\circle{15}}
		\put(255,45){\circle{15}}
		\put(225,-15){\circle{15}}
		\put(195,45){\circle{15}}
		\put(195,-15){\circle*{15}}
		\put(190,70){$t$}
		\put(195,-45){\circle{15}}
		\put(165,15){\circle*{15}}
		\put(135,15){\circle{15}}
		\put(165,-15){\circle{15}}
		\put(165,45){\circle{15}}
		
		\thinlines
		\put(225,-105){\vector(0,1){20}} \put(220,-118){$r_2$}
		\put(165,-105){\vector(0,1){20}} \put(160,-118){$r_7$}
		
		\put(250,-118){$r_5$} \put(255,-105){\vector(0,1){20}}
		\put(195,-105){\vector(0,1){20}} \put(190,-118){$r_6$}
		\put(315,-15){\vector(-1,0){20}} \put(318,-17){$r_4$}
		\put(315,15){\vector(-1,0){20}} \put(318,13){$r_1$}
		\put(315,45){\vector(-1,0){20}} \put(318,43){$r_3$}
		\end{Large}
		\end{picture}
		\vskip 4.5 cm
		\caption{Here we depict part of the configuration in the case II for $\D s=-2$.}
		\label{fig:caso-2}
	\end{center}
\end{figure}
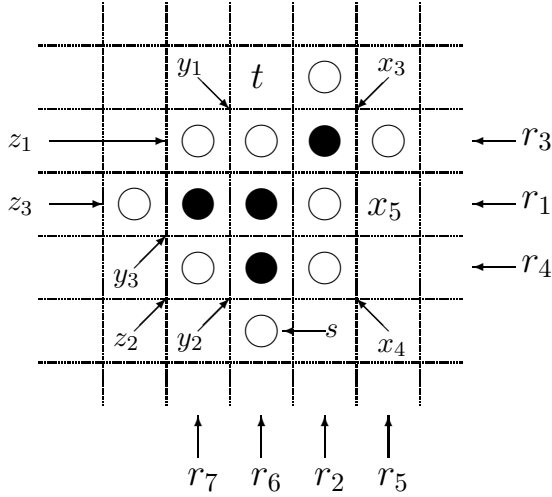

First, we consider the case (a). If $s(\h)=s^*-1$, by Lemma \ref{ulteriore}{\it (i)} we conclude that it is impossible to leave $\cB$. If $s(\h)\geq s^*$, since by (\ref{5.2}) $p_2(\h)=p_2(\bar\h)=l_2^*$, by Lemma \ref{ulteriore}{\it (ii)} we conclude that it is impossible to leave $\cB$.

In the case (b), by (\ref{5.2}), we get $p_2(\h)=p_2(\bar\h)\leq l_2^*-1$, thus it is impossible to leave $\cB$. 

\item [-)] If the line becoming active is $r_4$ we conclude by symmetry with a similar argument.

\item [-)] If $r_3$ and $r_4$ do not become active, in both cases (a) and (b) we have $p_2(\h)\leq p_2(\bar\h)-1\leq l_2^*-1$ and thus it is impossible to leave $\cB$.
\ei
\ei

\noindent
{\bf Case III.} There are two lines becoming inactive and no line becoming active. By Lemma \ref{lemma1} the lines that can become inactive are two among $r_1$, $r_3$, $r_4$, $r_6$ and $r_7$. If at least one horizontal line becomes inactive, in both cases (a) and (b) we get $p_2(\h)\leq p_2(\bar\h)-1\leq l_2^*-1$ and thus $\h\in\cB$. It remains to analyze the case in which the two lines that become inactive are both vertical: $r_6$ and $r_7$. Since $r_7$ must become inactive, it is necessary to have site $y_3$ occupied and sites $z_1$, $z_2$ and $z_3$ empty (see Figure \ref{fig:rappresentazione}).

First, we consider the case (a). We note that the sites $y_1$ and $t$ can not be both occupied, otherwise the line $r_6$ can not become inactive. For the same reason also the sites $y_2$ and $s$ can not be both occupied. Furthermore, if $y_1$ is occupied then $x_3$ is empty and similarly if $y_2$ is occupied then $x_4$ is empty. Thus, if either $\bar\h_{cl}$ is connected, or both $\bar\h_{cl}$ and $\h_{cl}$ are not connected, we get $v(\h)\geq v(\bar\h)-2$. Indeed, if the sites $y_1$ and $y_2$ are both empty, the particles in the sites $x_1$ and $y_3$ compose a cluster in $\bar\h$, thus with the move we get $v(\h)\geq v(\bar\h)$: in particular $v(\h)\geq v(\bar\h)-2$. If the sites $y_1$ and $y_2$ are both occupied we note that the sites $z_1$ and $z_2$ are the unique vacancies in $\bar\h$ that are not vacancies in $\h$, since the sites $t$, $x_3$, $s$ and $x_4$ are necessarily empty for what we have already observed. Thus we get $v(\h)\geq v(\bar\h)-2$. If only one site among $y_1$ and $y_2$ is occupied, then only one vacancy in $\bar\h$ is not a vacancy in $\h$, i.e., either $z_1$ if the site $y_1$ is occupied or $z_2$ if the site $y_2$ is occupied. Since $p_2(\h)=p_2(\bar\h)$, we get $\h\in\cB$, indeed

$$v(\h)\geq v(\bar\h)-2\geq p_{max}(\bar\h)-3=p_1(\h)-1=p_{max}(\h)-1.$$

\noindent
If $\bar\h_{cl}$ is not connected and $\h_{cl}$ is connected, since the moved particle and the particle in $y_3$ are free in $\h$, we deduce that at least one of the clusters in $\bar\h_{cl}$ must intersect $r_1$: $\bar\h_{cl}\cap\{r_1\setminus\{x_1,y_3\}\}\neq\emptyset$. This implies that $g_2'(\bar\h)\geq1$. 
By (a) and $s(\bar\h)\geq s^*+1$, we have that the circumscribed rectangle of $\bar\h$ is $\cR(2l_2^*+k,l_2^*)$ for any $k\geq0$. Since $k\geq0$, $0<\d<1$ and $\e\ll U_2$, we get $H(\bar\h)>\G$, indeed
\be{}
\ba{lll}
H(\bar\h)\geq H(\cR(2l_2^*+k,l_2^*))+\e(2l_2^*+k-1)+U_1=\\

=U_1l_2^*+2U_2l_2^*+kU_2-\e(2(l_2^*)^2+kl_2^*)+2\e l_2^*+\e k-\e +U_1>\G\Leftrightarrow \\

\Leftrightarrow \ \e(1+k(1-\d)-\d)>0.
\ea
\ee

In the case (b), since $p_2(\h)=p_2(\bar\h)$, we get $p_2(\h)\leq l_2^*-1$ and thus $\h\in\cB$.

\epr

\subsection{Proof of Proposition \lowercase{\ref{lemma3}}}
\label{dim3}

\bpr
By Lemma \ref{lemma3.2}{\it (i)} there are seven cases corresponding to $\D s=-1,0,1,2,3,4,5$. We analyze separately each case.

\medskip
\noindent
{\bf Case $\D s=-1$.} We note that in order to obtain $\D s=-1$ we have the four following possibilities:

\begin{description}
	\item [Case I.] four lines becoming inactive and three lines becoming active;
	\item [Case II.] three lines becoming inactive and two lines becoming active;
	\item [Case III.] two lines becoming inactive and one line becoming active;
	\item [Case IV.] one line becoming inactive and no line becoming active.
\end{description}

\noindent
{\bf Case I.} By Lemma \ref{lemma1} the lines $r_2$ and $r_5$ can not become inactive, thus the lines that become inactive are four among $r_1$, $r_3$, $r_4$, $r_6$ and $r_7$. Every choice of four of these lines includes at least two horizontal lines becoming inactive and then at most one horizontal line becoming active. Thus by Lemma \ref{decrescita} we get $p_2(\h)\leq p_2(\bar\h)-1\leq l_2^*-2\leq l_2^*-1$, so it is impossible to leave $\cB$.

\medskip
\noindent
{\bf Case II.} Again by Lemma \ref{lemma1} we know that lines $r_2$ and $r_5$ can not become inactive, thus the lines that become inactive are three among $r_1$, $r_3$, $r_4$, $r_6$ and $r_7$. We distinguish the following cases:

\bi
\item[$\bullet$] If the line $r_1$ becomes inactive, by Lemma \ref{lemma0}{\it (i)} no line can become active and thus this case is not admissible.
\item[$\bullet$] If the line $r_1$ becomes active, by Lemma \ref{lemma0}{\it (ii)} the only lines that can become inactive are $r_3$ and $r_4$. Since we require three deactivating lines, we deduce that this case is not admissible.
\item [$\bullet$] If the line $r_1$ does not become active nor inactive, by Lemma \ref{lemma1} the lines becoming inactive must be three among $r_3$, $r_4$, $r_6$ and $r_7$ and the lines becoming active must be two among $r_2$, $r_3$, $r_4$ and $r_5$. Every choice of three of these deactivating lines includes at least one horizontal line, so at most one horizontal line is becoming active. For this reason we get $p_2(\h)\leq p_2(\bar\h)\leq l_2^*-1$ and thus it is impossible to leave $\cB$.
\ei

\noindent
{\bf Case III.} There are two lines becoming inactive and one line becoming active. If the activating line is vertical, by Lemma \ref{decrescita} we get $p_2(\h)\leq p_2(\bar\h)\leq l_2^*-1$ and thus it is impossible to leave $\cB$, so we can assume for the rest of this case that the activating line is horizontal. If at least one of the deactivating lines is horizontal we conclude as before. Indeed $p_2(\h)\leq p_2(\bar\h)\leq l_2^*-1$, so it is impossible to leave $\cB$. Thus we can reduce the proof to the case in which the two deactivating lines are both vertical and the activating line is horizontal. Since by Lemma \ref{lemma1} the lines $r_2$ and $r_5$ can not become inactive, the lines that become inactive are necessarily $r_6$ and $r_7$. Since $r_7$ becomes inactive, the site $y_3$ must be occupied and the sites $z_1$, $z_2$ and $z_3$ must be empty (see Figure \ref{fig:caso-1}). Thus the particle in $x_1$ is not free in $\bar\h$ and the line $r_1$ cannot become active, so the activating line is either $r_3$ or $r_4$. By Lemma \ref{ult2}{\it (i)} we deduce that it is impossible to leave $\cB$.

\bigskip
\setlength{\unitlength}{0.8pt}
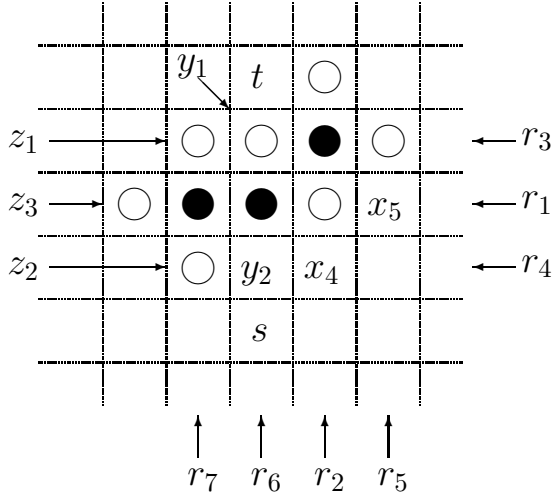
\begin{figure}
	\begin{center}
		\begin{picture}(380,60)(-30,40)
		
		\thinlines 
		\qbezier[150](90,0)(190,0)(290,0)
		\qbezier[150](90,30)(190,30)(290,30) \qbezier[150](90,60)(190,60)(290,60)
		\qbezier[150](90,90)(190,90)(290,90)
		\qbezier[150](90,-30)(190,-30)(290,-30)
		\qbezier[150](90,-60)(190,-60)(290,-60)
		\qbezier[150](240,-80)(240,15)(240,110)
		\qbezier[150](270,-80)(270,15)(270,110)
		\qbezier[150](120,-80)(120,15)(120,110)
		\qbezier[150](150,-80)(150,15)(150,110)
		\qbezier[150](180,-80)(180,15)(180,110)
		\qbezier[150](210,-80)(210,15)(210,110)
		\put(195,15){\circle*{15}}
		\put(225,15){\circle{15}}
		\begin{Large}
		\put(95,15){\vector(1,0){25}}
		\put(75,12){$z_3$}
		\put(95,45){\vector(1,0){55}}
		\put(75,42){$z_1$}
		\put(95,-15){\vector(1,0){55}}
		\put(75,-18){$z_2$}
		\put(165,75){\vector(1,-1){15}}
		\put(155,78){$y_1$}
		\put(245,10){$x_5$}
		\put(225,45){\circle*{15}}
		\put(225,75){\circle{15}}
		\put(255,45){\circle{15}}
		\put(215,-20){$x_4$}
		\put(195,45){\circle{15}}
		\put(185,-20){$y_2$}
		\put(190,70){$t$}
		\put(190,-50){$s$}
		\put(165,15){\circle*{15}}
		\put(135,15){\circle{15}}
		\put(165,-15){\circle{15}}
		\put(165,45){\circle{15}}
		
		\thinlines
		\put(225,-105){\vector(0,1){20}} \put(220,-118){$r_2$}
		\put(165,-105){\vector(0,1){20}} \put(160,-118){$r_7$}
		
		\put(250,-118){$r_5$} \put(255,-105){\vector(0,1){20}}
		\put(195,-105){\vector(0,1){20}} \put(190,-118){$r_6$}
		\put(315,-15){\vector(-1,0){20}} \put(318,-17){$r_4$}
		\put(315,15){\vector(-1,0){20}} \put(318,13){$r_1$}
		\put(315,45){\vector(-1,0){20}} \put(318,43){$r_3$}
		\end{Large}
		\end{picture}
		\vskip 4.5 cm
		\caption{Here we depict part of the configuration in the case III for $\D s=-1$.}
		\label{fig:caso-1}
	\end{center}
\end{figure}

\medskip
\noindent
{\bf Case IV.} We have one deactivating line and no activating line. Thus $p_2(\h)\leq p_2(\bar\h)-1\leq l_2^*-2\leq l_2^*-1$, so it is impossible to leave $\cB$.

\medskip
\noindent
{\bf Case $\D s=0$.} We observe that in order to obtain $\D s=0$ we have the four following possibilities:

\begin{description}
	\item [Case I.] three lines becoming inactive and three lines becoming active;
	\item [Case II.] two lines becoming inactive and two lines becoming active;
	\item [Case III.] one line becoming inactive and one line becoming active;
	\item [Case IV.] no line becoming inactive and no line becoming active.
\end{description}
\noindent
{\bf Case I.} Since by Lemma \ref{lemma1} lines $r_2$ and $r_5$ cannot become inactive, the lines becoming inactive are three among $r_1$, $r_3$, $r_4$, $r_6$ and $r_7$. Again we distinguish the following cases:
\bi
\item[$\bullet$] If the line $r_1$ becomes inactive, by Lemma \ref{lemma0}{\it (i)} no line can become active: this case is not admissible.
\item[$\bullet$] If the line $r_1$ becomes active, by Lemma \ref{lemma0}{\it (ii)} the only lines that can become inactive are $r_3$ and $r_4$. Since we require three deactivating lines, we deduce that this case is not admissible.
\item [$\bullet$] If the line $r_1$ does not become active nor inactive, then the lines becoming inactive are three among $r_3$, $r_4$, $r_6$ and $r_7$ and the lines becoming active are three among $r_2$, $r_3$, $r_4$ and $r_5$. Every choice of three of these deactivating lines includes at least one horizontal line, so at most one horizontal line is becoming active. Thus $p_2(\h)\leq p_2(\bar\h)\leq l_2^*-1$, so it is impossible to leave $\cB$.
\ei

\noindent
{\bf Case II.} We have two deactivating and two activating lines. If the two activating lines are vertical, then we get $p_2(\h)\leq p_2(\bar\h)\leq l_2^*-1$, so in this case it is impossible to leave $\cB$. Thus we can reduce the proof to the case in which at least one of the activating lines is horizontal. Again we consider the following cases:

\bi
\item [$\bullet$] If the line $r_1$ becomes inactive, again from Lemma \ref{lemma0}{\it (i)} we deduce that this case is not admissible, since no line can become active.
\item [$\bullet$] If the line $r_1$ becomes active, again from Lemma \ref{lemma0}{\it (ii)} we deduce that the two deactivating lines are $r_3$ and $r_4$. Thus we get $p_2(\h)=p_2(\bar\h)-1\leq l_2^*-1$: this implies that it is impossible to leave $\cB$.
\item [$\bullet$] If the line $r_1$ does not become active nor inactive, the deactivating lines are two among $r_3$, $r_4$, $r_6$, $r_7$ and the activating lines are two among $r_2$, $r_3$, $r_4$ and $r_5$. If at least one deactivating line is horizontal, we have that at most one activating line is horizontal, thus we get $p_2(\h)\leq p_2(\bar\h)\leq l_2^*-1$: this implies that $\h\in\cB$. Thus we analyze the case in which the deactivating lines are both vertical: $r_6$ and $r_7$. We focus on the case in which at least one activating line is horizontal, i.e., the line $r_3$ and/or $r_4$, otherwise $p_2(\h)=p_2(\bar\h)\leq l_2^*-1$ and thus it is impossible to leave $\cB$. If only one horizontal line becomes active, by Lemma \ref{ult2}{\it (i)} we conclude that $\h\in\cB$. In the case in which two horizontal lines become active, if $p_2(\bar\h)\leq l_2^*-2$ or $p_2(\bar\h)=l_2^*-1$ and $s(\h)=s^*-1$, then by Lemma \ref{ult2}{\it (ii)} we conclude that $\h\in\cB$. By Remark \ref{remark2} we know that $s(\bar\h)=s(\h)\geq s^*-1$. Hence we are left to consider the case $s(\bar\h)\geq s^*$. We obtain that the circumscribed rectangle of $\bar\h$ is $\cR(2l_2^*+k-1,l_2^*-1)$ for any $k\geq1$. Recalling that $\G=U_1l_2^*+2U_2l_2^*+U_1-U_2-2\e(l_2^*)^2+3\e l_2^*-2\e$, since $n(\bar\h)\geq2$, $k\geq1$, $\d<1$ and $\e\ll U_2$, we get $H(\bar\h)>\G$. Indeed
\be{}
H(\bar\h)\geq H(\cR(2l_2^*+k-1,l_2^*-1))+2\D>\G\Leftrightarrow \ 2U_2>\e(1+k(\d-1)).
\ee

\ei

\noindent
{\bf Case III.} We have one activating and one deactivating line. If the activating line is vertical we get $p_2(\h)\leq p_2(\bar\h)\leq l_2^*-1$ and we deduce that it is impossible to leave $\cB$. Now we assume that the activating line is horizontal. If also the deactivating line is horizontal, we get $p_2(\h)=p_2(\bar\h)\leq l_2^*-1$, thus it is impossible to leave $\cB$. We can reduce the proof to the case in which the activating line is horizontal ($r_1$ or $r_3$ or $r_4$) and the deactivating line is vertical ($r_6$ or $r_7$). We distinguish the following cases:

\bi
\item [$\bullet$] If the line $r_1$ becomes active, by Lemma \ref{lemma0}{\it (ii)} no vertical line can become inactive, thus this case is not admissible.

\item [$\bullet$] If the line $r_3$ ($r_4$ respectively) becomes active, by Lemma \ref{ult2}{\it (i)} we deduce that $\h\in\cB$, since by hypotheses $p_2(\bar\h)\leq l_2^*-1$.

\ei
\noindent
{\bf Case IV.} We have no deactivating and no activating line, thus we get $p_2(\h)=p_2(\bar\h)\leq l_2^*-1$, so it is impossible to leave $\cB$.

\medskip
\noindent
{\bf Case $\D s=1$.} We note that in order to obtain $\D s=1$ we have the four following possibilities:

\begin{description}
	\item [Case I.] one activating line and no deactivating line;
	\item [Case II.] two activating lines and one deactivating line;
	\item [Case III.] three activating lines and two deactivating lines;
	\item [Case IV.] four activating lines and three deactivating lines.
\end{description}
\noindent
{\bf Case I.} We have only to consider the case in which the activating line is horizontal, otherwise we get $p_2(\h)=p_2(\bar\h)\leq l_2^*-1$: thus $\h\in\cB$ and it is impossible to leave $\cB$. If the activating line is $r_i$, with $i\in{\{3,4\}}$, by Lemma \ref{ult2}{\it (i)} we deduce that it is impossible to leave $\cB$. Otherwise the horizontal line becoming active is $r_1$. We distinguish the following two cases:

\bi
\item [(a)] If the site $x_5$ is empty, by Lemma \ref{ulteriore}{\it (i),(ii)} we can conclude that it is impossible to leave $\cB$.
\item [(b)] If the site $x_5$ is occupied, in order that $r_1$ becomes active the particle in $x_5$ should be free in $\bar\h$. This implies that $n(\bar\h)\geq2$ ($x_1$ and $x_5$ contain a free particle). Since $p_2(\bar\h)=l_2^*-1$, by Remark \ref{remark2} we deduce that the circumscribed rectangle of $\bar\h$ is $\cR(2l_2^*+k-2,l_2^*-1)$ for any $k\geq0$. Since $k\geq0$, $\e\ll U_2$ and $0<\d<1$, we get $H(\bar\h)>\G$. Indeed
\be{}
\ba{lll}
H(\bar\h)\geq H(\cR(2l_2^*+k-2,l_2^*-1))+2\D=\\

=U_1l_2^*+2U_2l_2^*+kU_2-\e(2(l_2^*)^2-4l_2^*+kl_2^*-k+2)+U_1-2\e>\G\Leftrightarrow\\

\Leftrightarrow \ 2U_2>\e(2-\d+k(\d-1)).
\ea
\ee
\ei
\noindent
{\bf Case II.} We have two activating lines and one deactivating line. If no horizontal line becomes active, we get $p_2(\h)\leq p_2(\bar\h)\leq l_2^*-1$, thus $\h\in\cB$. Now we consider the case in which at least one horizontal line becomes active. In particular, the relevant cases are the followings:

\bi
\item [$\bullet$] An horizontal line becomes active and no horizontal line becomes inactive.

\bi
\item [-)] If $r_1$ is the horizontal line becoming active, by Lemma \ref{lemma0}{\it (ii)} no vertical line can be deactivated and thus this case is not admissible.
\item [-)] If $r_3$ (respectively $r_4$) is the horizontal line becoming active, by Lemma \ref{ult2}{\it (i)} we get $\h\in\cB$.

\ei

\item [$\bullet$] Two horizontal lines become active, thus $n(\bar\h)\geq2$. If $p_2(\bar\h)\leq l_2^*-2$ or $p_2(\bar\h)=l_2^*-1$ and $s(\bar\h)=s^*-2$ that implies $s(\h)=s^*-1$, we get $\h\in\cB$. By Remark \ref{remark2} we deduce that $s(\bar\h)\geq s^*-2$, so we are left to consider the case $s(\bar\h)\geq s^*-1$ that implies $s(\h)\geq s^*$. Thus the circumscribed rectangle of $\bar\h$ is $\cR(2l_2^*+k-2,l_2^*-1)$ for any $k\geq1$. Since $k\geq1$, $\d<1$ and $\e\ll U_2$, we get $H(\bar\h)>\G$. Indeed
\be{5.3}
\ba{lll}
H(\bar\h)\geq H(\cR(2l_2^*+k-2,l_2^*-1))+2\D=\\

=U_1l_2^*+2U_2l_2^*+kU_2-\e(2(l_2^*)^2+kl_2^*-k-4l_2^*+2)+U_1-2\e>\G\Leftrightarrow\\

\Leftrightarrow \ 2U_2>\e(2-\d+k(\d-1)).
\ea
\ee

\ei 

\noindent
{\bf Case III.} We have three activating lines and two deactivating lines. The lines becoming active are three among $r_1$, $r_2$, $r_3$, $r_4$, $r_5$ and the lines becoming inactive are two among $r_1$, $r_3$, $r_4$, $r_6$ and $r_7$. Again we distinguish the following cases:

\bi
\item [$\bullet$] If the line $r_1$ becomes active, by Lemma \ref{lemma0}{\it (ii)} the two deactivating lines are necessarily $r_3$ and $r_4$, that are both horizontal. Thus we get $p_2(\h)=p_2(\bar\h)-1\leq l_2^*-2\leq l_2^*-1$, so we conclude that it is impossible to leave $\cB$.
\item [$\bullet$] If the line $r_1$ becomes inactive, by Lemma \ref{lemma0}{\it (i)} no line can become active and thus this case is not admissible.
\item [$\bullet$] If the line $r_1$ does not become active nor inactive, the lines becoming active are three among $r_2$, $r_3$, $r_4$, $r_5$ and the lines becoming inactive are two among $r_3$, $r_4$, $r_6$ and $r_7$. If no horizontal line becomes active, we get $p_2(\h)=p_2(\bar\h)\leq l_2^*-1$ and thus it is impossible to leave $\cB$, otherwise we consider the following cases:

\bi
\item [-)] Two horizontal lines become active, i.e., the lines $r_3$ and $r_4$ become active, thus $n(\bar\h)\geq2$ (in the sites $x_3$ and $x_4$ there are free particles). If $p_2(\bar\h)\leq l_2^*-2$ or $p_2(\bar\h)=l_2^*-1$ and $s(\bar\h)=s^*-2$ that implies $s(\h)=s^*-1$, by Lemma \ref{ult2}{\it (ii)} we get $\h\in\cB$. By Remark \ref{remark2} we deduce that $s(\bar\h)\geq s^*-2$, so we are left to consider the case $s(\bar\h)\geq s^*-1$ that implies $s(\h)\geq s^*$. Thus by (\ref{5.3}) we get $H(\bar\h)\geq H(\cR(2l_2^*+k-2,l_2^*-1))+2\D>\G$ for any $k\geq1$.

\item [-)] Only one horizontal line becomes active, that is either $r_3$ or $r_4$. By Lemma \ref{ult2}{\it (i)}, we get $\h\in\cB$.
\ei
\ei

\noindent
{\bf Case IV.} We have four activating lines and three deactivating lines. The lines becoming active are four among $r_1$, $r_2$, $r_3$, $r_4$, $r_5$ and the lines becoming inactive are three among $r_1$, $r_3$, $r_4$, $r_6$ and $r_7$. Again we distinguish the following cases:

\bi
\item [$\bullet$] If the line $r_1$ becomes active, by Lemma \ref{lemma0}{\it (ii)} at most two lines ($r_3$ and/or $r_4$) can become inactive, thus this case is not admissible.
\item [$\bullet$] If the line $r_1$ becomes inactive, by Lemma \ref{lemma0}{\it (i)} no line becomes active, thus this case is not admissible.
\item [$\bullet$] If the line $r_1$ does not become active nor inactive, the lines becoming active are $r_2$, $r_3$, $r_4$, $r_5$ and lines becoming inactive are only $r_6$ and $r_7$: this means that $\D s=2$, which is in contradiction with $\D s=1$.
\ei

\medskip
\noindent
{\bf Case $\D s=2$.} We observe that in order to obtain $\D s=2$ we have the three following possibilities:

\begin{description}
	\item [Case I.] four activating lines and two deactivating lines;
	\item [Case II.] three activating lines and one deactivating line;
	\item [Case III.] two activating lines and no deactivating line.
\end{description}
\noindent
{\bf Case I.} By Lemma \ref{lemma1} the lines $r_6$, $r_7$ cannot become active and the lines $r_2$, $r_5$ cannot become inactive. Thus the activating lines are four among $r_1$, $r_2$, $r_3$, $r_4$, $r_5$ and the deactivating lines are two among $r_1$, $r_3$, $r_4$, $r_6$ and $r_7$. We distinguish the following cases:

\bi
\item [$\bullet$] If the line $r_1$ becomes active, by Lemma \ref{lemma0}{\it (ii)} the lines $r_3$ and $r_4$ become inactive, that are both horizontal. Thus we get $p_2(\h)=p_2(\bar\h)-1\leq l_2^*-2\leq l_2^*-1$, so it is impossible to leave $\cB$.
\item [$\bullet$] If the line $r_1$ becomes inactive, by Lemma \ref{lemma0}{\it (i)} no line can become active, thus this case is not admissible.
\item [$\bullet$] If the line $r_1$ does not become active nor inactive, then the four activating lines are $r_2$, $r_3$, $r_4$, $r_5$ and the two deactivating lines are $r_6$ and $r_7$.
By Remark \ref{remark2} we deduce that $s(\bar\h)\geq s^*-3$ that implies $s(\h)\geq s^*-1$. If $p_2(\bar\h)\leq l_2^*-2$ or $p_2(\bar\h)=l_2^*-1$ and $s(\h)=s^*-1$, by Lemma \ref{ult2}{\it (ii)} we deduce that it is impossible to leave $\cB$. If $s(\h)\geq s^*$, that implies $s(\bar\h)\geq s^*-2$, the circumscribed rectangle of $\bar\h$ is $\cR(2l_2^*+k-3,l_2^*-1)$ for any $k\geq1$. Since $k\geq1$, $\d<1$ and $\e\ll U_2$, we get $H(\bar\h)>\G$. Indeed
\be{5.4}
\ba{lll}
H(\bar\h)\geq H(\cR(2l_2^*+k-3,l_2^*-1))+2\D=\\

=U_1l_2^*+2U_2l_2^*+kU_2-\e(2(l_2^*)^2+kl_2^*-k-5\e l_2^*+3)+U_1-U_2-2\e>\G\Leftrightarrow\\
\Leftrightarrow \ 2U_2>\e(3-2\d+k(\d-1)). 
\ea
\ee
\ei 
\noindent
{\bf Case II.} We have three activating lines and one deactivating line. Again we distinguish the following cases:

\bi
\item [$\bullet$] If the line $r_1$ becomes active, by Lemma \ref{lemma0}{\it (ii)} we know that either $r_3$ or $r_4$ becomes inactive, that are both horizontal. If $p_2(\bar\h)\leq l_2^*-2$, we get $p_2(\h)\leq p_2(\bar\h)+1\leq l_2^*-1$, thus it is impossible to leave $\cB$. In the case $p_2(\bar\h)=l_2^*-1$ that implies $p_2(\h)=l_2^*-1$, by Lemma \ref{ulteriore}{\it (i),(ii)} we can conclude that it is impossible to leave $\cB$.
\item [$\bullet$] If the line $r_1$ becomes inactive, by Lemma \ref{lemma0}{\it (i)} no line can become active and thus this case is not admissible.
\item [$\bullet$] If the line $r_1$ does not become active nor inactive, thus the activating lines are three among $r_2$, $r_3$, $r_4$, $r_5$ and the deactivating line is one among $r_3$, $r_4$, $r_6$ and $r_7$. Every choice of three of these activating lines includes at least one horizontal line. Thus the relevant cases are the followings:

\bi
\item [-)] If one horizontal line becomes active, i.e., either $r_3$ or $r_4$, by Lemma \ref{ult2}{\it (i)} we deduce that it is impossible to leave $\cB$. 
\item [-)] If two horizontal lines become active, we get $n(\bar\h)\geq2$. By Remark \ref{remark2} we deduce that $s(\bar\h)\geq s^*-3$ that implies $s(\h)\geq s^*-1$. If $p_2(\bar\h)\leq l_2^*-2$ or $p_2(\bar\h)=l_2^*-1$ and $s(\h)=s^*-1$, by Lemma \ref{ult2}{\it (ii)} we deduce that it is impossible to leave $\cB$. If $s(\h)\geq s^*$, that implies $s(\bar\h)\geq s^*-2$, by (\ref{5.4}) we get $H(\bar\h)\geq H(\cR(2l_2^*+k-3,l_2^*-1))+2\D>\G$ for any $k\geq1$.
\ei 
\ei
\noindent
{\bf Case III.} We have two activating lines and no deactivating line. If no horizontal line becomes active, we get $p_2(\h)=p_2(\bar\h)\leq l_2^*-1$ and thus it is impossible to leave $\cB$. Thus we can reduce our analysis to the case in which at least one horizontal line becomes active. We distinguish the following cases:

\bi
\item [$\bullet$] If the line $r_1$ becomes active, 
we distinguish the following cases:
\bi
\item [-)] If two horizontal lines become active, i.e., the line $r_1$ and one among $r_3$ and $r_4$, by Remark \ref{remark2}, we deduce that $s(\bar\h)\geq s^*-3$ that implies $s(\h)\geq s^*-1$. If $p_2(\bar\h)\leq l_2^*-2$ or $p_2(\bar\h)=l_2^*-1$ and $s(\h)=s^*-1$, by Lemma \ref{ult2}{\it (ii)} we deduce that it is impossible to leave $\cB$. If $s(\h)\geq s^*$, that implies $s(\bar\h)\geq s^*-2$, by (\ref{5.4}) we get $H(\bar\h)\geq H(\cR(2l_2^*+k-3,l_2^*-1))+2\D>\G$ for any $k\geq1$.
\item [-)] If only one horizontal line become active, i.e., the line $r_1$, we get $p_2(\h)=p_2(\bar\h)+1$. If $p_2(\bar\h)\leq l_2^*-2$, we get $p_2(\h)\leq l_2^*-1$ and thus it is impossible to leave $\cB$. Suppose now that the site $x_5$ is empty. If $p_2(\bar\h)=l_2^*-1$ then $p_2(\h)=l_2^*$, so by Lemma \ref{ulteriore}{\it (i),(ii)} we conclude that it is impossible to leave $\cB$. Hence we are left to consider the case $x_5$ occupied. Thus $x_5$ must contain a free particle, otherwise the line $r_1$ is already active in $\bar\h$. If $p_2(\bar\h)=l_2^*-1$ then $p_2(\h)=l_2^*$, so by (\ref{5.4}) we get $H(\bar\h)\geq H(\cR(2l_2^*+k-3,l_2^*-1))+2\D>\G$ for any $k\geq1$.
\ei
\item [$\bullet$] If the line $r_1$ does not become active, the horizontal lines that can become active are $r_3$ and $r_4$. We distinguish the following cases:
\bi
\item [-)] If only one horizontal line become active, i.e., either $r_3$ or $r_4$, by Lemma \ref{ult2}{\it (i)} we get $\h\in\cB$.

\item [-)] If two horizontal lines become active, i.e., the lines $r_3$ and $r_4$, by Remark \ref{remark2} we deduce that $s(\bar\h)\geq s^*-3$ that implies $s(\h)\geq s^*-1$. If $p_2(\bar\h)\leq l_2^*-2$ or $p_2(\bar\h)=l_2^*-1$ and $s(\h)=s^*-1$, by Lemma \ref{ult2}{\it (ii)} we deduce that it is impossible to leave $\cB$. If $s(\h)\geq s^*$, that implies $s(\bar\h)\geq s^*-2$, by (\ref{5.4}) we get $H(\bar\h)\geq H(\cR(2l_2^*+k-3,l_2^*-1))+2\D>\G$ for any $k\geq1$.
\ei
\ei

\medskip
\noindent
{\bf Case $\D s=3$.} We note that in order to obtain $\D s=3$ we have the three following possibilities:

\begin{description}
	\item [Case I.] five activating lines and two deactivating lines;
	\item [Case II.] four activating lines and one deactivating line;
	\item [Case III.] three activating lines and no deactivating line.
\end{description}
\noindent
{\bf Case I.} Since by Lemma \ref{lemma1} the lines $r_6$ and $r_7$ cannot become active, thus the lines becoming active are $r_1$, $r_2$, $r_3$, $r_4$, $r_5$ and the lines becoming inactive are $r_6$ and $r_7$. Since $r_1$ becomes active, by Lemma \ref{lemma0}{\it (ii)} the only lines that can become inactive are $r_3$ and $r_4$, which is in contradiction with this case. 

\medskip
\noindent
{\bf Case II.} We have four activating lines and one deactivating line. Since by Lemma \ref{lemma1} the lines $r_6$ and $r_7$ cannot become active, thus we have that at least one horizontal line must become active. We distinguish the following cases:

\bi
\item [$\bullet$] If the line $r_1$ becomes active, by Lemma \ref{lemma0}{\it (ii)} the deactivating line is one among $r_3$ and $r_4$, that are both horizontal. Thus we deduce that $p_2(\h)\leq p_2(\bar\h)+1$. If only one horizontal line becomes active we get $p_2(\h)=p_2(\bar\h)\leq l_2^*-1$ and thus it is impossible to leave $\cB$. If two horizontal lines become active, i.e., $r_1$ and one line among $r_3$ and $r_4$, we get $p_2(\h)=p_2(\bar\h)+1$. If $p_2(\bar\h)\leq l_2^*-2$, we get $p_2(\h)\leq l_2^*-1$ and thus it is impossible to leave $\cB$. If $p_2(\bar\h)=l_2^*-1$ then $p_2(\h)=l_2^*$, so by Lemma \ref{ulteriore}{\it (i),(ii)} we conclude that $\h\in\cB$.
\item [$\bullet$] If the line $r_1$ becomes inactive, by Lemma \ref{lemma0}{\it (i)} no line can become active and thus this case is not admissible.
\item [$\bullet$] If the line $r_1$ does not become active nor inactive, we deduce that one horizontal line becomes active and one horizontal becomes inactive. Thus $p_2(\h)=p_2(\bar\h)\leq l_2^*-1$, so it is impossible to leave $\cB$.
\ei

\noindent
{\bf Case III.} We have three activating lines and no deactivating line. Since by Lemma \ref{lemma1} the lines $r_6$ and $r_7$ can not become active, we have that at least one horizontal line must become active. We distinguish the following cases:

\bi
\item [$\bullet$] If the line $r_1$ becomes active, we distinguish the following cases:

\bi
\item [-)] If only one horizontal line becomes active, we get $p_2(\h)=p_2(\bar\h)+1$. If $p_2(\bar\h)\leq l_2^*-2$, we get $p_2(\h)\leq l_2^*-1$, thus it is impossible to leave $\cB$. Suppose now that the site $x_5$ is empty. If $p_2(\bar\h)=l_2^*-1$ then $p_2(\h)=l_2^*$, so by Lemma \ref{ulteriore}{\it (i),(ii)} we conclude that it is impossible to leave $\cB$. Hence we are left to consider the case $x_5$ occupied. Thus it contains a free particle in $\bar\h$, otherwise the line $r_1$ is already active in $\bar\h$. If $p_2(\bar\h)=l_2^*-1$ then $p_2(\h)=l_2^*$, so by Remark \ref{remark2} we deduce that the circumscribed rectangle of $\bar\h$ is $\cR(2l_2^*+k-5,l_2^*-1)$ for any $k\geq1$. Since $n(\bar\h)\geq2$, $k\geq1$, $\d<1$ and $\e\ll U_2$, we get $H(\bar\h)>\G$. Indeed
\be{}
\ba{lll}
H(\bar\h)\geq H(\cR(2l_2^*+k-5,l_2^*-1))+2\D=\\

=U_1l_2^*+2U_2l_2^*+U_1+kU_2-3U_2-\e(2(l_2^*)^2+kl_2^*-k-7l_2^*+7)>\G\Leftrightarrow\\

\Leftrightarrow \ 2U_2>\e(5-4\d+k(\d-1)).
\ea
\ee

\item [-)] If two horizontal lines become active, i.e., the line $r_1$ and one among $r_3$ and $r_4$, by Remark \ref{remark2} we deduce that $s(\bar\h)\geq s^*-4$ that implies $s(\h)\geq s^*-1$. If $p_2(\bar\h)\leq l_2^*-2$ or $p_2(\bar\h)=l_2^*-1$ and $s(\h)=s^*-1$, by Lemma \ref{ult2}{\it (ii)} we deduce that it is impossible to leave $\cB$. If $s(\h)\geq s^*$, that implies $s(\bar\h)\geq s^*-3$, we obtain that the circumscribed rectangle of $\bar\h$ is $\cR(2l_2^*+k-4, l_2^*-1)$ for any $k\geq1$. Since $n(\bar\h)\geq2$, $k\geq1$, $\d<1$ and $\e\ll U_2$, we get $H(\bar\h)>\G$. Indeed
\be{5.5}
H(\bar\h)\geq H(\cR(2l_2^*+k-4,l_2^*-1))+2\D>\G\Leftrightarrow \ 2U_2>\e(4-3\d+k(\d-1)). 	
\ee

\item [-)] If three horizontal lines become active, we get $p_2(\h)=p_2(\bar\h)+3$. If $p_2(\bar\h)\leq l_2^*-4$, we get $p_2(\h)\leq l_2^*-1$, thus it is impossible to leave $\cB$. If $p_2(\bar\h)=l_2^*-3$ then $p_2(\h)=l_2^*$, so by Lemma \ref{ulteriore}{\it (i),(ii)} we conclude that $\h\in\cB$. By Remark \ref{remark2} we know that $s(\bar\h)\geq s^*-4$. If $s(\bar\h)=s^*-4$ then $s(\h)=s^*-1$, thus by Lemma \ref{ulteriore}{\it (i)} we conclude that it is impossible to leave $\cB$. If $s(\bar\h)\geq s^*-3$, with $p_2(\bar\h)=l_2^*-2$ that implies $p_2(\h)=l_2^*+1$, we get $H(\bar\h)\geq H(\cR(2l_2^*+k-3,l_2^*-2))+3\D>\G$ for any $k\geq1$, by a direct computation, since $n(\bar\h)\geq3$. We refer to the Appendix for the explicit computation. 
\ei

\item [$\bullet$] If the line $r_1$ does not become active, the horizontal lines that can become active are $r_3$ and $r_4$. Again we distinguish the following cases:
\bi
\item [-)] If only one horizontal line becomes active, by Lemma \ref{ult2}{\it (i)}, we get $\h\in\cB$.
\item [-)] If two horizontal lines become active, i.e., the lines $r_3$ and $r_4$, we get $p_2(\h)=p_2(\bar\h)+2$. By Remark \ref{remark2} we know that $s(\bar\h)\geq s^*-4$, that implies $s(\h)\geq s^*-1$. If $p_2(\bar\h)\leq l_2^*-2$ or $p_2(\bar\h)=l_2^*-1$ and $s(\h)=s^*-1$, by Lemma \ref{ult2}{\it (ii)} we obtain that it is impossible to leave $\cB$. If $s(\bar\h)\geq s^*-3$, that implies $s(\h)\geq s^*$, with $p_2(\bar\h)=l_2^*-1$ and then $p_2(\h)=l_2^*+1$, by (\ref{5.5}) we get $H(\bar\h)\geq H(\cR(2l_2^*+k-4,l_2^*-1))+2\D>\G$ for any $k\geq1$, since $n(\bar\h)\geq2$.
\ei
\ei

\medskip
\noindent
{\bf Case $\D s=4$.}  We observe that in order to obtain $\D s=4$ we have the two following possibilities:

\begin{description}
	\item [Case I.] five activating lines and one deactivating line;
	\item [Case II.] four activating lines and no deactivating line.
\end{description}

\noindent
{\bf Case I.} The five activating lines are $r_1$, $r_2$, $r_3$, $r_4$ and $r_5$, thus we get $p_2(\h)=p_2(\bar\h)+3$ and $n(\bar\h)\geq3$. If $p_2(\bar\h)\leq l_2^*-4$, we get $p_2(\h)\leq l_2^*-1$ and thus $\h\in\cB$. If $p_2(\bar\h)=l_2^*-3$ then $p_2(\h)=l_2^*$, so by Lemma \ref{ulteriore}{\it (i),(ii)} we deduce that it is impossible to leave $\cB$. By Remark \ref{remark2} we know that $s(\bar\h)\geq s^*-5$. If $s(\bar\h)=s^*-5$ then $s(\h)=s^*-1$, thus by Lemma \ref{ulteriore}{\it (i)} we conclude that it is impossible to leave $\cB$. Hence we are left to the case $s(\bar\h)\geq s^*-4$, with $p_2(\bar\h)=l_2^*-1$ and $p_2(\bar\h)=l_2^*-2$. In both cases, we directly get $H(\bar\h)>\G$ for any $k\geq1$. We refer to the Appendix for the explicit computations. 

\medskip
\noindent
{\bf Case II.} We have four activating lines and no deactivating line. We observe that there are at least two horizontal lines becoming active, so at least one line among $r_3$ and $r_4$ must become active. If $p_2(\bar\h)\leq l_2^*-3$, we get $p_2(\h)\leq p_2(\bar\h)+2\leq l_2^*-1$ and thus it is impossible to leave $\cB$. If $l_2^*-2\leq p_2(\bar\h)\leq l_2^*-1$ and $p_2(\h)=l_2^*$, by Lemma \ref{ulteriore}{\it (i),(ii)} we conclude that it is impossible to leave $\cB$. By Remark \ref{remark2} we know that $s(\bar\h)\geq s^*-5$. Thus we analyze separately the case $s(\bar\h)=s^*-5$ and $s(\bar\h)\geq s^*-4$. If $s(\bar\h)=s^*-5$ then $s(\h)=s^*-1$, so by Lemma \ref{ulteriore}{\it (i)} we deduce that it is impossible to leave $\cB$. If $s(\bar\h)\geq s^*-4$, we have to analyze the cases $p_2(\bar\h)=l_2^*-2$ and $p_2(\bar\h)=l_2^*-1$. Thus we get $H(\bar\h)>\G$ with a direct computation. We refer to the Appendix for the explicit computations.

\medskip
\noindent
{\bf Case $\D s=5$.} We note that the unique way to obtain $\D s=5$ is with five activating lines and no deactivating line. By Lemma \ref{lemma1} the lines becoming active are necessarily $r_1$, $r_2$, $r_3$, $r_4$ and $r_5$: three are horizontal and two are vertical. Thus we get $p_2(\h)=p_2(\bar\h)+3$. If $p_2(\bar\h)\leq l_2^*-4$, we get $p_2(\h)\leq l_2^*-1$, so it is impossible to leave $\cB$. If $p_2(\bar\h)=l_2^*-3$ then $p_2(\h)=l_2^*$, so by Lemma \ref{ulteriore}{\it (i),(ii)} we can conclude that it is impossible to leave $\cB$. For the remaining cases, by Remark \ref{remark2} we know that $s(\bar\h)\geq s^*-6$, that implies $s(\h)\geq s^*-1$. We analyze separately the cases $s(\bar\h)=s^*-6$ and $s(\bar\h)\geq s^*-5$. If $s(\bar\h)=s^*-6$ then $s(\h)=s^*-1$, thus by Lemma \ref{ulteriore}{\it (i)} we can conclude that $\h\in\cB$. If $s(\bar\h)\geq s^*-5$ with a direct computation we get $H(\bar\h)>\G$, since $n(\bar\h)\geq4$. We refer to the Appendix for the explicit computations. 
\epr

\subsection{Proof of Proposition \lowercase{\ref{lemma4}}}
\label{dim4}

\bpr
Let $(\bar\h,\h)$ be a move with $\bar\h\in\cB$ and $\D s$ be its corresponding variation of $s$; by Lemma \ref{lemma3.2}{\it (i)} we have to consider only the cases corresponding to $\D s=3,4,5$.

\bigskip
\noindent
{\bf Case $\D s=3$.} Let $\D s=3$; by Remark \ref{remark2}{\it (i)} we consider only the case $s(\bar\h)\geq s^*-4$. For the convenience of the proof, we distinguish the following cases:

\bi
\item [(a)] $s(\bar\h)=s^*-4$;
\item [(b)] $s^*-3\leq s(\bar\h)\leq s^*-1$;
\item [(c)] $s(\bar\h)\geq s^*$.
\ei

First we consider the case (a). If $s(\bar\h)=s^*-4$, then $s(\h)=s^*-1$. If $p_2(\h)\leq l_2^*-1$, we directly get $\h\in\cB$. If $p_2(\h)\geq l_2^*$, by Lemma \ref{lemma3.2}{\it (ii)} we obtain

$$v(\h)\geq3p_{min}(\h)-6\geq p_{min}(\h)-1\Leftrightarrow \ 2p_{min}(\h)\geq5.$$

\noindent
Since $p_{min}(\h)\geq4$, by (\ref{defB}) we deduce that in the case (a) it is impossible to leave $\cB$.

By Lemma \ref{lemma3.2}{\it (ii)} we have $n(\bar{\h})\geq2$. The circumscribed rectangle of $\bar\h$ is $\cR(2l_2^*-k-x,l_2^*+k)$ for any $k\geq0$, where in the case (b) $2\leq x\leq4$, and in the case (c) $x \leq 1$. In the case (b) we obtain
\be{}
\ba{ll}
H(\bar\h)\geq H(\cR(2l_2^*-k-x,l_2^*+k))+2\D>\G\Leftrightarrow \\

\Leftrightarrow \ \e k^2+k[U_1-U_2-\e l_2^*+x\e]+U_1+3U_2-xU_2+x\e l_2^*-3\e l_2^*>0.
\ea
\ee
\noindent
Since $k\geq0$, $x\geq2$, $\d<1$ and $\e\ll U_2$ we obtain $H(\bar{\h})>\G$, indeed

\be{5.10}
\ba{ll}
U_1-U_2-\e l_2^*+x\e U_1-2U_2-\e\d+x\e>\e(x-\d)\geq \e(2-\d)>\e>0 \\

U_1+3U_2-xU_2+x\e l_2^*-3\e l_2^*=U_1+x\e\d-3\e\d\geq U_1-\e\d\gg U_1-U_2>0.
\ea
\ee

In the case (c), since $s(\bar\h)\geq s^*$, by (\ref{defB}) we deduce that $p_2(\bar\h)=l_2^*$, so in this case $k=0$ and $v(\bar\h)\geq p_{max}(\bar\h)-1$. Since $x\leq1$, $1-\d>0$ and $U_2\gg\e$, we get $H(\bar\h)>\G$. Indeed
\be{}
\ba{ll}
H(\bar\h)\geq H(\cR(2l_2^*-x,l_2^*))+\e(2l_2^*-x-1)+2\D>\G\Leftrightarrow\\

\Leftrightarrow \ 2U_2+U_1>\e(1+x(1-\d)).
\ea
\ee

This concludes the proof in the case $\D s=3$.

\bigskip
\noindent
{\bf Case $\D s=4$.} Let $\D s=4$; by Remark \ref{remark2} we consider only the case $s(\bar\h)\geq s^*-5$. For the convenience of the proof, we distinguish the following cases:

\bi
\item [(a)] $s(\bar\h)=s^*-5$;
\item [(b)] $s^*-4\leq s(\bar\h)\leq s^*-1$;
\item [(c)] $s(\bar\h)\geq s^*$.
\ei

First, we consider the case (a). If $s(\bar\h)=s^*-5$, we get $s(\h)=s^*-1$. If $p_2(\h)\leq l_2^*-1$, we directly get $\h\in\cB$. If $p_2(\h)\geq l_2^*$, by Lemma \ref{lemma3.2}(ii) we obtain

$$v(\h)\geq4p_{min}(\h)-7\geq p_{min}(\h)-1\Leftrightarrow \ p_{min}(\h)\geq2.$$

\noindent
Since $p_{min}(\h)\geq4$, by (\ref{defB}) we obtain that in the case (a) it is impossible to leave $\cB$.

By Lemma \ref{lemma3.2}(ii) we have $n(\bar{\h})\geq3$. The circumscribed rectangle of $\bar\h$ is $\cR(2l_2^*-k-x,l_2^*+k)$ for any $k\geq0$, where in the case (b) $2\leq x\leq5$, and in the case (c) $x\leq1$. In the case (b), the proof is analogue to the one for $\D s=3$, so we get $H(\bar\h)\geq H(\cR(2l_2^*-k-x,l_2^*+k))+3\D>\G$ and we refer to the Appendix for the explicit computation.

In the case (c), since $s(\bar\h)\geq s^*$, by (\ref{defB}) we deduce that $p_2(\bar\h)=l_2^*$, so in this case $k=0$ and $v(\bar\h)\geq p_{max}(\bar\h)-1$. Thus we get $H(\bar\h)\geq H(\cR(2l_2^*-x,l_2^*))+\e(2l_2^*-x-1)+3\D>\G$ and we refer to the Appendix for the explicit computation. This concludes the proof in the case $\D s=4$.

\bigskip
\noindent
{\bf Case $\D s=5$.} Let $\D s=5$; by Remark \ref{remark2} we consider only the case $s(\bar\h)\geq s^*-6$. For the convenience of the proof, we distinguish the following cases:

\bi
\item [(a)] $s(\bar\h)=s^*-6$;
\item [(b)] $s^*-5\leq s(\bar\h)\leq s^*-1$;
\item [(c)] $s(\bar\h)\geq s^*$.
\ei

First we consider case (a). If $s(\bar\h)=s^*-6$, then $s(\h)=s^*-1$. If $p_2(\h)\leq l_2^*-1$, we get $\h\in\cB$. If $p_2(\h)\geq l_2^*$, by Lemma \ref{lemma3.2}(ii) we obtain

$$v(\h)\geq5p_{min}(\h)-8\geq p_{min}(\h)-1\Leftrightarrow \ 4p_{min}(\h)\geq7.$$

\noindent
Since $p_{min}(\h)\geq4$, by (\ref{defB}) we deduce that in the case (a) it is impossible to leave $\cB$.

By Lemma \ref{lemma3.2}(ii) we have $n(\bar{\h})\geq4$. The circumscribed rectangle of $\bar\h$ is $\cR(2l_2^*-k-x,l_2^*+k)$ for any $k\geq0$, where in the case (b) $2\leq x\leq6$, and in the case (c) $x\leq1$. In the case (b), the proof is analogue to the one for $\D s=3$, so we get $H(\bar\h)\geq H(\cR(2l_2^*-k-x,l_2^*+k))+4\D>\G$ and we refer to the Appendix for the explicit computation.

In the case (c), since $s(\bar\h)\geq s^*$, by definition (\ref{defB}) we deduce that $p_2(\bar\h)=l_2^*$, so in this case $k=0$, and $v(\bar\h)\geq p_{max}(\bar\h)-1$. Thus we get $H(\bar\h)\geq H(\cR(2l_2^*-x,l_2^*))+\e(2l_2^*-x-1)+4\D>\G$ and we refer to the Appendix for the explicit computation. This concludes the proof in the case $\D s=5$.

\epr

\subsection{Proof of Proposition \lowercase{\ref{lemma5}}}
\label{dim5}
\bpr
We will prove that if $\D s=-1$ and $(\bar{\h},\h)\in\partial\cB$, then $\max{\{H(\bar{\h}),H(\h)\}}\geq\G$. Moreover, we will identify in which configurations the maximum is equal to $\G$ to prove Theorem \ref{3.1}{\it (i)}. By Remark \ref{remark2} we can only consider the case $s(\bar\h)\geq s^*$, i.e., $s(\bar\h)=s^*+k$ for any $k\geq0$. By (\ref{defB}) and $s(\bar\h)\geq s^*$, we have 

\be{sistemap2} \left\{
\ba{lll}
v(\bar{\h})&\geq& p_{max}(\bar{\h})-1 \\           
p_2(\bar{\h})&=& l_2^*.
\ea
\right.
\ee
Thus we obtain that the circumscribed rectangle of $\bar{\h}$ is $\cR(2l_2^*+k-1,l_2^*)$, indeed $p_1(\bar{\h})=s(\bar\h)-l_2^*=2l_2^*+k-1$. From (\ref{energeta}) we have
\be{}
H(\bar{\h})=H(\cR(2l_2^*+k-1,l_2^*))+\e v(\bar{\h})+U_1g_2'(\bar\h)+U_2g_1'(\bar{\h})+\D n(\bar{\h}). \ee

\noindent
For $p_{min}(\bar\h)\geq 4$ we consider the following cases:

\begin{description}
	\item [Case I.] $n(\bar{\h})\geq 1;$
	\item [Case II.] $g_1'(\bar{\h})+g_2'(\bar{\h})\geq1$ (and $n(\bar{\h})=0$);
	\begin{itemize}
		\item [A.] $g_2'(\bar{\h})=1;$
		\item [B.] $g_1'(\bar{\h})=1;$
		\item [C.] either $g_1'(\bar{\h})=1$ and $g_2'(\bar{\h})=1$, or $g_2'(\bar\h)\geq2$;
		\item [D.] $g_1'(\bar{\h})\geq 2$ and $g_2'(\bar{\h})=0$.
	\end{itemize}
	\item [Case III.] $n(\bar{\h})=0$ and $g_1'(\bar{\h})+g_2'(\bar{\h})=0.$
\end{description}

\setlength{\unitlength}{1.1pt}
\begin{figure}
	\begin{picture}(400,90)(0,30)
	\thinlines
	\put(20,50){\line(1,0){50}}
	\put(20,80){\line(1,0){50}}
	\put(20,50){\line(0,1){30}}
	\put(70,50){\line(0,1){30}}
	\thinlines  
	\put(40,60){$(1)$}
	\begin{Huge}
	\put(75,60){$\longrightarrow$} \end{Huge}
	\begin{footnotesize} \put(85,55){$\D$} \end{footnotesize}
	\thinlines
	\put(115,50){\line(1,0){50}}
	\put(115,80){\line(1,0){50}}
	\put(115,50){\line(0,1){30}}
	\put(165,50){\line(0,1){30}}
	\put(160,85){\line(0,1){5}}
	\put(155,85){\line(1,0){5}}
	\put(155,85){\line(0,1){5}}
	\put(155,90){\line(1,0){5}}
	\thinlines 
	\put(135,60){$(2)$}
	\begin{Huge}
	\put(170,60){$\longrightarrow$}
	\end{Huge}
	\begin{footnotesize}
	\put(180,55){$-U_2$}
	\end{footnotesize}
	
	\thinlines
	\put(210,50){\line(1,0){50}}
	\put(260,50){\line(0,1){30}}
	\put(210,50){\line(0,1){30}}
	\put(210,80){\line(1,0){40}}
	\put(255,80){\line(0,1){5}}
	\put(255,80){\line(1,0){5}}
	\put(250,80){\line(0,1){5}}
	\put(250,85){\line(1,0){5}}
	\thinlines 
	\put(230,60){$(3)$}
	\begin{Huge}
	\put(265,60){$\longrightarrow$}
	\end{Huge}
	\begin{footnotesize}
	\put(275,55){$+U_2$}
	\end{footnotesize}
	
	\thinlines
	\put(305,50){\line(1,0){50}}
	\put(305,50){\line(0,1){30}}
	\put(305,80){\line(1,0){40}}
	\put(345,80){\line(0,1){5}}
	\put(345,85){\line(1,0){10}}
	\put(355,50){\line(0,1){25}}
	\put(350,75){\line(1,0){5}}
	\put(350,75){\line(0,1){5}}
	\put(350,80){\line(1,0){5}}
	\put(355,80){\line(0,1){5}}
	\thinlines 
	\put(325,60){$(4)$}
	\begin{Huge}
	\put(360,60){$\dashrightarrow$}
	\end{Huge}
	\begin{footnotesize}
	\put(365,55){$0,..,0$}
	\end{footnotesize}

	\thinlines
	\put(20,0){\line(1,0){50}}
	\put(20,0){\line(0,1){30}}
	\put(70,0){\line(0,1){5}}
	\put(20,30){\line(1,0){40}}
	\put(60,30){\line(0,1){5}}
	\put(60,35){\line(1,0){10}}
	\put(65,5){\line(1,0){5}}
	\put(65,5){\line(0,1){5}}
	\put(65,10){\line(1,0){5}}
	\put(70,10){\line(0,1){25}}
	\thinlines 
	\put(37,13){$(5)$}
	\begin{Huge}
	\put(75,10){$\longrightarrow$}
	\end{Huge}
	\begin{footnotesize}
	\put(83,5){$-U_2$}
	\end{footnotesize}

	\thinlines
	\put(115,0){\line(1,0){45}}
	\put(115,0){\line(0,1){30}}
	\put(160,0){\line(0,1){5}}
	\put(115,30){\line(1,0){40}}
	\put(160,5){\line(1,0){5}}
	\put(165,5){\line(0,1){30}}
	\put(155,35){\line(1,0){10}}
	\put(155,30){\line(0,1){5}}
	\thinlines 
	\put(135,13){$(6)$}
	\begin{Huge}
	\put(170,10){$\longrightarrow$}
	\end{Huge}
	\begin{footnotesize}
	\put(180,5){$+U_1$}
	\end{footnotesize}
	
	\thinlines
	\put(210,0){\line(1,0){45}}
	\put(210,0){\line(0,1){30}}
	\put(210,30){\line(1,0){35}}
	\put(245,30){\line(0,1){5}}
	\put(245,35){\line(1,0){5}}
	\put(250,30){\line(0,1){5}}
	\put(250,30){\line(1,0){5}}
	\put(255,30){\line(0,1){5}}
	\put(255,35){\line(1,0){5}}
	\put(255,0){\line(0,1){5}}
	\put(255,5){\line(1,0){5}}
	\put(260,5){\line(0,1){30}}
	\thinlines 
	\put(230,13){$(7)$}
	\begin{Huge}
	\put(265,10){$\longrightarrow$}
	\end{Huge}
	\begin{footnotesize}
	\put(275,5){$-U_1$}
	\end{footnotesize}

	\thinlines
	\put(305,0){\line(1,0){45}}
	\put(305,0){\line(0,1){30}}
	\put(305,30){\line(1,0){35}}
	\put(340,30){\line(0,1){5}}
	\put(340,35){\line(1,0){10}}
	\put(350,30){\line(0,1){5}}
	\put(350,30){\line(1,0){5}}
	\put(350,0){\line(0,1){5}}
	\put(350,5){\line(1,0){5}}
	\put(355,5){\line(0,1){25}}
	\thinlines 
	\put(325,13){$(8)$}
	\begin{Huge}
	\put(360,10){$\longrightarrow$}
	\end{Huge}
	\begin{footnotesize}
	\put(370,5){$+U_2$}
	\end{footnotesize}

	\thinlines
	\put(20,-50){\line(1,0){45}}
	\put(20,-50){\line(0,1){30}}
	\put(20,-20){\line(1,0){35}}
	\put(55,-20){\line(0,1){5}}
	\put(55,-15){\line(1,0){15}}
	\put(65,-50){\line(0,1){5}}
	\put(65,-45){\line(1,0){5}}
	\put(70,-45){\line(0,1){20}}
	\put(65,-25){\line(1,0){5}}
	\put(65,-25){\line(0,1){5}}
	\put(65,-20){\line(1,0){5}}
	\put(70,-20){\line(0,1){5}}
	\thinlines 
	\put(37,-37){$(9)$}
	\begin{Huge}
	\put(75,-40){$\dashrightarrow$}
	\end{Huge}
	\begin{tiny}
	\put(75,-45){$0,..,0,-U_2$}
	\end{tiny}

	\thinlines
	\put(115,-50){\line(1,0){45}}
	\put(115,-50){\line(0,1){30}}
	\put(115,-20){\line(1,0){35}}
	\put(150,-20){\line(0,1){5}}
	\put(150,-15){\line(1,0){15}}
	\put(160,-50){\line(0,1){10}}
	\put(160,-40){\line(1,0){5}}
	\put(165,-40){\line(0,1){25}}
	\thinlines 
	\put(130,-37){$(10)$}
	\begin{Huge}
	\put(170,-40){$\dashrightarrow$}
	\end{Huge}
	\begin{footnotesize}
	\put(177,-45){$+U_1$}
	\end{footnotesize}

	\thinlines
	\put(210,-50){\line(1,0){45}}
	\put(210,-50){\line(0,1){30}}
	\put(210,-20){\line(1,0){15}}
	\put(225,-20){\line(0,1){5}}
	\put(225,-15){\line(1,0){5}}
	\put(230,-15){\line(0,-1){5}}
	\put(230,-20){\line(1,0){5}}
	\put(235,-20){\line(0,1){5}}
	\put(235,-15){\line(1,0){25}}
	\put(255,-50){\line(0,1){30}}
	\put(255,-20){\line(1,0){5}}
	\put(260,-20){\line(0,1){5}}
	
	\thinlines 
	\put(225,-35){$(11)$}
	\begin{Huge}
	\put(265,-40){$\dashrightarrow$}
	\end{Huge}
	\begin{tiny}
	\put(265,-45){$-U_1+U_1$}
	\end{tiny}

	\thinlines
	\put(305,-50){\line(1,0){45}}
	\put(305,-50){\line(0,1){30}}
	\put(305,-20){\line(1,0){25}}
	\put(330,-20){\line(0,1){5}}
	\put(330,-15){\line(1,0){20}}
	\put(350,-50){\line(0,1){35}}
	\put(355,-20){\line(1,0){5}}
	\put(360,-20){\line(0,1){5}}
	\put(355,-15){\line(1,0){5}}
	\put(355,-20){\line(0,1){5}}
	\put(363,-18){\line(1,0){2}}
	\put(366,-18){\line(1,0){2}}
	\put(370,-18){\line(1,0){2}}
	\put(372,-17){\line(0,1){2}}
	\put(372,-13){\line(0,1){2}}
	\put(372,-10){\line(0,1){2}}
	\put(370,-8){\line(-1,0){2}}
	\put(367,-8){\line(-1,0){2}}
	\put(364,-8){\line(-1,0){2}}
	\put(361,-8){\line(-1,0){2}}
	\put(358,-8){\line(-1,0){2}}
	\put(355,-8){\line(-1,0){2}}
	\put(352,-8){\line(-1,0){2}}
	\put(349,-8){\line(-1,0){2}}
	\put(346,-8){\line(-1,0){2}}
	\put(343,-8){\line(-1,0){2}}
	\put(340,-8){\line(-1,0){2}}
	\put(337,-8){\line(-1,0){2}}
	\put(334,-8){\line(-1,0){2}}
	\put(331,-8){\line(-1,0){2}}
	\put(328,-8){\line(-1,0){2}}
	\put(326,-9){\line(0,-1){2}}
	\put(326,-11){\vector(0,-1){5}}
	\thinlines 
	\put(320,-35){$(12)$}
	
	\end{picture}
	\vskip 3.5 cm
	\caption{Transition column to row for $\cR(2l_2^*-1,l_2^*-1)$: the configuration (1) has energy equal to $\G-\D+U_2-U_1$ and thus the configurations (7), (11) and (12) have energy equal to $\G$. Between the configuration (11) and (12) the protuberance along the horizontal side is attached to the bar decreasing the energy by $U_1$, afterwards the protuberance on the right column is detached from the cluster increasing the energy by $U_1$. In (12) we indicate with a dashed arrow the movement of the free particle until it connects to the cluster that decreases the energy by $U_1+U_2$.}
	\label{fig:columntorow}
\end{figure}
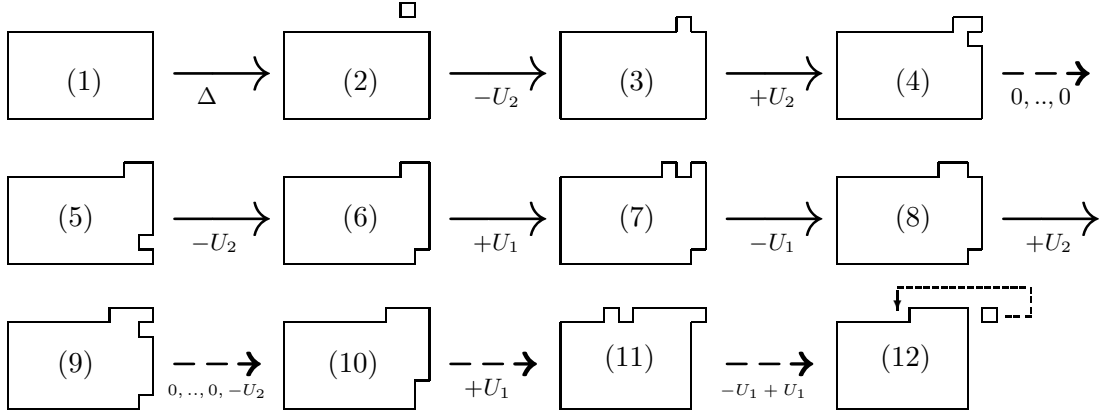

\medskip
We will prove that in cases I, II-C and II-D we have $H(\bar{\h})>\G$; in case II-A we have $H(\bar{\h})\geq\G$ and in cases II-B and III either $\h\in\cB$ or $\h\notin\cB$ and $H(\h)>\G$.

\medskip
\noindent
{\bf Case I.} Since $n(\bar{\h})\geq 1$, by (\ref{energeta}) and (\ref{sistemap2}) we obtain
\be{}
\ba{lll}
H(\bar{\h})\geq H(\cR(2l_2^*+k-1,l_2^*))+\e v(\bar{\h})+\D\geq\\

\geq U_1l_2^*+2U_2l_2^*+kU_2-U_2-\e (2(l_2^*)^2+kl_2^*-l_2^*)+\e (2l_2^*+k-2)+U_1+U_2-\e=\\

=U_1l_2^*+2U_2l_2^*+kU_2-2\e (l_2^*)^2-k\e l_2^*+3\e l_2^*+k\e-3\e +U_1.
\ea
\ee

We recall that $\G=U_1l_2^*+2U_2l_2^*+U_1-U_2-2\e (l_2^*)^2+3\e l_2^*-2\e$. Since $k\geq0$, $\e\ll U_2$ and $\d<1$, we have $U_2>\e(1-k(1-\d))$, that implies $H(\bar{\h})>\G$.

\setlength{\unitlength}{1.1pt}
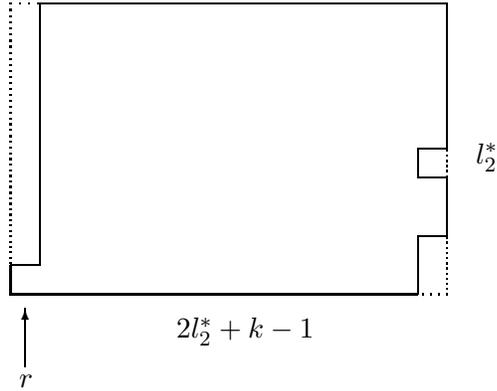
\begin{figure}
	\begin{picture}(400,90)(0,30)
		\thinlines
		\qbezier[51](130,0)(210,0)(280,0)
		\qbezier[51](130,100)(210,100)(280,100)
		\qbezier[51](130,0)(130,45)(130,100)
		\qbezier[51](280,0)(280,45)(280,90)
		\thinlines
		\put(130,0){\line(1,0){140}}
		\put(140,100){\line(1,0){140}}
		\put(140,10){\line(0,1){90}}
		\put(130,0){\line(0,1){10}}
		\put(130,10){\line(1,0){10}}
		\put(270,0){\line(0,1){20}}
		\put(270,20){\line(1,0){10}}
		\put(280,20){\line(0,1){20}}
		\put(280,40){\line(-1,0){10}}
		\put(270,40){\line(0,1){10}}
		\put(270,50){\line(1,0){10}}
		\put(280,50){\line(0,1){50}}
		\thinlines  
		\put(290,45){$l_2^*$}
		\put(187,-15){$2l_2^*+k-1$}
		\put(135,-25){\vector(0,1){20}}
		\put(133,-32){$r$}
		
	\end{picture}
	\vskip 2.5 cm
	\caption{Possible representation for $\bar\h$ in the case $n(\bar\h)=0$, $g_1'(\bar\h)=1$ and $g_2'(\bar\h)=0$.}
	\label{fig:g1'}
\end{figure}

\bigskip
\noindent
{\bf Case II.A} Since $g_2'(\bar{\h})=1$, $k\geq0$ and $\d<1$, using (\ref{energeta}) we obtain
\be{}
\ba{lll}
H(\bar{\h})=H(\cR(2l_2^*+k-1,l_2^*))+\e v(\bar{\h})+U_1\geq\\

\geq U_1l_2^*+2U_2l_2^*+kU_2-U_2-\e (2(l_2^*)^2+kl_2^*-l_2^*)+\e (2l_2^*+k-2)+U_1\geq\G\Leftrightarrow\\

\Leftrightarrow \ k\e(1-\d)\geq0.
\ea
\ee

Our goal is also to emphasize in which pairs $(\bar\h,\h)$ we have $\max\{H(\bar\h),H(\h)\}=\G$. This is the case if $k=0$ and $\bar{\h}$ is a configuration such that $g_2'(\bar{\h})=1$ and $v(\bar{\h})=2l_2^*+k-2=p_{max}(\bar{\h})-1$, namely the configurations in $\cP_1$ (as number $(7)$ and $(11)$ in Figure \ref{fig:columntorow}). Note that if $\bar\h\in\cP_1$, then $\h\in\cB$. Indeed with the move (see transition from configuration $(7)$ to $(8)$ in Figure \ref{fig:columntorow}), we get $p_2(\h)=p_2(\bar\h)$, $s(\h)=s(\bar\h)$ and $v(\h)=v(\bar\h)$. If $v(\bar{\h})>p_{max}(\bar{\h})-1$ or $k\geq1$, we get $H(\bar{\h})>\G.$

\medskip
\noindent
{\bf Case II.B} We have $g_1'(\bar{\h})=1$. By Remark \ref{remarkulteriore} we know that no line can become active. If $p_2(\h)=p_2(\bar\h)-1=l_2^*-1$, we get $\h\in\cB$. Thus we have to consider only the case in which $\D s=-1$ is obtained by a vertical line becoming inactive. Since by Lemma \ref{lemma1} the lines $r_2$ and $r_5$ can not become inactive, the deactivating line is one among $r_6$ and $r_7$. We distinguish the case $\bar\h_{cl}$ connected or not.

If $\bar\h_{cl}$ is connected we note that $g_1'(\bar\h)=1$ is given by two protuberances at distance strictly bigger than one on the shorter side, that is the vertical one. Indeed by definition (\ref{defB}) in the case $s\geq s^*$, we deduce that $p_1(\bar\h)>p_2(\bar\h)$ (see Figure \ref{fig:g1'}). Clearly, there could be more vacancies in $\bar\h$, but what is relevant to obtain $\D s=-1$ is that there must exist a vertical line $r$ such that $r\cap\bar\h_{cl}$ is a single site $a$, such that $a$ becomes free in $\h$, indeed this is the only admissible operation with $\D s=-1$. Since $\bar\h\in\cB$ and $s(\bar\h)\geq s^*$, we note that $v(\bar\h)\geq p_{max}(\bar\h)-1=2l_2^*+k-2$. Since with the move we lose $(l_2^*-1)$ vacancies, we get 

$$v(\h)=v(\bar\h)-(l_2^*-1)\geq l_2^*+k-1.$$

If $k=0$, we get $v(\h)\geq l_2^*-1=p_{min}(\h)-1$ and, since $s(\h)=s^*-1$, we deduce that $\h\in\cB$.

If $k\geq1$, we get
\be{}
\ba{lll}
H(\h)\geq H(\cR(2l_2^*+k-2,l_2^*))+\e(l_2^*+k-1)+\D=\\

=U_1l_2^*+2U_2l_2^*-\e(2(l_2^*)^2+kl_2^*-2l_2^*)+\e l_2^*+U_1-U_2+kU_2+k\e-2\e>\G\Leftrightarrow\\

\Leftrightarrow \ k\e(1-\d)>0, \quad \hbox{always since } \d<1.
\ea
\ee

If $\bar\h_{cl}$ is not connected, we note that $g_1'(\bar\h)=1$ can be obtained with two protuberances at distance strictly bigger than one on the vertical side (as before) or with two distinct clusters that intersect the same vertical line. In the first case we can conclude with a similar argument as before. In the latter case, in order to obtain $\D s=-1$ with the move, we have to move a particle in such a way it becomes free in $\h$ and the vertical line $r$ where it lies in $\bar\h$ becomes inactive. In this way we deduce that $r\cap\bar\h_{cl}$ consists of a single non-empty site $a$ and such a particle has to become free in $\g$. From this we deduce that the particle in $a$ is the moving particle and again in $\h$ we lose $(l_2^*-1)$ vacancies. Thus the situation for $\bar\h$ is the same as before: for $k=0$, we get $\h\in\cB$, and for $k\geq1$, we get $H(\h)>\G$.

\setlength{\unitlength}{1.5pt}
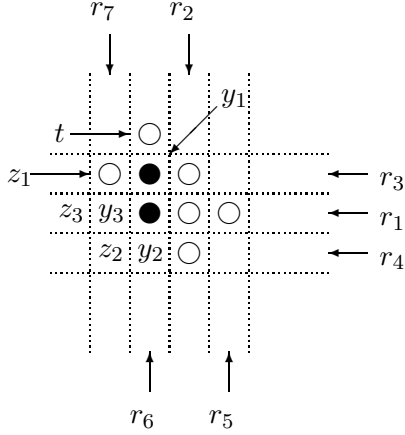
\begin{figure}
	\begin{picture}(380,60)(-30,40)
		\thinlines 
		\qbezier[51](90,0)(125,0)(160,0)
		\qbezier[51](90,10)(125,10)(160,10)
		\qbezier[51](90,20)(125,20)(160,20)
		\qbezier[51](90,30)(125,30)(160,30)
		\qbezier[51](110,-20)(110,15)(110,50)
		\qbezier[51](100,-20)(100,15)(100,50)
		\qbezier[51](120,-20)(120,15)(120,50)
		\qbezier[51](130,-20)(130,15)(130,50)
		\qbezier[51](140,-20)(140,15)(140,50) 
		\put(125,5){\circle{5}}
		\put(125,15){\circle{5}} 
		\put(125,25){\circle{5}}
		\put(115,15){\circle*{5}} 
		\put(115,25){\circle*{5}}
		\put(102,14){$y_3$}
		\put(102,4){$z_2$}
		\put(112,4){$y_2$}
		\put(115,35){\circle{5}}
		\put(105,25){\circle{5}}
		\put(85,25){\vector(1,0){15}}
		\put(79,23){$z_1$}
		\put(95,35){\vector(1,0){15}}
		\put(91,33){$t$}
		\put(132,42){\vector(-1,-1){12}}
		\put(133,43){$y_1$}
		\put(92,14){$z_3$}
		\put(135,15){\circle{5}}
		\thinlines
		\put(125,60){\vector(0,-1){10}} \put(120,65){$r_2$}
		\put(105,60){\vector(0,-1){10}} \put(100,65){$r_7$}
		\put(135,-30){\vector(0,1){10}} \put(130,-38){$r_5$}
		\put(115,-30){\vector(0,1){10}} \put(110,-38){$r_6$}
		\put(170,5){\vector(-1,0){10}} \put(173,2){$r_4$}
		\put(170,15){\vector(-1,0){10}} \put(173,12){$r_1$}
		\put(170,25){\vector(-1,0){10}} \put(173,22){$r_3$}
		
	\end{picture}
	\vskip 4. cm
	\caption{Here we depict part of the configuration in the case $n(\bar\h)=0$, $g_1'(\bar\h)+g_2'(\bar\h)=0$ and $y_1$ occupied.}
	\label{fig:r6prima}
\end{figure}

\medskip
\noindent
{\bf Case II.C} By (\ref{energeta}), the energy increases by a quantity $U_1+z>U_1$, where $z=U_2$ if $g_1'(\bar\h)=1$ and $g_2'(\bar\h)=1$, or $z\geq2U_1$ if $g_2'(\bar\h)\geq2$. Thus with a similar reasoning as in the case II.A, we get $H(\bar{\h})>\G$.

\medskip
\noindent
{\bf Case II.D} We argue in a similar way as in the case II-B. 

\bigskip
{\bf Case III.} We have $s(\bar{\h})\geq s^*$, with $n(\bar{\h})=0$ and $g_1'(\bar{\h})+g_2'(\bar{\h})=0$. We will prove that if $\h\notin\cB$, then $H(\h)>\G$. We recall that the circumscribed rectangle of $\bar\h$ is $\cR(2l_2^*+k-1,l_2^*)$, for any $k\geq0$, and, since $\bar\h\in\cB$, $v(\bar\h)\geq p_{max}(\bar\h)-1=2l_2^*+k-1$. By Remark \ref{remark1} it is impossible to activate lines and thus $\D s=-1$ is obtained by a unique line becoming inactive (see 	Figure \ref{fig:rappresentazione}). Since by Lemma \ref{lemma1} the lines $r_2$ and $r_5$ can not become inactive, we analyze separately the case in which $r_i$ is the line becoming inactive, with $i\in{\{1,3,4,6,7\}}$. If the line becoming inactive is either $r_1$ or $r_3$ or $r_4$, by Lemma \ref{decrescita} we know that $p_2(\h)=p_2(\bar\h)-1=l_2^*-1$, thus $\h\in\cB$. Hence we are left to consider the case in which the line becoming inactive is either $r_6$ or $r_7$:

\bi

\item $r_6$ is the line becoming inactive. Since $\D s=-1$, the sites $x_3$, $x_4$ and $x_5$ must be empty, otherwise some of the lines $r_2$, $r_4$, $r_4$ and $r_5$ are active in $\h$, so that they must be active also in $\bar\h$ (see Figure \ref{fig:rappresentazione}). This implies that $n(\bar{\h})\geq1$ or $\bar{\h}$ is not monotone, which is in contradiction with the assumptions $n(\bar{\h})=0$ and $g_1'(\bar{\h})+g_2'(\bar{\h})=0$. We distinguish the following cases:
\bi
\item [-)] If the site $y_1$ is occupied in $\bar{\h}$, since $r_6\cap\h_{cl}=\emptyset$, we necessarily have that the site $y_1$ must contain a free particle in $\h$. Thus the sites $z_1$ and $t$ are empty (see Figure \ref{fig:r6prima}). Since $\bar\h$ is monotone by assumption, along the line $r_3$ there must be a unique particle in $\bar\h_{cl}$, that is in the site $y_1$. This implies that with the move also the line $r_3$ becomes inactive, which is in contradiction with $\D s=-1$.

\item [-)] If the site $y_2$ is occupied, by symmetry we can conclude with a similar argument as before.
\ei

Hence we are left to analyze the case $y_1$ and $y_2$ empty, that implies $y_3$ occupied in order that $r_6$ becomes inactive. Thus $r_6\cap\bar{\h}_{cl}$ consists in the moving particle, otherwise the line $r_6$ can not become inactive with the move. What is relevant is that $\bar\h_{cl}$ has exactly $p_2(\bar{\h})-1=l_2^*-1$ vacancies along the line $r_6$ (see Figure \ref{fig:r6ultima} on the left hand-side). If there ore other $m$ vacancies in $\bar\h$ additional to the ones along the line $r_6$, there must be $m\geq l_2^*+k-1$. Indeed, in order that $\bar{\h}\in\cB$, we require $v(\bar{\h})\geq p_{max}(\bar{\h})-1$, so we get

$$v(\bar{\h})=m+l_2^*-1\geq p_{max}(\bar{\h})-1=2l_2^*+k-2\Leftrightarrow \ m\geq l_2^*+k-1.$$

Since $x_3$, $x_4$ and $x_5$ are empty, the moved particle is free in $\h$, i.e., $n(\h)\geq1$ (see Figure \ref{fig:r6ultima} on the right hand-side). We distinguish the case $k=0$ and $k\geq1$. If $k=0$, we obtain $\h\in\cB$, because $v(\h)=m\geq l_2^*-1=p_{min}(\bar{\h})-1$. If $k\geq1$, we obtain
\be{sistema}\left\{ \ba{lll}
s(\h)&=&l_2^*+2l_2^*+k-2\geq3l_2^*-1=s^* \\
p_2(\h)&=&l_2^* \\
v(\h)&=&m\geq l_2^*+k-1.
\ea
\right.
\ee

\setlength{\unitlength}{1.1pt}
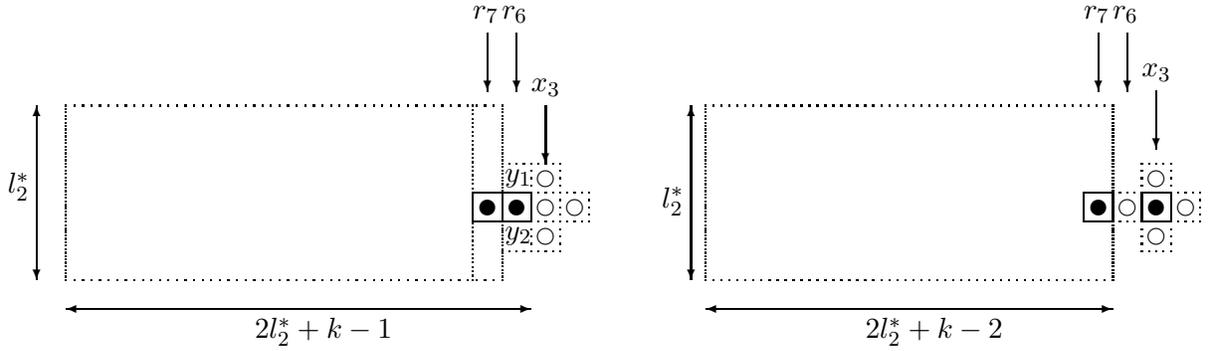
\begin{figure}
	\begin{picture}(380,60)(-30,40)
	
	\thinlines 
	\qbezier[51](-10,0)(-10,30)(-10,60)
	\qbezier[51](-10,0)(65,0)(140,0)
	\qbezier[51](-10,60)(65,60)(140,60)
	\qbezier[8](140,10)(150,10)(160,10)
	\qbezier[8](140,40)(150,40)(160,40)
	\qbezier[12](150,10)(150,25)(150,40)
	\qbezier[8](150,30)(160,30)(170,30)
	\qbezier[12](160,10)(160,25)(160,40)
	\qbezier[8](150,20)(160,20)(170,20)
	\qbezier[4](170,20)(170,25)(170,30)
	\qbezier[4](150,40)(155,40)(160,40)
	\qbezier[40](130,0)(130,30)(130,60)
	\qbezier[51](140,0)(140,30)(140,60)
	\put(141,34){$y_1$}
	\put(141,14){$y_2$}
	\put(140,20){\line(0,1){10}}
	\put(130,20){\line(1,0){20}}
	\put(150,20){\line(0,1){10}}
	\put(150,30){\line(-1,0){20}}
	\put(130,20){\line(0,1){10}}
	\put(145,85){\vector(0,-1){20}}
	\put(140,90){$r_6$}
	\put(135,85){\vector(0,-1){20}}
	\put(130,90){$r_7$}
	\put(155,60){\vector(0,-1){20}}
	\put(150,65){$x_3$}
	\put(135,25){\circle*{5}}
	\put(145,25){\circle*{5}}
	\put(155,25){\circle{5}} 
	\put(165,25){\circle{5}}
	\put(155,35){\circle{5}} 
	\put(155,15){\circle{5}}
	\put(-20,30){\vector(0,1){30}}
	\put(-20,30){\vector(0,-1){30}}
	\put(-30,30){$l_2^*$}
	\put(70,-10){\vector(1,0){80}}
	\put(70,-10){\vector(-1,0){80}}
	\put(55,-20){$2l_2^*+k-1$}

	\thinlines
	\qbezier[51](210,0)(210,30)(210,60)
	\qbezier[51](210,0)(280,0)(350,0)
	\qbezier[51](350,0)(350,30)(350,60)
	\qbezier[51](210,60)(280,60)(350,60)
	\qbezier[12](350,30)(365,30)(380,30)
	\qbezier[12](350,20)(365,20)(380,20)
	\qbezier[4](380,20)(380,25)(380,30)
	\qbezier[4](360,40)(365,40)(370,40)
	\qbezier[4](360,10)(365,10)(370,10)
	\qbezier[12](360,10)(360,25)(360,40)
	\qbezier[12](370,10)(370,25)(370,40)
	\put(370,20){\line(0,1){10}}
	\put(360,20){\line(0,1){10}}
	\put(360,30){\line(1,0){10}}
	\put(360,20){\line(1,0){10}}
	\put(365,25){\circle*{5}}
	\put(355,25){\circle{5}}
	\put(375,25){\circle{5}}
	\put(345,25){\circle*{5}}
	\put(340,20){\line(1,0){10}}
	\put(340,20){\line(0,1){10}}
	\put(350,20){\line(0,1){10}}
	\put(340,30){\line(1,0){10}}
	\put(365,15){\circle{5}}
	\put(365,35){\circle{5}}
	\put(205,35){\vector(0,1){25}}
	\put(205,35){\vector(0,-1){35}}
	\put(195,25){$l_2^*$}
	\put(290,-10){\vector(1,0){60}}
	\put(290,-10){\vector(-1,0){80}}
	\put(265,-20){$2l_2^*+k-2$}
	\put(355,85){\vector(0,-1){20}}
	\put(350,90){$r_6$}
	\put(345,85){\vector(0,-1){20}}
	\put(340,90){$r_7$}
	\put(365,65){\vector(0,-1){20}}
	\put(360,70){$x_3$}
	\end{picture}
	\vskip 3. cm 
	\caption{Possible configurations for $\bar{\h}$ (on the left hand-side) and the corresponding for $\h$ (on the right hand-side) in the case that $r_6$ becomes inactive.}
	\label{fig:r6ultima}
\end{figure}

Since $l_2^*\gg1$, we have that $v(\h)<p_{max}(\h)-1$, thus $\h\notin\cB$. By (\ref{sistema}) and (\ref{energeta}), we get $H(\h)>\G$, indeed
\be{5.6}
\ba{lll}
H(\h)\geq H(\cR(2l_2^*+k-2,l_2^*))+\e m+\D\geq\\

\geq U_1l_2^*+2U_2l_2^*+kU_2-\e(2(l_2^*)^2+kl_2^*-2l_2^*)+\e l_2^*+k\e+U_1-U_2-2\e>\G\Leftrightarrow\\

\Leftrightarrow \e k(1-\d)>0, \quad \hbox{always since } k\geq1 \hbox{ and } \d<1.
\ea
\ee

\item $r_7$ is the line becoming inactive. We deduce that the site $y_3$ is occupied and the sites $z_1$, $z_2$ and $z_3$ are empty. Furthermore, since $\bar\h$ is monotone and $n(\bar\h)=0$, we obtain that the site $x_5$ is empty. Similarly, since by assumptions $n(\bar\h)=0$ and $g_1'(\bar\h)=g_2'(\bar\h)=0$, we deduce that the sites $x_3$ and $x_4$ can not be both occupied. Moreover, the sites $x_3$ and $x_4$ can not be both empty, otherwise the line $r_1$ becomes inactive, which is in contradiction with $\D s=-1$ (see Figure \ref{fig:r7}, in which we assume without loss of generality $x_3$ occupied and $x_4$ empty). Since $\bar\h$ is monotone and $n(\bar\h)=0$, $r_1\cap\bar\h_{cl}$ consists in the two sites $x_1$ and $y_3$. The circumscribed rectangle of $\bar\h$ is $\cR(2l_2^*+k-1,l_2^*)$ for any $k\geq0$, thus by Lemma \ref{decrescita} we have that the circumscribed rectangle of $\h$ is $\cR(2l_2^*+k-2,l_2^*)$ for any $k\geq0$. If either $\bar\h_{cl}$ is connected or both $\bar\h_{cl}$ and $\h_{cl}$ are not connected, at most five vacancies in $\bar\h$ are not vacancies in $\h$: the sites $z_1$, $z_2$, $a$, $b$ and $c$. Thus we obtain $v(\h)\geq v(\bar\h)-5\geq 2l_2^*+k-7$. If $k=0$ we have $s(\bar\h)=s^*-1$. Since $U_2\gg\e$, we obtain that $\h\in\cB$, indeed

$$v(\h)\geq 2l_2^*-7\geq p_{min}(\h)-1=l_2^*-1\Leftrightarrow \ l_2^*\geq6.$$

\bigskip
\setlength{\unitlength}{0.8pt}
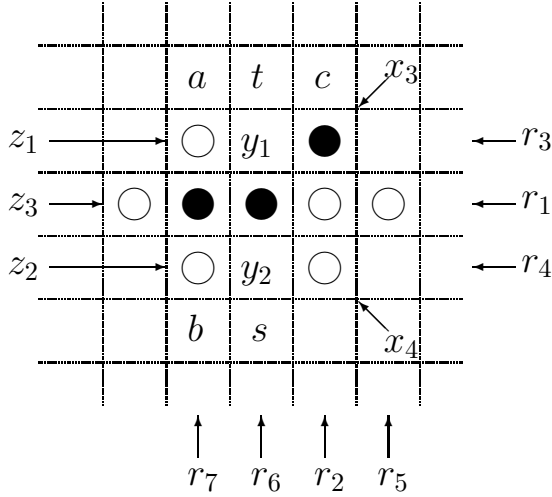
\begin{figure}
	\begin{center}
		\begin{picture}(380,60)(-30,40)
		
		\thinlines 
		\qbezier[150](90,0)(190,0)(290,0)
		\qbezier[150](90,30)(190,30)(290,30) \qbezier[150](90,60)(190,60)(290,60)
		\qbezier[150](90,90)(190,90)(290,90)
		\qbezier[150](90,-30)(190,-30)(290,-30)
		\qbezier[150](90,-60)(190,-60)(290,-60)
		\qbezier[150](240,-80)(240,15)(240,110)
		\qbezier[150](270,-80)(270,15)(270,110)
		\qbezier[150](120,-80)(120,15)(120,110)
		\qbezier[150](150,-80)(150,15)(150,110)
		\qbezier[150](180,-80)(180,15)(180,110)
		\qbezier[150](210,-80)(210,15)(210,110)
		\put(195,15){\circle*{15}}
		\put(225,15){\circle{15}}
		\begin{Large}
		\put(95,15){\vector(1,0){25}}
		\put(75,12){$z_3$}
		\put(95,45){\vector(1,0){55}}
		\put(75,42){$z_1$}
		\put(95,-15){\vector(1,0){55}}
		\put(75,-18){$z_2$}
		\put(255,-45){\vector(-1,1){14}}
		\put(253,-55){$x_4$}
		\put(255,75){\vector(-1,-1){14}}
		\put(253,75){$x_3$}
		\put(225,45){\circle*{15}}
		\put(225,-15){\circle{15}}
		\put(255,15){\circle{15}}
		\put(185,40){$y_1$}
		\put(185,-20){$y_2$}
		\put(190,70){$t$}
		\put(220,70){$c$}
		\put(160,70){$a$}
		\put(160,-50){$b$}
		\put(190,-50){$s$}
		\put(165,15){\circle*{15}}
		\put(135,15){\circle{15}}
		\put(165,-15){\circle{15}}
		\put(165,45){\circle{15}}
		
		\thinlines
		\put(225,-105){\vector(0,1){20}} \put(220,-118){$r_2$}
		\put(165,-105){\vector(0,1){20}} \put(160,-118){$r_7$}
		
		\put(250,-118){$r_5$} \put(255,-105){\vector(0,1){20}}
		\put(195,-105){\vector(0,1){20}} \put(190,-118){$r_6$}
		\put(315,-15){\vector(-1,0){20}} \put(318,-17){$r_4$}
		\put(315,15){\vector(-1,0){20}} \put(318,13){$r_1$}
		\put(315,45){\vector(-1,0){20}} \put(318,43){$r_3$}
		\end{Large}
		\end{picture}
		\vskip 4.5 cm
		\caption{Here we depict part of the configuration in the case in which $r_7$ becomes inactive for $\D s=-1$.}
		\label{fig:r7}
	\end{center}
\end{figure}

If $k\geq1$, we obtain
\be{sistema2}\left\{ \ba{lll}
n(\h)&\geq&1 \\
v(\h)&\geq&2l_2^*+k-7\geq2l_2^*-6\geq l_2^*.
\ea
\right.
\ee

By (\ref{sistema2}), we get $H(\h)\geq H(\cR(2l_2^*+k-2,l_2^*))+\e(2l_2^*-6)+\D\geq H(\cR(2l_2^*+k-2,l_2^*))+\e l_2^*+\D>\G$, where the last inequality follows by (\ref{5.6}).

The case $\bar\h_{cl}$ not connected and $\h_{cl}$ connected is not admissible. Indeed, if $\bar\h$ is not connected, along the line $r_6$ there must be at least other two particles in addiction to the moving one, otherwise $r_6$ becomes inactive with the move, which is in contradiction with $\D s=-1$. In this way, since the site $r_1\cap r_6$ is empty in $\h$ and the configuration $\bar\h$ is monotone, we deduce that $\h_{cl}$ is non connected.

\ei

\epr

\br{-1}
Let $(\bar\h,\h)$ be a move with $\bar\h\in\cB$ and let $\D s=-1$ be its corresponding variation of $s$; if $H(\bar\h)=\G$ then  $\bar\h\in\cP_1$ and $\h\in\cB$. 

\er

\bpr
We can justify this remark with the following argument. In the proof of Proposition \ref{lemma5} we have underlined which configurations $\bar{\h}\in\cB$ have energy equal to $\G$: now we want to analyze which of them are the saddles such that $\h\notin\cB$. Starting from one of these possible configurations for $\bar{\h}$, we analyze the moves with minimal cost from an energetical point of view. What is relevant is that the energy of the configuration must not exceed the value $\G$. Thus we can move the protuberance along the vertical side: this is a zero cost move and does not allow the exit from $\cB$. Another operation that we can make is the connection between the two protuberances of the configuration (see the transition from (7) to (8) in Figure \ref{fig:columntorow}), since the energy decreases by a quantity $U_1$. These configurations have $p_{max}(\bar{\h})-1$ vacancies and thus belong to $\cB$. Moreover, we could move the particles that lie along the last column by an upward position (see the transitions from (8) to (10) in Figure \ref{fig:columntorow}): all together these operations have an energy cost equal to zero. More precisely, when we move the first particle the energy increases by $U_2$, then the subsequent moves do not change the energy until we move the last particle and it decreases by $U_2$. Let $\h'$ such a configuration (as (10) in Figure \ref{fig:columntorow}), so we have $H(\h')=\G-U_1+U_2+0+..+0-U_2=\G-U_1<\G$, with $v(\h')=p_{max}(\h')-1=2l_2^*-2$. Thus the configuration $\h'$ belongs to $\cB$. The real saddles $\bar\h_k$, depicted in $(12)$ in Figure \ref{fig:columntorow}, that allow to arrive to a configuration $\h_k\notin\cB$, are the ones where the free particle moves along the dashed arrow, while $\h_k$ is the configuration where the free particle is attached to the cluster. Note that in this last step of the path described $\D s=0$, which contradicts $\D s=-1$. We get $n(\bar\h_k)=1$, $p_1(\bar\h_k)=2l_2^*-2$, $p_2(\bar\h_k)=l_2^*$ and $v(\bar\h_k)=l_2^*-1$, thus $\bar\h_k\in\cP_2\subset\cB$. Furthermore, we get $n(\h_k)=0$, $p_1(\h_k)=2l_2^*-2$, $p_2(\h_k)=l_2^*$ and $v(\h_k)=l_2^*-2$, thus $\h_k\notin\cB$.
\epr

\subsection{Proof of Proposition \lowercase{\ref{lemma6}}}
\label{dim6}

\bpr
Let $\D s=0$ and $(\bar\h,\h)\in\partial\cB$ the exiting move from $\cB$. By Remark \ref{remark2} we consider only the case $s(\bar\h)\geq s^*-1$. By definition (\ref{defB}) and $s(\bar\h)\geq s^*-1$, we deduce that only the two following cases are admissible:

\bi
\item [(a)] $s(\bar{\h})=s^*-1$ and $v(\bar\h)\geq p_{min}(\bar\h)-1$;
\item [(b)] $s(\bar{\h})\geq s^*$, $p_2(\bar\h)=l_2^*$ and $v(\bar\h)\geq p_{max}(\bar\h)-1$.
\ei

If $p_2(\h)\leq l_2^*-1$, in both cases (a) and (b) we get $\h\in\cB$. Hence we are left to consider $p_2(\h)\geq l_2^*$. Let $\D v=v(\h)-v(\bar\h)$ the variation of the number of vacancies with the move. If $\D v>-1$ we have:

\bi
\item [(a)] $v(\h)>v(\bar{\h})-1\geq p_{min}(\bar{\h})-2\geq p_{min}(\h)-2$ that implies $v(\h)\geq p_{min}(\h)-1$, thus by (\ref{defB}) we get $\h\in\cB$.
\item [(b)] $v(\h)>v(\bar{\h})-1\geq p_{max}(\bar{\h})-2\geq p_{max}(\h)-2$ that implies $v(\h)\geq p_{max}(\h)-1$, thus, if $p_2(\h)=l_2^*$, by (\ref{defB}) we get $\h\in\cB$. Hence we are left to analyze the case $p_2(\h)\geq l_2^*+1$; by Remark \ref{crescita} we know that $l_2^*+1\leq p_2(\h)\leq l_2^*+3$, since at most three horizontal lines can become active in order that $\D s=0$. For each horizontal line that becomes active we have a free particle in $\bar\h$, thus we get $1\leq n(\bar\h)\leq3$. We note that the circumscribed rectangle of $\bar\h$ is $\cR(2l_2^*+k-1,l_2^*)$ for any $k\geq0$. Since $k\geq0$, $\d<1$ and $U_2\gg\e$, we get $H(\bar\h)>\G$. Indeed by (\ref{energeta}) we have
\be{}
\ba{lll}
H(\bar\h)\geq H(\cR(2l_2^*+k-1,l_2^*))+\e(2l_2^*+k-2)+\D=\\

U_1l_2^*+2U_2l_2^*+U_1+kU_2-\e(2(l_2^*)^2+kl_2^*-l_2^*)+2\e(l_2^*)+k\e-3\e>\G\Leftrightarrow\\

\Leftrightarrow \ U_2>\e(1+k(\d-1)).
\ea
\ee

\ei

\noindent
Thus we can consider both cases (a) and (b) with the further condition that $\D v\leq-1$. Since the number of vacancies can decrease only if either $p_1(\h)-p_1(\bar{\h})=p_2(\bar{\h})-p_2(\h)\neq0$ or $p_1(\h)-p_1(\bar{\h})=p_2(\bar{\h})-p_2(\h)=0$ but $|\h_{cl}|-|\bar{\h}_{cl}|>0$ that implies that in both cases (a) and (b) we have:

\bi
\item [(i)] either $n(\bar{\h})\geq1$ or
\item [(ii)] $n(\bar{\h})=0$ and, by Remark \ref{remark1}, $g_1'(\bar{\h})+g_2'(\bar{\h})\geq1$ and $r_2$ becomes active bringing at least $p_{min}(\h)-1$ vacancies in $\h.$
\ei

The cases (b-i) and (b-ii) can be treated as in point $\D s=-1$. We refer to the Appendix for explicit computations.

The case (a) is compatible only with case (i), since in the case (a-ii) we have that $\h\in\cB$, indeed
\be{} \left\{
\ba{lll}
v(\h)&\geq&p_{min}(\h)-1 \\
s(\h)&=&s(\bar{\h})=s^*-1.
\ea
\right.
\ee

\noindent
Hence we are left to analyze the case (a-i), with $p_2(\bar{\h})\geq l_2^*$. We note that the circumscribed rectangle of $\bar{\h}$ is $\cR(2l_2^*-k-2,l_2^*+k)$ for any $k\geq0$. Since $U_1-2U_2>0$ and $\d-k-3<-2$, we have that $H(\bar{\h})\geq\G$. Indeed
\be{}
\ba{llll}
H(\bar{\h})\geq H(\cR(2l_2^*-k-2,l_2^*+k))+\e(l_2^*+k-1)+\D=\\

=U_1l_2^*+2U_2l_2^*+U_1-U_2+k(U_1-U_2)-2\e(l_2^*)^2-\e kl_2^*+\e k^2+3\e l_2^*+3\e k-2\e>\G\\

\Leftrightarrow\e k^2+k[U_1-U_2+3\e-\e l_2^*]>0 \\

\Leftrightarrow \ \frac{U_1-2U_2}{\e}>\d-k-3.
\ea
\ee

\noindent
In particular, we obtain $H(\bar{\h})=\G$ if $k=0$, $n(\bar\h)=1$ and $v(\bar{\h})=l_2^*-1$, otherwise $H(\bar{\h})>\G$. We note that $H(\bar\h)=\G$ for every $\bar\h\in\cP_2$ (see definition in (\ref{defP2})).

\epr

\subsection{Proof of Proposition \lowercase{\ref{lemma7}}}
\label{dim7}

\bpr 
Let $\D s=1$; by Remark \ref{remark2} we consider only the case $s(\bar\h)\geq s^*-2$. We distinguish the following three cases:

\bi
\item [(a)] $s(\bar\h)=s^*-2$;
\item [(b)] $s(\bar\h)=s^*-1$;
\item [(c)] $s(\bar\h)\geq s^*$.
\ei

First, we consider the case (a). If $s(\bar{\h})=s^*-2$, then $s(\h)=s^*-1$. If $p_2(\h)\leq l_2^*-1$, by (\ref{defB}) we get $\h\in\cB$. If $p_2(\h)\geq l_2^*$ and $\h\notin\cB$, we get $v(\h)<p_{min}(\h)-1$, so by Lemma \ref{lemma3.2}{\it (iii)} we have $n(\bar{\h})\geq2$. Furthermore, the circumscribed rectangle of $\bar{\h}$ is $\cR(2l_2^*-k-3,l_2^*+k)$ for any $k\geq0$. Thus by (\ref{energeta}) we obtain
\be{5**}
\ba{lll}
H(\bar{\h})\geq H(\cR(2l_2^*-k-3,l_2^*+k))+2\D=\\

=U_1l_2^*+kU_1+2U_2l_2^*-kU_2+2U_1-U_2-2\e(l_2^*)^2-\e kl_2^*+\e k^2+3\e l_2^*+3 \e k-2\e>\G\\

\Leftrightarrow \ \e k^2+k[U_1-U_2-\e l_2^*+3\e]+U_1>0.
\ea
\ee

\noindent
Since $k\geq0$, $\e>0$, $U_1>0$ and $\d<1$, we get $H(\bar{\h})>\G$, using

\be{5*}
U_1-U_2-\e l_2^*+3\e=U_1-2U_2-\e\d+3\e>\e(2-\d)>\e>0.
\ee

\noindent
Hence we are left to consider the cases (b) and (c). Again we consider two possibilities:

\begin{itemize}
	\item [(i)] either $n(\bar{\h})\geq1$ or
	\item [(ii)] $n(\bar{\h})=0$ and, by Remark \ref{remark1}, $g_1'(\bar{\h})+g_2'(\bar{\h})\geq1$ and $r_2$ becomes active bringing at least $p_{min}(\h)-1$ vacancies in $\h.$
\end{itemize}

In the case (b-i) we conclude as in the Proposition \ref{lemma6} case (a-i) and we refer to the Appendix for the explicit computation. In particular, we get $H(\bar\h)=\G$ if $\bar\h\in\cP_2$.

In the cases (c-i) and (c-ii), we get $H(\bar{\h})\geq\G$ as in Proposition \ref{lemma5} and we refer to the Appendix for the explicit computations.

We have only to discuss the case (b-ii). In this case the only line that can become active is $r_2$ and $\bar{\h}_{cl}$ is not connected. This also implies that there are no deactivating line with the move. So we have to consider separately the case in which the move is horizontal or vertical.

\medskip
\noindent
{\bf Case (h).} Suppose first that the move is horizontal, i.e., $r_1$ is an horizontal line. Since the line that becomes active is only $r_2$, at least one among the sites $x_3$, $x_4$ and $x_5$ must be occupied, otherwise $r_2$ does not change its behavior. If the site $x_3$ is occupied, since $r_2$ must be inactive in $\bar{\h}$, it contains in $\bar\h$ a free particle, that implies the case (i). Hence the site $x_3$ is empty and by symmetry also the site $x_4$ is empty: this implies that the site $x_5$ is occupied (see Figure \ref{fig:casoh} on the left hand-side). We deduce that $g_2'(\bar{\h})\geq1$. All the sites along the line $r_2$ are empty, otherwise either $n(\bar\h)\geq1$ or $r_2$ is active in $\bar{\h}$. Since $s(\bar{\h})=s^*-1$ and $\bar{\h}\in\cB$, we get $v(\bar{\h})\geq p_{min}(\bar{\h})-1$. We note that the circumscribed rectangle of $\bar{\h}$ is $\cR(2l_2^*-k-2,l_2^*+k)$ for any $k\geq0$. 

\bi
\item [$\bullet$] If $g_2'(\bar{\h})\geq2$, we get
\be{}
\ba{lll}
H(\bar{\h})\geq H(\cR(2l_2^*-k-2,l_2^*+k))+\e(l_2^*+k-1)+2U_1=\\

=U_1l_2^*+k(U_1-U_2)+2U_2l_2^*+2U_1-2U_2-2\e(l_2^*)^2-k\e l_2^*+\e k^2+3\e k-\e>\G\\

\Leftrightarrow \ \e k^2+k[U_1-U_2-\e l_2^*+3\e]+U_1-U_2+\e>0.
\ea
\ee

\noindent
Analougusly to (\ref{5**}), by (\ref{5*}) and $U_1-U_2+\e>0$, we get $H(\bar{\h})>\G$.

\item [$\bullet$] If $g_1'(\bar{\h})\geq1$, by (\ref{5*}) we get $H(\bar{\h})>\G$. Indeed
\be{5.7}
\ba{lll}
H(\bar{\h})\geq H(\cR(2l_2^*-k-2,l_2^*+k))+\e(l_2^*+k-1)+U_1+U_2=\\

=U_1l_2^*+2U_2l_2^*+k(U_1-U_2)+U_1-U_2-2\e(l_2^*)^2-\e kl_2^*+\e k^2+3\e l_2^*+3\e k-\e>\G\\

\Leftrightarrow \ \e k^2+k[U_1-U_2-\e l_2^*+3\e]+\e>0.
\ea
\ee
\ei

It remains to analyze the case $g_2'(\bar{\h})=1$ and $g_1'(\bar{\h})=0$. In this case it is necessary to analyze in more detail the geometry of the configuration $\bar{\h}$.

\bi
\item [-] If the moving particle has in $\bar{\h}$ at least one vertical and one horizontal bond that connects it to other particles ($y_4$ and one among $y_1$ and $y_2$, see Figure \ref{fig:casoh} on the left hand-side). Thus $\D H:=H(\h)-H(\bar{\h})\geq U_2$, since in the move at least one vertical bond is lost and the horizontal bond with $y_3$ is recovered with the one with $x_5$. This implies by (\ref{5.7}) that $H(\h)\geq H(\bar{\h})+U_2\geq H(\cR(2l_2^*-k-2,l_2^*+k))+\e(l_2^*+k-1)+U_1+U_2>\G$.

\setlength{\unitlength}{1.1pt}
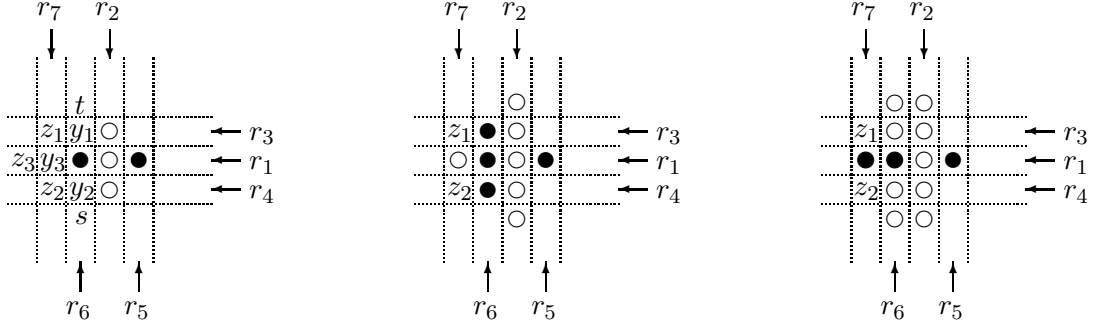
\begin{figure}
	\begin{picture}(380,60)(-30,40)
	\thinlines \qbezier[51](0,0)(35,0)(70,0)
	\qbezier[51](0,10)(35,10)(70,10) \qbezier[51](0,20)(35,20)(70,20)
	\qbezier[51](0,30)(35,30)(70,30)
	\qbezier[51](10,-20)(10,15)(10,50)
	\qbezier[51](20,-20)(20,15)(20,50)
	\qbezier[51](30,-20)(30,15)(30,50)
	\qbezier[51](40,-20)(40,15)(40,50)
	\qbezier[51](50,-20)(50,15)(50,50)
	\put(21,23){$y_1$}
	\put(35,15){\circle{5}}
	\put(35,25){\circle{5}}
	\put(45,15){\circle*{5}}
	\put(25,15){\circle*{5}}
	\put(23,31){$t$}
	\put(21,3){$y_2$}
	\put(11,23){$z_1$}
	\put(11,13){$y_3$}
	\put(23,-7){$s$}
	\put(11,3){$z_2$}
	\put(1,13){$z_3$}
	\put(35,5){\circle{5}}
	\thinlines
	\put(35,60){\vector(0,-1){10}} \put(30,65){$r_2$}
	\put(15,60){\vector(0,-1){10}} \put(10,65){$r_7$}
	\put(40,-38){$r_5$}\put(45,-30){\vector(0,1){10}}
	\put(25,-30){\vector(0,1){10}} \put(20,-38){$r_6$}
	\put(80,5){\vector(-1,0){10}} \put(83,2){$r_4$}
	\put(80,15){\vector(-1,0){10}} \put(83,12){$r_1$}
	\put(80,25){\vector(-1,0){10}} \put(83,22){$r_3$}
	\thinlines 
	\qbezier[51](140,0)(175,0)(210,0)
	\qbezier[51](140,10)(175,10)(210,10)
	\qbezier[51](140,20)(175,20)(210,20)
	\qbezier[51](140,30)(175,30)(210,30)
	\qbezier[51](160,-20)(160,15)(160,50)
	\qbezier[51](150,-20)(150,15)(150,50)
	\qbezier[51](170,-20)(170,15)(170,50)
	\qbezier[51](180,-20)(180,15)(180,50)
	\qbezier[51](190,-20)(190,15)(190,50) 
	\put(175,5){\circle{5}}
	\put(175,15){\circle{5}} 
	\put(175,25){\circle{5}}
	\put(165,15){\circle*{5}} 
	\put(165,25){\circle*{5}}
	\put(155,15){\circle{5}}
	\put(151,3){$z_2$}
	\put(165,5){\circle*{5}}
	\put(175,35){\circle{5}}
	\put(175,-5){\circle{5}}
	\put(151,23){$z_1$}
	\put(185,15){\circle*{5}}
	\thinlines
	\put(175,60){\vector(0,-1){10}} \put(170,65){$r_2$}
	\put(155,60){\vector(0,-1){10}} \put(150,65){$r_7$}
	\put(185,-30){\vector(0,1){10}} \put(180,-38){$r_5$}
	\put(165,-30){\vector(0,1){10}} \put(160,-38){$r_6$}
	\put(220,5){\vector(-1,0){10}} \put(223,2){$r_4$}
	\put(220,15){\vector(-1,0){10}} \put(223,12){$r_1$}
	\put(220,25){\vector(-1,0){10}} \put(223,22){$r_3$}
	\thinlines
	\qbezier[51](280,0)(315,0)(350,0)
	\qbezier[51](280,10)(315,10)(350,10)
	\qbezier[51](280,20)(315,20)(350,20)
	\qbezier[51](280,30)(315,30)(350,30)
	\qbezier[51](300,-20)(300,15)(300,50)
	\qbezier[51](290,-20)(290,15)(290,50)
	\qbezier[51](310,-20)(310,15)(310,50)
	\qbezier[51](320,-20)(320,15)(320,50)
	\qbezier[51](330,-20)(330,15)(330,50)
	
	\put(315,15){\circle{5}} 
	\put(295,15){\circle{5}}
	\put(305,15){\circle{5}} 
	\put(315,5){\circle{5}}
	\put(315,25){\circle{5}}
	\put(305,35){\circle{5}}
	\put(291,23){$z_1$}
	
	\put(305,15){\circle*{5}} 
	\put(295,15){\circle*{5}}
	\put(291,3){$z_2$}
	\put(325,15){\circle*{5}}
	\put(315,35){\circle{5}}
	\put(315,-5){\circle{5}}
	\put(305,5){\circle{5}} 
	\put(305,25){\circle{5}} 
	\put(305,-5){\circle{5}}
	\thinlines
	\put(315,60){\vector(0,-1){10}} \put(310,65){$r_2$}
	\put(325,-30){\vector(0,1){10}} \put(300,-38){$r_6$}
	\put(305,-30){\vector(0,1){10}} \put(320,-38){$r_5$}
	\put(360,5){\vector(-1,0){10}} \put(363,2){$r_4$}
	\put(360,15){\vector(-1,0){10}} \put(363,12){$r_1$}
	\put(360,25){\vector(-1,0){10}} \put(363,22){$r_3$}
	\put(295,60){\vector(0,-1){10}} \put(290,65){$r_7$}
	
	\end{picture}
	\vskip 3. cm
	\caption{Here we depict part of the possible configurations in the case in which the move is horizontal.}
	\label{fig:casoh}
\end{figure}

\item [-] If the moving particle has in $\bar{\h}$ two vertical bonds, that implies the sites $y_1$ and $y_2$ occupied and the site $y_3$ empty (see Figure \ref{fig:casoh} in the middle). By assumptions $g_2'(\bar{\h})=1$ and $g_1'(\bar{\h})=0$, $v(\bar\h)\geq p_1(\bar\h)+p_2(\bar\h)-2=3l_2^*-4$, since in the lines $r_1$, $r_3$ we have overall at least $(p_1(\bar{\h})-1)$ vacancies in $\bar{\h}$. Hence we obtain 
\be{5.8}
\ba{lll}
H(\bar{\h})\geq H(\cR(2l_2^*-k-2,l_2^*+k))+\e(3l_2^*-4)+U_1=\\

U_1l_2^*+2U_2l_2^*+k(U_1-U_2)+U_1-2U_2-2\e(l_2^*)^2-\e kl_2^*+\e k^2+5\e l_2^*+2\e k-4\e>\G\\

\Leftrightarrow \ \e k^2+k[U_1-U_2-\e l_2^*+2\e]-U_2+2\e l_2^*-2\e>0.
\ea
\ee

\noindent
Since $\e\ll U_2$, we get $H(\bar{\h})>\G$. Indeed
\be{}
\ba{ll}
-U_2+2\e l_2^*-2\e=U_2+2\e(\d-1)>U_2-2\e>0,\\

U_1-U_2-\e l_2^*+2\e=U_1-2U_2-\e\d+2\e>\e(1-\d)>0.
\ea
\ee

\item [-] If the moving particle has in $\bar{\h}$ only one horizontal bond, that implies $y_3$ occupied and the sites $y_1$ and $y_2$ empty (see Figure \ref{fig:casoh} on the right hand-side). We observe that $r_6\cap\bar{\h}_{cl}$ consists in the moving particle, otherwise we obtain an horizontal non monotonicity or $g_2'(\bar{\h})\geq2$, which are in contradiction with $g_1'(\bar{\h})=0$ or $g_2'(\bar{\h})=1$ respectively. Thus the line $r_6$ becomes inactive with the move, against $\D s=1$.

\item[-] If the moving particle has in $\bar{\h}$ only a vertical bond, that implies $y_3$ empty and one site among $y_1$ and $y_2$ occupied (see Figure \ref{fig:casov} on the left hand-side). We assume that such a particle is in the site $r_6\cap r_3$ without loss of generality, indeed the argument is analogue if the particle is in the site $r_6\cap r_4$. We observe that the site $s$ and those under it along the line $r_6$ are empty, otherwise we obtain a configuration with an horizontal non monotonicity. The sites next to the left to $y_3$ and next to the right to $x_3$ are empty by the assumptions $g_2'(\bar{\h})=1$ and $g_1'(\bar{\h})=0$. As before, the lines $r_1$ and $r_3$ bring overall at least $p_1(\bar{\h})+p_2(\bar{\h})-2=3l_2^*-4$ vacancies, so by (\ref{5.8}) we get $H(\bar{\h})>\G$.

\ei
The proof is completed in this horizontal case.

\setlength{\unitlength}{1.1pt}
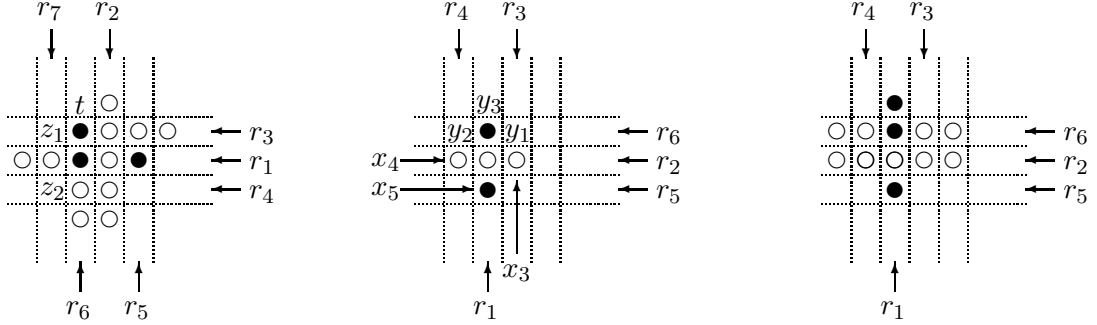
\begin{figure}
	\begin{picture}(380,60)(-30,40)
	\thinlines \qbezier[51](0,0)(35,0)(70,0)
	\qbezier[51](0,10)(35,10)(70,10) \qbezier[51](0,20)(35,20)(70,20)
	\qbezier[51](0,30)(35,30)(70,30)
	\qbezier[51](10,-20)(10,15)(10,50)
	\qbezier[51](20,-20)(20,15)(20,50)
	\qbezier[51](30,-20)(30,15)(30,50)
	\qbezier[51](40,-20)(40,15)(40,50)
	\qbezier[51](50,-20)(50,15)(50,50)
	\put(25,25){\circle*{5}}
	\put(35,15){\circle{5}}
	\put(35,25){\circle{5}}
	\put(45,15){\circle*{5}}
	\put(45,25){\circle{5}}
	\put(55,25){\circle{5}}
	\put(35,35){\circle{5}}
	\put(35,-5){\circle{5}}
	\put(25,15){\circle*{5}}
	\put(23,31){$t$}
	\put(25,5){\circle{5}}
	\put(11,23){$z_1$}
	\put(15,15){\circle{5}}
	\put(25,-5){\circle{5}}
	\put(11,3){$z_2$}
	\put(5,15){\circle{5}}
	\put(35,5){\circle{5}}
	\thinlines
	\put(35,60){\vector(0,-1){10}} \put(30,65){$r_2$}
	\put(15,60){\vector(0,-1){10}} \put(10,65){$r_7$}
	\put(40,-38){$r_5$}\put(45,-30){\vector(0,1){10}}
	\put(25,-30){\vector(0,1){10}} \put(20,-38){$r_6$}
	\put(80,5){\vector(-1,0){10}} \put(83,2){$r_4$}
	\put(80,15){\vector(-1,0){10}} \put(83,12){$r_1$}
	\put(80,25){\vector(-1,0){10}} \put(83,22){$r_3$}
	\thinlines 
	\qbezier[51](140,0)(175,0)(210,0)
	\qbezier[51](140,10)(175,10)(210,10)
	\qbezier[51](140,20)(175,20)(210,20)
	\qbezier[51](140,30)(175,30)(210,30)
	\qbezier[51](160,-20)(160,15)(160,50)
	\qbezier[51](150,-20)(150,15)(150,50)
	\qbezier[51](170,-20)(170,15)(170,50)
	\qbezier[51](180,-20)(180,15)(180,50)
	\qbezier[51](190,-20)(190,15)(190,50) 
	\put(175,15){\circle{5}} 
	\put(165,15){\circle{5}} 
	\put(165,25){\circle*{5}}
	\put(155,15){\circle{5}}
	\put(165,5){\circle*{5}}
	\put(151,23){$y_2$}
	\put(171,23){$y_1$}
	\put(135,5){\vector(1,0){25}}
	\put(125,3){$x_5$}
	\put(135,15){\vector(1,0){15}}
	\put(125,13){$x_4$}
	\put(175,-17){\vector(0,1){25}}
	\put(161,33){$y_3$}
	\put(170,-25){$x_3$}
	\thinlines
	\put(175,60){\vector(0,-1){10}} \put(170,65){$r_3$}
	\put(155,60){\vector(0,-1){10}} \put(150,65){$r_4$}
	\put(165,-30){\vector(0,1){10}} \put(160,-38){$r_1$}
	\put(220,5){\vector(-1,0){10}} \put(223,2){$r_5$}
	\put(220,15){\vector(-1,0){10}} \put(223,12){$r_2$}
	\put(220,25){\vector(-1,0){10}} \put(223,22){$r_6$}
	\thinlines
	\qbezier[51](280,0)(315,0)(350,0)
	\qbezier[51](280,10)(315,10)(350,10)
	\qbezier[51](280,20)(315,20)(350,20)
	\qbezier[51](280,30)(315,30)(350,30)
	\qbezier[51](300,-20)(300,15)(300,50)
	\qbezier[51](290,-20)(290,15)(290,50)
	\qbezier[51](310,-20)(310,15)(310,50)
	\qbezier[51](320,-20)(320,15)(320,50)
	\qbezier[51](330,-20)(330,15)(330,50)
	
	\put(315,15){\circle{5}} 
	\put(295,15){\circle{5}}
	\put(305,15){\circle{5}} 
	\put(315,25){\circle{5}}
	\put(305,35){\circle*{5}}
	
	\put(305,15){\circle{5}} 
	\put(295,15){\circle{5}}
	\put(295,25){\circle{5}}
	\put(285,15){\circle{5}}
	\put(285,25){\circle{5}}
	\put(325,15){\circle{5}}
	\put(325,25){\circle{5}}
	\put(305,5){\circle*{5}} 
	\put(305,25){\circle*{5}} 
	\thinlines
	\put(315,60){\vector(0,-1){10}} \put(310,65){$r_3$}
	\put(300,-38){$r_1$}
	\put(305,-30){\vector(0,1){10}} 
	\put(360,5){\vector(-1,0){10}} \put(363,2){$r_5$}
	\put(360,15){\vector(-1,0){10}} \put(363,12){$r_2$}
	\put(360,25){\vector(-1,0){10}} \put(363,22){$r_6$}
	\put(295,60){\vector(0,-1){10}} \put(290,65){$r_4$}
	
	\end{picture}
	\vskip 3. cm
	\caption{Here we depict on the left hand-side the remaining case for the horizontal line; in the middle and on the right hand-side we depict part of the configuration if the move is vertical.}
	\label{fig:casov}
\end{figure}

\medskip
\noindent
{\bf Case (v).} Suppose now that the move is vertical, i.e., $r_1$ is a vertical line. Since the unique line that must become active is $r_2$, at least one site among $x_3$, $x_4$ and $x_5$ must be occupied. As in the case (h), we deduce that $x_5$ is occupied, and $x_3$ and $x_4$ are empty (see Figure \ref{fig:casov} in the middle). We deduce $g_1'(\bar{\h})\geq1$.

\bi
\item [-] If the moving particle has in $\bar{\h}$ at least one vertical and one horizontal bond, that implies $y_3$ occupied and one site among $y_1$ and $y_2$ occupied (see Figure \ref{fig:casov} in the middle). Thus $\D H=H(\h)-H(\bar{\h})\geq U_1$, since in the move two horizontal bonds are lost (because $x_3$ and $x_4$ are empty) and the vertical bond is recovered with the one with $x_5$. Thus by (\ref{5.7}) we get
\be{}
H(\h)\geq H(\bar{\h})+U_1\geq H(\cR(2l_2^*-k-2,l_2^*+k))+\e(l_2^*+k-1)+U_1+U_2>\G.
\ee

\item [-] If the moving particle has in $\bar{\h}$ at least two horizontal bonds connecting it to other particles, then $\D H=H(\h)-H(\bar{\h})\geq -U_2+2U_1$, since we lose two horizontal bonds because $x_3$ and $x_4$ are empty, and we recover a vertical bond with the particle in the site $x_5$. Thus by (\ref{5.7}) we get
\be{}
H(\h)\geq H(\bar{\h})+2U_1-U_2\geq H(\cR(2l_2^*-k-2,l_2^*+k))+\e(l_2^*+k-1)+2U_1>\G.
\ee

\item [-] If the moving particle has in $\bar{\h}$ only one vertical bond connecting it to other particles, the site $y_3$ must be occupied and the sites $y_1$ and $y_2$ must be empty. We observe that $r_6\cap\bar\h_{cl}$ consists in the moving particle, indeed either $n(\bar\h)\geq1$ or $\bar\h$ is not monotone, which is in contradiction with $g_2'(\bar{\h})=0$. The arguments used are similar to those used for the corresponding in case (h). In this situation the line $r_6$ becomes inactive, which contradicts $\D s=1$.

\item [-] If the moving particle has in $\bar{\h}$ only one horizontal bond, that implies $y_3$ empty and one site among $y_1$ and $y_2$occupied. Thus $\D H=H(\h)-H(\bar{\h})=U_1-U_2$, since we lose a vertical bond and we recover an horizontal one.
\bi
\item [$\bullet$] If $g_1'(\bar{\h})=1$, since $v(\bar{\h})\geq p_1(\bar{\h})+p_2(\bar{\h})-2$ (analogue reasoning used for the corresponding in the case (h)), by (\ref{5.8}) we get
\be{}
H(\h)=H(\bar{\h})+U_1-U_2\geq H(\cR(2l_2^*-k-2,l_2^*+k))+\e(3l_2^*-4)+U_1>\G.
\ee

\item [$\bullet$] If $g_2'(\bar{\h})\geq2$, by (\ref{5.7}) we get
\be{}
H(\h)\geq H(\bar{\h})+U_1-U_2\geq H(\cR(2l_2^*-k-2,l_2^*+k))+\e(l_2^*+k-1)+U_1-U_2+2U_2>\G.
\ee

\ei
\ei
The proof is completed in this vertical case, so it is concluded for the case $\D s=1$.

\epr

\subsection{Proof of Proposition \lowercase{\ref{lemma8}}}
\label{dim8}
\bpr
Let $\D s=2$; by Remark \ref{remark2} we consider $s(\bar{\h})\geq s^*-3$. We distinguish the following cases:

\bi
\item [(a)] $s(\bar\h)=s^*-3$;
\item [(b)] $s(\bar\h)=s^*-2$;
\item [(c)] $s(\bar\h)=s^*-1$;
\item [(d)] $s(\bar\h)\geq s^*$.
\ei

First, we consider the case (a). If $s(\bar{\h})=s^*-3$, then $s(\bar{\h})=s^*-1$. If $p_2(\h)\leq l_2^*-1$, we get $\h\in\cB$. If $p_2(\h)\geq l_2^*$, by Lemma \ref{lemma3.2}{\it (ii)} we have 

$$v(\h)\geq2p_{min}(\h)-5\geq p_{min}(\h)-1\Leftrightarrow \ p_{min}(\h)\geq4.$$

\noindent
Since $p_{min}(\h)\geq4,$ in this case it is impossible to leave $\cB$.

By Lemma \ref{lemma3.2}{\it (ii)} we know that $n(\bar{\h})\geq1$. In the case (c) we conclude as in the case ($\D$ s=0.a-i) and we refer to the Appendix for the explicit computations.

In the case (d), we get $H(\bar{\h})>\G$ as in the case ($\D$ s=-1.Case I) and we refer to the Appendix for the explicit computations.

It remains to analyze the case (b). If $s(\bar{\h})=s^*-2$ and $n(\bar{\h})\geq2$, by (\ref{5*}) and $k\geq0$, we get $H(\bar{\h})>\G$. Indeed
\be{}
\ba{ll}
H(\bar{\h})\geq H(\cR(2l_2^*-k-3,l_2^*+k))+2\D>\G\Leftrightarrow\\

\Leftrightarrow \ \e k^2+k[U_1-U_2-\e l_2^*+3\e]+U_1>0.
\ea
\ee

\noindent
It remains to analyze the case $s(\bar{\h})=s^*-2$ and $n(\bar{\h})=1$. If the unique free particle is the moving particle, we can not have $\D s=2$, indeed the lines $r_3$, $r_4$ and $r_5$ can not be activated and in order to have $\D s=2$ the lines $r_1$ and $r_2$ must become active. This implies that the sites $x_3$, $x_4$ and $x_5$ must be empty, but then $\D s=0$, which contradicts $\D s=2$.

If $g_1'(\bar{\h})+g_2'(\bar{\h})\geq1$, by (\ref{5*}) and $k\geq0$, we get $H(\bar{\h})>\G$. Indeed
\be{5.9}
\ba{ll}
H(\bar{\h})\geq H(\cR(2l_2^*-k-3,l_2^*+k))+U_2+\D>\G\Leftrightarrow\\

\Leftrightarrow \ \e k^2+k[U_1-U_2-\e l_2^*+3\e]+\e>0.
\ea
\ee

\noindent
It remains only to analyze the case $s(\bar{\h})=s^*-2$, $n(\bar{\h})\geq1$ and $g_1'(\bar{\h})=g_2'(\bar{\h})=0$. We disinguish two cases:

\bi
\item[1)]
the free particle is in the site $x_i$, with $i\in\{3,4\}$, and the lines that become active are $r_2$ and $r_i$. Due to $g'_1(\bar\h)=g'_2(\bar\h)=0$, the site $x_5$ is empty and the site $\{x_3, x_4\}\backslash \{x_i\}$ is empty.

\item[2)]
the free particle is in the site $x_5$, the lines that become active are $r_2$ and $r_5$ and the sites $x_3$ and $x_4$ are empty.
\ei

\setlength{\unitlength}{1.3pt}
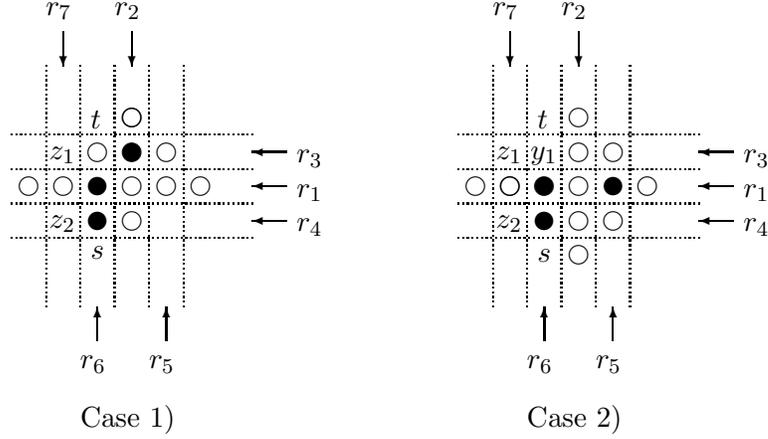
\begin{figure}
	\begin{picture}(380,60)(-30,40)
	\thinlines 
	\qbezier[51](50,0)(85,0)(120,0)
	\qbezier[51](50,10)(85,10)(120,10) \qbezier[51](50,20)(85,20)(120,20)
	\qbezier[51](50,30)(85,30)(120,30)
	\qbezier[51](60,-20)(60,15)(60,50)
	\qbezier[51](70,-20)(70,15)(70,50)
	\qbezier[51](80,-20)(80,15)(80,50)
	\qbezier[51](90,-20)(90,15)(90,50)
	\qbezier[51](100,-20)(100,15)(100,50)
	\put(75,25){\circle{5}}
	\put(85,15){\circle{5}}
	\put(85,25){\circle{5}}
	\put(95,15){\circle{5}}
	\put(95,25){\circle{5}}
	\put(85,25){\circle*{5}}
	\put(85,35){\circle{5}}
	\put(105,15){\circle{5}}
	\put(85,35){\circle{5}}
	\put(75,15){\circle*{5}}
	\put(73,32){$t$}
	\put(75,5){\circle*{5}}
	\put(61,23){$z_1$}
	\put(65,15){\circle{5}}
	\put(73,-6){$s$}
	\put(61,3){$z_2$}
	\put(55,15){\circle{5}}
	\put(85,5){\circle{5}}
	\thinlines
	\put(85,60){\vector(0,-1){10}} \put(80,65){$r_2$}
	\put(65,60){\vector(0,-1){10}} \put(60,65){$r_7$}
	\put(90,-38){$r_5$} \put(95,-30){\vector(0,1){10}}
	\put(75,-30){\vector(0,1){10}} \put(70,-38){$r_6$}
	\put(130,5){\vector(-1,0){10}} \put(133,2){$r_4$}
	\put(130,15){\vector(-1,0){10}} \put(133,12){$r_1$}
	\put(130,25){\vector(-1,0){10}} \put(133,22){$r_3$}
	\put(70,-55){$\rm \hbox{Case 1)}$}
	
	\thinlines
	\qbezier[51](180,0)(215,0)(250,0)
	\qbezier[51](180,10)(215,10)(250,10)
	\qbezier[51](180,20)(215,20)(250,20)
	\qbezier[51](180,30)(215,30)(250,30)
	\qbezier[51](200,-20)(200,15)(200,50)
	\qbezier[51](190,-20)(190,15)(190,50)
	\qbezier[51](210,-20)(210,15)(210,50)
	\qbezier[51](220,-20)(220,15)(220,50)
	\qbezier[51](230,-20)(230,15)(230,50)
	
	\put(215,15){\circle{5}} 
	\put(195,15){\circle{5}}
	\put(205,15){\circle{5}} 
	\put(215,25){\circle{5}}
	\put(203,32){$t$}
	\put(203,-7){$s$}
	
	\put(205,15){\circle*{5}} 
	\put(195,15){\circle{5}}
	\put(185,15){\circle{5}}
	\put(225,15){\circle*{5}}
	\put(225,5){\circle{5}}
	\put(215,5){\circle{5}}
	\put(215,-5){\circle{5}}
	\put(215,35){\circle{5}}
	\put(235,15){\circle{5}}
	\put(225,25){\circle{5}}
	\put(205,5){\circle*{5}} 
	\put(201,23){$y_1$} 
	\put(191,23){$z_1$} 
	\put(191,3){$z_2$}
	\thinlines
	\put(215,60){\vector(0,-1){10}} \put(210,65){$r_2$}
	\put(200,-38){$r_6$}
	\put(205,-30){\vector(0,1){10}} 
	\put(260,5){\vector(-1,0){10}} \put(263,2){$r_4$}
	\put(260,15){\vector(-1,0){10}} \put(263,12){$r_1$}
	\put(260,25){\vector(-1,0){10}} \put(263,22){$r_3$}
	\put(195,60){\vector(0,-1){10}} \put(190,65){$r_7$}
	\put(220,-38){$r_5$} \put(225,-30){\vector(0,1){10}}
	
	\put(200,-55){$\rm \hbox{Case 2)}$}
	\end{picture}
	\vskip 4.5 cm
	\caption{Here we depict part of the configuration if $\D s=2$ for the cases 1) and 2).}
	\label{fig:casos=2}
\end{figure}

In both cases if the moving particle has in $\bar\h$ at least one vertical and one horizontal bond connecting it to other particles, by (\ref{5.9}) we get $H(\h)>\G$. Indeed

$$\D H\geq U_2\Rightarrow \ H(\h)\ge H(\bar\h)+U_2\geq H(\cR(2l_2^*-k-3,l_2^*+k))+\D+U_2>\G.$$

If the moving particle has in $\bar\h$ either two bonds orthogonal to the move, or only one vertical or only one
horizontal bond connecting it to other particles, then it is
impossible to leave $\cB$, indeed in this case there exists a line $r$ ($r=r_1$ or $r=r_6$) such that its intersection with $\bar\h_{cl}$ is only the moving particle. If $r=r_6$, this line becomes inactive after the move, which is in contradiction with $\D s=2$. If $r=r_1$, we analyze in detail the possible configuration for $\bar{\h}$ in both cases 1) and 2). First, we consider the case 1). We assume $i=3$ without loss of generality. Thus in the site $x_3$ there is a free particle in $\bar\h$. We have that the site $y_2$ must be occupied, otherwise $n(\bar{\h})\geq2$. The site $x_4$ must be empty, since the line $r_2$ must become active with the move, so it must be inactive in $\bar\h$. The sites $x_5$, $y_3$ and $z_3$ are empty, since we are in the case $r_1\cap\bar\h_{cl}$ consists in the moving particle (see Figure \ref{fig:casos=2} on the left hand-side). At least one site among $z_2$ and $s$ must be occupied, otherwise the line $r_6$ becomes inactive after the move, which contradicts $\D s=2$. Thus in $\bar\h$ there is a monotone cluster (eventually a finite non connected union of monotone clusters) attached to the moving particle: suppose that it has $m\geq0$ vacancies. Since $p_2(\h)=p_2(\bar{\h})+1$ and $p_2(\bar{\h})\geq4$, we get $\h\in\cB$. Indeed

$$v(\h)\geq m+(p_2(\h)-2)+(p_1(\h)-1)\geq 3+p_1(\h)-1=p_1(\h)+2\geq p_1(\h)-1=p_{max}(\h)-1.$$

\noindent
In the case 2), we have that in the site $x_5$ there is a free particle. Furtermore, all the sites along the line $r_2$ are empty, otherwise either we have $n(\bar{\h})\geq2$ or $r_2$ is already active in $\bar{\h}$. The moving particle can not be free, so it has at least one vertical bond. Thus at least one site among $y_1$ and $y_2$ must be occupied: without loss of generality we assume $y_2$ occupied (see Figure \ref{fig:casos=2} on the right hand-side). Similarly to the case 1), since also in this case $r_1\cap\bar\h_{cl}$ consists in the moving particle, we deduce that the sites $y_3$ and $z_3$ are empty. Furthermore, in $\bar\h$ there is a monotone cluster (or a finite union of clusters) attached to the moving particle: again we suppose that it has $m\geq0$ vacancies. Since $p_2(\h)=p_2(\bar{\h})\geq4$, we get $\h\in\cB$. Indeed

$$v(\h)\geq m+(p_1(\h)-2)+(p_1(\h)-1)\geq p_1(\h)+1\geq p_1(\h)-1=p_{max}(\h)-1.$$

\epr

\subsection{Proof of Lemmas}
\label{dimlemmi}
In this subsection we report the proof of Lemmas given in subsection \ref{lemmi}.

\bigskip
\begin{proof*}{\bf of Lemma \ref{decrescita}}

	\noindent
	Let $r$ be a line that becomes inactive with the move. Since $r\cap\bar\h_{cl}\neq\emptyset$ and $r\cap\h_{cl}=\emptyset$, we note that only free particles could be present along the line $r$ in $\h$. Recalling (\ref{proie1}) for the definition of the projections $p_1$ and $p_2$, we note that they depends only on the clusterized part of the configuration, thus in $\h$ the line $r$ does not contribute to $p_1(\h)$ and $p_2(\h)$. This implies that $p_2(\h)=p_2(\bar\h)-1$ if $r$ is an horizontal line and $p_1(\h)=p_1(\bar\h)-1$ if $r$ is vertical.
	
\end{proof*}

\medskip
\begin{proof*}{\bf of Lemma \ref{lemma1}.}

	\noindent
	First we analyze the line $r_2$. We argue by contradiction: suppose that $r_2$ becomes inactive with the move. Thus the moved particle must be free in $\h$: this implies that the sites $x_3$, $x_4$ and $x_5$ must be empty. Since $r_2$ must be active in $\bar{\h}$, at least one particle above the site $x_3$ or under $x_4$ must be in $\bar{\h}_{cl}$. If $|r_2\cap\bar\h_{cl}|=1$, we indicate with $a$ the site occupied by such particle (represented on the left hand-side in Figure \ref{fig:r25}) and we can suppose without loss of generality that such particle is above $x_3$. Let $a_1$, $a_2$ and $a_3$ the nearest neighbors of that particle as in Figure \ref{fig:r25}. Since $r_2$ must be active in $\bar\h$, we have necessarily that at least one among $a_1$, $a_2$ and $a_3$ must contain a particle. The move does not involve the site $a$ and its nearest neighbors, so we conclude that it is not possible that the line $r_2$ becomes inactive. If $|r_2\cap\bar{\h}_{cl}|\geq2$, again the line $r_2$ does not become inactive. Thus we have proved that $r_2$ also remains inactive in $\h$ and thus it can not become inactive.

	\setlength{\unitlength}{1.1pt}
	\begin{figure}
		\begin{picture}(380,60)(-30,40)
		\thinlines 
		\qbezier[6](75,25)(82.5,25)(90,25)
		\qbezier[12](75,10)(90,10)(105,10)
		\qbezier[6](105,-5)(105,2.5)(105,10)
		\qbezier[12](75,-5)(90,-5)(105,-5)
		\qbezier[6](75,-20)(82.5,-20)(90,-20)
		\qbezier[18](75,10)(75,35)(75,60)
		\qbezier[6](75,-20)(75,-13.5)(75,-5)
		\qbezier[30](90,-20)(90,20)(90,60) 
		\qbezier[6](75,45)(82.5,45)(90,45)
		\qbezier[6](90,60)(97.5,60)(105,60)
		\qbezier[6](105,60)(105,67.5)(105,75)
		\qbezier[6](90,75)(97.5,75)(105,75)
		\qbezier[6](90,75)(90,82.5)(90,90)
		\qbezier[6](75,75)(75,82.5)(75,90)
		\qbezier[6](75,90)(82.5,90)(90,90)
		\qbezier[6](60,60)(67.5,60)(75,60)
		\qbezier[6](60,60)(60,67.5)(60,75)
		\qbezier[6](60,75)(67.5,75)(75,75)
		
		\put(77.5,80){$a_1$}
		\put(92.5,65){$a_2$}
		\put(62.5,65){$a_3$}
		\put(105,90){\vector(-1,-1){13}}
		\put(105,90){$a$}
		\put(111,17.5){\vector(-1,0){20}}
		\put(111,15){$x_3$}
		\put(60,-5){\line(1,0){15}}
		\put(60,-5){\line(0,1){15}}
		\put(75,-5){\line(0,1){15}}
		\put(75,10){\line(-1,0){15}}
		\put(82.5,2.5){\circle{7.5}}
		\put(82.5,-12.5){\circle{7.5}} 
		\put(82.5,17.5){\circle{7.5}}
		\put(67.5,2.5){\circle*{7.5}}
		\put(97.5,2.5){\circle{7.5}}
		\put(82.5,52.5){\circle{7.5}}
		\put(82.5,67.5){\circle*{7.5}}
		\put(90,60){\line(-1,0){15}}
		\put(75,60){\line(0,1){15}}
		\put(75,75){\line(1,0){15}}
		\put(90,75){\line(0,-1){15}}
		\put(82.5,-32.5){\vector(0,1){10}} 
		\put(77.5,-38.5){$r_2$}

		\thinlines 
		\qbezier[6](290,25)(297.5,25)(305,25)
		\qbezier[12](275,10)(290,10)(305,10)
		\qbezier[6](305,-5)(305,2.5)(305,10)
		\qbezier[12](275,-5)(290,-5)(305,-5)
		\qbezier[6](275,-20)(282.5,-20)(290,-20)
		\qbezier[18](305,10)(305,35)(305,60)
		\qbezier[6](275,-20)(275,-13.5)(275,-5)
		\qbezier[30](290,-20)(290,20)(290,60) 
		\qbezier[6](290,45)(297.5,45)(305,45)
		\qbezier[6](305,60)(312.5,60)(320,60)
		\qbezier[6](320,60)(320,67.5)(320,75)
		\qbezier[6](305,75)(312.5,75)(320,75)
		\qbezier[6](305,75)(305,82.5)(305,90)
		\qbezier[6](290,75)(290,82.5)(290,90)
		\qbezier[6](290,90)(297.5,90)(305,90)
		\qbezier[6](275,60)(282.5,60)(290,60)
		\qbezier[6](275,60)(275,67.5)(275,75)
		\qbezier[6](275,75)(282.5,75)(290,75)
		\qbezier[6](275,10)(275,17.5)(275,25)
		\qbezier[6](275,25)(282.5,25)(290,25)
		
		\put(292.5,80){$a_1$}
		\put(277.5,65){$a_3$}
		\put(307.5,65){$a_2$}
		\put(320,90){\vector(-1,-1){13}}
		\put(320,90){$a$}
		\put(260,-5){\line(1,0){15}}
		\put(260,-5){\line(0,1){15}}
		\put(275,-5){\line(0,1){15}}
		\put(275,10){\line(-1,0){15}}
		\put(282.5,2.5){\circle{7.5}}
		\put(278,-15){$x_4$} 
		\put(278,15){$x_3$}
		\put(297.5,17.5){\circle{7.5}}
		\put(267.5,2.5){\circle*{7.5}}
		\put(297.5,2.5){\circle{7.5}}
		\put(297.5,52.5){\circle{7.5}}
		\put(297.5,67.5){\circle*{7.5}}
		\put(305,60){\line(-1,0){15}}
		\put(290,60){\line(0,1){15}}
		\put(290,75){\line(1,0){15}}
		\put(305,75){\line(0,-1){15}}
		\put(297.5,-32.5){\vector(0,1){10}} 
		\put(292.5,-38.5){$r_5$}
		\end{picture}
		\vskip 3.5 cm
		\caption{Possible representation of $\bar{\h}$ in the case that $r_2$ (on the left hand-side) and $r_5$ (on the right hand-side) could become inactive.}
		\label{fig:r25}
		
	\end{figure}
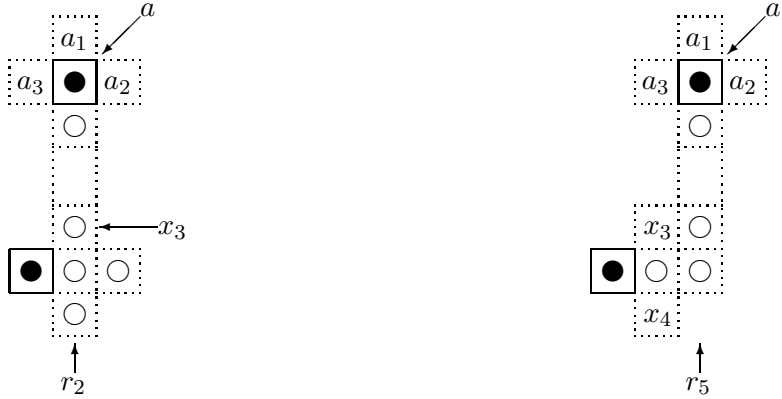

	If we consider line $r_5$, again we argue by contradiction: suppose that $r_5$ becomes inactive with the move. Thus the site $x_5$ must be empty, otherwise the line $r_5$ can not become inactive. First, we consider the case $|r_5\cap\bar{\h}_{cl}|=1$: suppose without loss of generality that such a particle is above the site $x_5$ and call it $a$ (see Figure \ref{fig:r25} on the right hand-side). Let $a_1$, $a_2$ and $a_3$ be the nearest neighbors of the particle in $a$. Since $r_5$ must be active in $\bar\h$, we have necessarily that at least one site among $a_1$, $a_2$ and $a_3$ must be occupied. The move does not involve the particle in $a$ and its neighbors, so we are able to conclude that the line $r_5$ can not become inactive. If $|r_5\cap\bar{\h}_{cl}|\geq2$, again the line $r_5$ does not become inactive. Thus we have proved that $r_5$ also remains inactive in $\h$ and thus it can not become inactive.

	Now we focus on lines $r_6$ and $r_7$: we have to prove that these lines can not become active. If we consider line $r_7$, we suppose by contradiction that $r_7$ becomes active with the move. Thus we have that line $r_7$ must be inactive in $\bar\h$, so the site $y_3$ must be empty. With the move of the particle from $x_1$ to $x_2$ the number of neighboring particles of $y_3$ decreases. Thus $r_7$ can not be active in $\h$. 
	
	In a similar way we suppose by contradiction that the line $r_6$ becomes active: this means that line $r_6$ is inactive in $\bar{\h}$, so the moving particle is free in $\bar\h$. Thus $r_6$ can not become active after the move. 
\end{proof*}

\medskip
\begin{proof*}{\bf of Lemma \ref{lemma0}.}

	\noindent
	First, we consider the point (i). If the line $r_1$ becomes inactive we have that the sites $x_3$, $x_4$ and $x_5$ must be empty. Thus the lines $r_2$, $r_3$, $r_4$ and $r_5$ can not become active. Furthermore, by Lemma \ref{lemma1} the lines $r_6$ and $r_7$ can not become active, thus no line can become active.
	
	Now we analyze point (ii). If the line $r_1$ becomes active we have that in the site $x_1$ there is a free particle in $\bar\h$ ($y_1$, $y_2$ and $y_3$ must be empty in $\bar\h$). Thus the line $r_6$ can not become inactive. Furthermore, by Lemma \ref{lemma1} the lines $r_2$ and $r_5$ can not become inactive. In order to have $r_7$ inactive in $\h$, the site $y_3$ must be occupied, but this contradicts the fact that in $x_1$ there is a free particle in $\bar\h$. Thus the only lines that can become inactive in this situation are $r_3$ and $r_4$. 
\end{proof*}

\medskip
\begin{proof*}{\bf of Lemma \ref{ulteriore}.}

	\noindent
	First, we suppose that the line becoming active is $r_3$ (respectively $r_4$). We consider now the case $\D s=-2$. By Remark \ref{remark2} we know that $s(\bar\h)\geq s^*+1$, with $\bar\h\in\cB$. By (\ref{defB}) for the case $s\geq s^*$, we get $p_2(\bar\h)=l_2^*$ and thus $p_1(\bar\h)=p_{max}(\bar\h)$. Since $l_2^*-5\leq p_2(\h)\leq l_2^*+5$, also for the configuration $\h$ we deduce that $p_1(\h)>p_2(\h)$. The line $r_3$ (resp. $r_4$) becomes active with the move bringing $p_1(\h)-1$ vacancies in $\h_{cl}$, since the unique particle along the line $r_3$ (resp. $r_4$) in $\h_{cl}$ is in the site $x_3$ (resp. $x_4$), otherwise the line $r_3$ (resp. $r_4$) is already active in $\bar\h$. We analyze separately the cases (i) and (ii). If $s(\h)=s^*-1$, since $p_1(\h)\geq p_{min}(\h)$ and $p_2(\h)\geq l_2^*$, we deduce that $v(\h)\geq p_{min}(\h)-1$: this implies that $\h\in\cB$. If $s(\h)\geq s^*$, since $p_2(\h)=l_2^*$ by assumption and $v(\h)\geq p_{1}(\h)-1=p_{max}(\h)-1$, we get $\h\in\cB$.
	
	Now we analyze the case $-1\leq\D s\leq5$. For each value of $\D s$, by Remark \ref{remark2} we know that $s(\h)\geq s^*-1$. Again the line $r_3$ (resp. $r_4$) becomes active with the move bringing $p_1(\h)-1$ vacancies in $\h$. We analyze separately the cases (i) and (ii). If $s(\h)=s^*-1$, since $p_1(\h)\geq p_{min}(\h)$ and $p_2(\h)\geq l_2^*$, we deduce that $v(\h)\geq p_{min}(\h)-1$: this implies that $\h\in\cB$. If $s(\h)\geq s^*$, since $p_2(\h)=l_2^*$ and thus $p_1(\h)>p_2(\h)$, we deduce that $p_1(\h)=p_{max}(\h)$, so we get $v(\h)\geq p_{max}(\h)-1$. Thus we obtain that $\h\in\cB$.
	
	Now we suppose that the line becoming active is $r_1$ and the site $x_5$ is empty. Since $r_1$ must become active, we know that in the site $x_1$ there is a free particle in $\bar\h$. The line $r_1$ becomes active with the move bringing $p_1(\h)-1$ vacancies in $\h$, since the unique particle along the line $r_1$ in $\h_{cl}$ is in the site $x_1$, otherwise the line $r_1$ is already active in $\bar\h$, since the site $x_5$ is empty by assumption. The proof proceed from now on in the same way as in the case in which $r_3$ (resp. $r_4$) is the horizontal line becoming active.
	
\end{proof*}

\medskip
\begin{proof*}{\bf of Lemma \ref{ult2}}

	\noindent
	First, we consider the case {\it (i)}. By assumption the line $r_1$ does not become active and by Lemma \ref{lemma0}{\it (i)} we deduce that $r_1$ does not become inactive, otherwise no line can become active. Since only one line among $r_3$ and $r_4$ becomes active, we get $p_2(\h)=p_2(\bar\h)+1$. If $p_2(\bar\h)\leq l_2^*-2$ then $p_2(\h)\leq l_2^*-1$, so we get $\h\in\cB$. If $p_2(\bar\h)=l_2^*-1$ then $p_2(\h)=l_2^*$, thus by Lemma \ref{ulteriore}{\it (i),(ii)} we conclude that $\h\in\cB$.
	
	Now we consider the case {\it (ii)}. Since only two horizontal lines become active, we get $p_2(\h)=p_2(\bar\h)+2$. If $p_2(\bar\h)\leq l_2^*-3$, we get $p_2(\h)\leq l_2^*-1$, thus $\h\in\cB$. If $p_2(\bar\h)=l_2^*-2$ then $p_2(\h)=l_2^*$, thus by Lemma \ref{ulteriore}{\it(i),(ii)} we conclude that $\h\in\cB$. If $p_2(\bar\h)=l_2^*-1$ and $s(\h)=s^*-1$, by Lemma \ref{ulteriore}{\it (i)} we deduce that $\h\in\cB$.
	
\end{proof*}

\section{Reduction}
\label{recurrence}

\subsection{Recurrence property}
\label{reduction}
The goal of this subsection is to prove the Proposition \ref{V^*}. In order to prove this, we adopt a strategy explained in \cite{MNOS}, the so-called {\it reduction}.

As an application of this technique, we apply \cite[Theorem 3.1]{MNOS} to $V=V^*$. Thus we obtain

\be{supV^*}
\b\longmapsto\sup_{\h\in\cX}\P\Bigl(\t^\h_{\cX_{V^*}}>e^{\b(V^*+\e)}\Bigl)\hbox{ is } \SES, \ee

for any $\e>0$ and for any $\b$ sufficiently large.

\noindent
Since in this case $\cX_{V^*}\subseteq\{\vuoto,\pieno\}$, this result states that the probability that the first hitting time to the set $\{\vuoto,\pieno\}$ from any state $\h\in\cX$ is arbitrarily large is super-exponentially small.

An analogous result can be derived with $V=\G$ and so we have that $\cX_{\G}=\cX^m=\{\vuoto\}$.

\medskip 
For the sequel we need some geometrical definitions. Let $\eta\in \cX$ given.
\bd{} A site $x\in{\Lambda}$ is \emph{connected trough empty (resp.\ full) sites} to $\partial^{-}\Lambda$ if there exists $x_1, \ \ldots, \ x_n$ a connected chain of nearest-neighbor empty (resp.\ full) sites, i.e., $x_1\in nn(x)$, $x_2\in nn(x_1)$, $\ldots$, $x_n\in nn(x_{n-1})$, $x_n\in \partial^{-}\Lambda$ and $\eta(x_1)=\eta(x_2)=\cdots=\eta(x_n)=0$ ($\eta(x_1)=\eta(x_2)=\cdots=\eta(x_n)=1$).
\ed
\bd{} An \emph{external corner} of a set $A\subset\Lambda$ is a site $x\notin A$ such that 
$$\sum_{y\in nn(x):(x,y)\in \Lambda_{0,h}^{*}} \chi_{A}(y)=1 \quad \mbox{and }\sum_{y\in nn(x):(x,y)\in \Lambda_{0,v}^{*}} \chi_{A}(y)=1,$$
where $\chi_{A}$ denotes the characteristic function of the set $A$.
\ed

\bd{} An \emph{internal corner} of a set $A\subset \Lambda$ is a site $x\in A$ such that 
$$\sum_{y\in nn(x):(x,y)\in \Lambda_{0,h}^{*}} \chi_{A}(y)=1 \quad \mbox{and } \sum_{y\in nn(x):(x,y)\in \Lambda_{0,v}^{*}} \chi_{A}(y)=1.$$
\ed

\medskip
\begin{figure}[h!]
	\begin{center}
		\includegraphics[width=14cm]{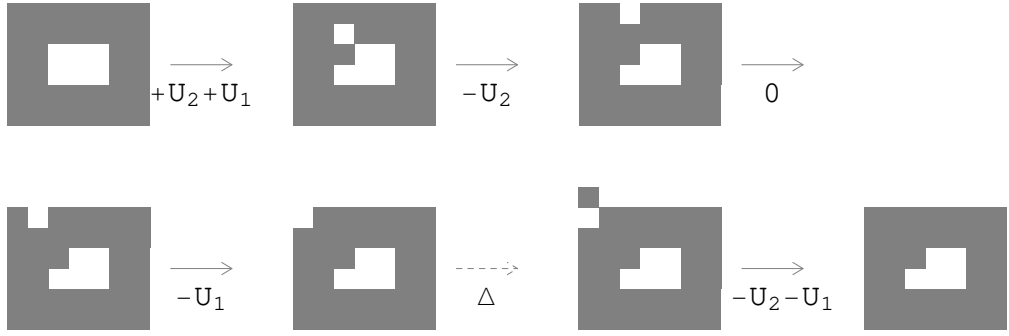}\\
		\caption{$(U_1+U_2)$-reduction of a rectangle with a hole.}
		\label{fig:reduction}
	\end{center}
\end{figure}

\noindent
Let $\eta^{ext}$ be the set of sites $x\in \Lambda_0$ such that $\eta(x)=1$ and $x$ is connected trough empty sites with $\partial^{-}\Lambda$.

\bp{U_1+U_2}\cite[Proposition 16]{NOS} Any configuration $\eta\in{\mathcal{X}_{U_1+U_2}}$ has no free particles and it has only rectangular clusters with minimal side larger than one.
\ep

\begin{proof*}{\bf of Proposition \ref{V^*}}
	
	Suppose that $\eta\in{\mathcal{X}_{U_1+U_2}}$ and $\eta\neq \underline{0}$, $\underline{1}$, so from the previous proposition $\eta$ has only rectangular clusters which are connected through empty sites to $\partial^{-}\Lambda$, i.e., $\eta^{ext}=\partial^{-}\eta$. \\
	Suppose now that a rectangular cluster of $\eta$ has a vertical subcritical side, i.e., $l_2\leq l_2^*-1$, then it is possible to reduce $\eta$ with the path described in Figure \ref{fig:growc} that removes a column of length $l_2$ with energy barrier $\Delta H$(remove column)$=U_1+U_2+\varepsilon(l_2-2)<2\Delta-U_1=\Delta H$(add column) (this is true if and only if we are in the subcritical side). Otherwise if any rectangle in $\eta$ has vertical supercritical sides ($l_2\geq l_2^{*}$), it is possible to reduce $\eta$ with the path described in Figure \ref{fig:growc} that adds a column with energy barrier $\Delta H$(add column)$=2\Delta-U_1$. Since $2\D-U_1<\D-U_2+U_1$ in the strongly anisotropic case, the proof is complete defining 
	$$V^*:=\max{\{U_1+U_2, \ U_1+U_2+\varepsilon(l_2-2), \ 2\Delta-U_1\}}=2\D-U_1<\Gamma.$$
\end{proof*}

We remark that Propositions \ref{U_1+U_2} and \ref{V^*} state the following inequalities:

$$0\leq U_1+U_2\leq V^*<\G.$$

Furthermore, note that from (\ref{supV^*}) and Proposition \ref{V^*} we obtain that from any configuration in $\cX$ the Kawasaki dynamics hits $\vuoto$ or $\pieno$ with an overwhelming probability in a time much less than the transition time. An analogous result is obtained in \cite{HNOS} for three-dimensional Kawasaki dynamics.

\subsection{Proof of Theorem \lowercase{\ref{t2}}}
\label{dimt2}
To prove Theorem \ref{t2} we need Proposition \cite[Theorem 3.2]{BNZ}, that represents the main property of cycles: with large probability every state in a cycle is visited by the process before the exit.

Using this result, to prove Theorem \ref{t2} it is sufficient to show the following:
\begin{enumerate}
	\item if $\eta$ is a rectangular configuration contained in $\mathcal{R}(2l_2^*-3,l_2^*-1)$, then there exists a cycle $\mathcal{C}_{\underline{0}}$ containing $\eta$ and $\underline{0}$ and not containing $\underline{1}$;
	\item if $\eta$ is a rectangular configuration containing $\mathcal{R}(2l_2^*-2,l_2^*)$, then there exists a cycle $\mathcal{C}_{\underline{1}}$ containing $\eta$ and $\underline{1}$ and not containing $\underline{0}$.
\end{enumerate}
We start by showing $1$. Let $\mathcal{C}_{\underline{0}}$ be the maximal connected set containing $\underline{0}$ such that

\noindent
$\max_{\eta'\in{\mathcal{C}_{\underline{0}}}}H(\eta')<\Gamma$. By definition $\mathcal{C}_{\underline{0}}$ is a cycle containing $\underline{0}$. It does not contain $\underline{1}$ since $\Phi(\underline{0}, \underline{1})=\Gamma$. We have only to prove that $\mathcal{C}_{\underline{0}}$ contains $\eta$. This can be easily obtained by constructing a path $\omega^{\eta,\underline{0}}$ going from $\eta$ to $\underline{0}$ keeping the energy less than $\Gamma$. $\omega^{\eta, \ \underline{0}}$ is obtained by erasing site by site, each column of $\eta$ and by showing that all the configurations of this path are in $\mathcal{C}_{\underline{0}}$. \\
More precisely, let $\eta=\{(x,y)\in{\mathbb{Z}^2}: \ x\in(n, n+l_1], \ y\in(m,m+l_2]\}\in{\mathcal{R}(l_1,l_2)}$ for some $n, \ m\in{\mathbb{Z}}$. Let $\{\bar{\omega}_i^{\eta, \underline{0}}\}_{i=0,..,l_1}$ be a path of rectangular configurations, starting from $\eta$ and ending in $\underline{0}$, given by 
\begin{equation}
\bar{\omega}_i^{\eta, \underline{0}}=\{(x, y):\ x\in(n, n+l_1-i], \ y\in(m, m+l_2]\}.
\end{equation}
To complete the construction we can use now the same idea applied in the definition of the reference path $\omega^*$: between every pair $\bar{\omega}_i^{\eta,\underline{0}}$, $\bar{\omega}_{i+1}^{\eta,\underline{0}}$ we can insert a sequence $\tilde{\omega}_{i,0}^{\eta,\underline{0}},..,\tilde{\omega}_{i,l_2}^{\eta,\underline{0}}$ such that $\tilde{\omega}_{i,0}^{\eta,\underline{0}}=\bar{\omega}_i^{\eta,\underline{0}}$ and for $j>0$, $\tilde{\omega}_{i,j}^{\eta,\underline{0}}$ is obtained by $\bar{\omega}_i^{\eta,\underline{0}}$ by erasing $j$ sites:
\begin{equation}
\tilde{\omega}_{i,j}^{\eta,\underline{0}}=\bar{\omega}_i^{\eta,\underline{0}}\setminus\{(x,y): \ x=n+l_1-i, \ y\in(m+l_2-j, m+l_2]\}.
\end{equation}
Again, as in the reference path $\omega^*$, the last interpolation consists in inserting between every pair of consecutive configurations in $\tilde{\omega}^{\eta,\underline{0}}$ a sequence of configurations with a free particle in a suitable sequence of sites going from the site previously occupied by the erased particle to $\partial\Lambda$. If $l_1\leq2l_2^*-3$ and $l_2\leq l_2^*$, we have $H(\mathcal{R}(l_1,l_2))\leq H(\mathcal{R}(2l_2^*-3, l_2^*))$. In our strong anisotropic case $U_1+U_2+\e(l_2^*-3)<2\Delta-U_1$, so for the path $\omega^{\eta,\underline{0}}$ obtained in this way we have
\begin{equation}
\max_{i}H(\omega_i^{\eta,\underline{0}})=\max_{l\in{[1,\ l_1]}}H(\mathcal{R}(l, l_2))+U_1+U_2+\e(l-l_2)\leq H(\mathcal{R}(2l_2^*-3, l_2^*))+U_1+U_2+\e(l_2^*-3)<\Gamma.
\end{equation}
If $l_1\leq 2l_2^*-1$ and $l_2\leq l_2^*-1$, we have $H(\mathcal{R}(l_1,l_2))\leq H(\mathcal{R}(2l_2^*-1, l_2^*-1))$. In our strong anisotropic case $U_1+U_2+\e(l_2^*-3)<\D-U_2+U_1$, so for the path $\omega^{\eta,\underline{0}}$ obtained in this way we have
\begin{equation}
\max_{i}H(\omega_i^{\eta,\underline{0}})=\max_{l\in{[1,\ l_1]}}H(\mathcal{R}(l, l_2))+U_1+U_2+\e(l_2-2)\leq H(\mathcal{R}(2l_2^*-1, l_2^*-1))+U_1+U_2+\e(l_2^*-3)<\Gamma.
\end{equation}

The proof of $2$ is similar. Let $\mathcal{C}_{\underline{1}}$ be the maximal connected set containing $\pieno$ such that 

\noindent 
$\max_{\h'\in{\cC_{\pieno}}}H(\h')<\G$. Again $\mathcal{C}_{\underline{1}}$ is by definition a cycle containing $\underline{1}$ and not containing $\underline{0}$. To prove that $\mathcal{C}_{\underline{1}}$ contains $\eta$ we define now a path $\omega^{\eta,\underline{1}}$ going from $\eta$ to $\underline{1}$ obtained by reaching rectangular configurations with $l_2=l_2^*$ or $l_1\geq L-1$ and, from there, following the path $\omega^*$. As before, it is easy to prove that all the configurations of this path have an energy smaller than $\Gamma$, so they are in $\mathcal{C}_{\underline{1}}$.

Going into details, let $\eta\in{\mathcal{R}(l_1, l_2)}$. First, we consider $3l_2^*-2\leq s\leq l_2^*+L-1$. If $l_2=l_2^*$ we can choose $\omega^{\eta,\underline{1}}$ as the part of the reference path $\omega^*$ going from $\eta$ to $\underline{1}$. If $l_2<l_2^*$ then first move columns to rows, until we obtain $l_2=l_2^*$. From there, follow the reference path $\omega^*$. If $l_2>l_2^*$ then first add columns until we reach $l_1=L-1$. The remaining part of the path follows $\omega^*$. Now we consider the case $s\geq l_2^*+L-1$. If $l_1\geq L-1$ we can choose $\omega^{\eta,\underline{1}}$ as the part of the reference path $\omega^*$ going from $\eta$ to $\underline{1}$. If $l_1<L-1$ we add columns until we have $l_1=L-1$ and then we follow the reference path $\omega^*$. For the path $\omega^{\eta,\underline{1}}$ obtained in this way, since $H(\mathcal{R}(l_1', l_2'))<H(\mathcal{R}(2l_2^*-1, l_2^*-1))$ for any $l_1'\geq 2l_2^*-2$ and $l_2'\geq l_2^*$, we obtain
\begin{equation}
\max_{i}H(\omega_i^{\eta,\underline{1}})\leq\max_{l_1'\geq 2l_2^*-2,l_2'\geq l_2^*}H(\mathcal{R}(l_1', l_2'))+2\Delta-U_1<\Gamma,
\end{equation}
so that $\omega_{i}^{\eta,\underline{1}}\in{\mathcal{C}_{\underline{1}}}$ for any $i$ and the proof of the theorem is complete.

\appendix
\section{Appendix}
We give explicit argument to complete the proof of Proposition \ref{lemma3}, considering the cases that were left in subsection \ref{dim3}, because the proofs are analogue to the ones discussed in that subsection.

\medskip
\noindent
{\bf Additional material for Subsection \ref{dim3}}

\medskip
\noindent
{\bf Case $\D s=3$.} We report the explicit computations for the case (III). Since $\D s=3$ by Remark \ref{remark2} we know that $s(\bar\h)\geq s^*-4$, but we consider only the case $s(\bar\h)\geq s^*-3$ (see Lemma \ref{lemma3}). If $p_2(\bar\h)=l_2^*-2$ the circumscribed rectangle of $\bar\h$ is $\cR(2l_2^*+k-3,l_2^*-2)$, for any $k\geq1$, and if $p_2(\bar\h)=l_2^*-1$ it is $\cR(2l_2^*+k-4,l_2^*-1)$, for any $k\geq1$. Since in the case we are analyzing there are three horizontal lines becoming active, we know that $n(\bar\h)\geq3$. If $p_2(\bar\h)=l_2^*-2$, since $k\geq1$, $\d<1$ and $\e\ll U_2$, we get $H(\bar\h)>\G$. Indeed
\be{}
\ba{lll}
H(\bar\h)\geq H(\cR(2l_2^*+k-3,l_2^*-2))+3\D=\\

U_1l_2^*+2U_2l_2^*+U_1+kU_2-\e(2(l_2^*)^2+kl_2^*-2k-7l_2^*+6)-3\e>\G\Leftrightarrow\\

\Leftrightarrow \ 5U_2>\e(7-4\d+k(\d-2)).
\ea
\ee

If $p_2(\bar\h)=l_2^*-1$, since $k\geq1$, $\d<1$ and $\e\ll U_2$, we get $H(\bar\h)>\G$. Indeed
\be{}
\ba{lll}
H(\bar\h)\geq H(\cR(2l_2^*+k-4,l_2^*-1))+3\D=\\

=U_1l_2^*+2U_2l_2^*+2U_1+kU_2-U_2-\e(2(l_2^*)^2+kl_2^*-k-6l_2^*+7)>\G\Leftrightarrow\\

\Leftrightarrow \ U_1+3U_2>\e(5-3\d+k(\d-1)).
\ea
\ee

\medskip
\noindent
{\bf Case $\D s=4$.} We report the explicit computations for the cases (I) and (II). Since $\D s=4$, by Remark \ref{remark2} we know that $s(\bar\h)\geq s^*-5$, but we consider only the case $s(\bar\h)\geq s^*-4$ (see Lemma \ref{lemma3}). If $p_2(\bar\h)=l_2^*-2$ the circumscribed rectangle of $\bar\h$ is $\cR(2l_2^*+k-4,l_2^*-2)$ with $k\geq1$ and if $p_2(\bar\h)=l_2^*-1$ it is $\cR(2l_2^*+k-5,l_2^*-1)$, for any $k\geq1$. Since in the case we are analyzing there are three horizontal lines becoming active, we know that $n(\bar\h)\geq3$. For the case (I), if $p_2(\bar\h)=l_2^*-2$, since $k\geq1$, $\d<1$ and $\e\ll U_2$, we get $H(\bar\h)>\G$. Indeed
\be{}
\ba{lll}
H(\bar\h)\geq H(\cR(2l_2^*+k-4,l_2^*-2))+3\D=\\

=U_1l_2^*+2U_2l_2^*+U_1-U_2+kU_2-\e(2(l_2^*)^2+kl_2^*-2k-8l_2^*+11)>\G\Leftrightarrow\\

\Leftrightarrow \ 5U_2>\e(9-5\d+k(\d-2)).
\ea
\ee
If $p_2(\bar\h)=l_2^*-1$, since $k\geq1$, $\d<1$ and $\e\ll U_2$, we get $H(\bar\h)>\G$. Indeed
\be{}
\ba{lll}
H(\bar\h)\geq H(\cR(2l_2^*+k-5,l_2^*-1))+3\D=\\

U_1l_2^*+2U_2l_2^*+2U_1-2U_2+kU_2-\e(2(l_2^*)^2+kl_2^*-k-7l_2^*+8)>\G\Leftrightarrow\\

\Leftrightarrow \ U_1+3U_2>\e(6-4\d+k(\d-1)).
\ea
\ee

For the case (II), since if $l_2^*-3\leq p_2(\bar\h)\leq l_2^*-1$ the only possibilities that we can have are $p_2(\h)\leq l_2^*-1$ or $p_2(\h)=l_2^*$, the cases that remain to analyze in detail are the followings:

\bi

\item $p_2(\bar\h)=l_2^*-2$ and $p_2(\h)>l_2^*$;
\item $p_2(\bar\h)=l_2^*-1$ and $p_2(\h)>l_2^*$.
\ei
\noindent
For the cases in which $p_2(\bar\h)=l_2^*-2$ and $p_2(\bar\h)=l_2^*-1$, the computations are exactly the same reported for the case (I).

\medskip
\noindent
{\bf Case $\D s=5$.} We report the explicit computations. Since $\D s=5$, by Remark \ref{remark2} we know that $s(\bar\h)\geq s^*-6$, but we consider only the case $s(\bar\h)\geq s^*-5$ (see Lemma \ref{lemma3}). The cases that remain to analyze in detail are the followings:

\bi
\item $p_2(\bar\h)=l_2^*-2$;
\item $p_2(\bar\h)=l_2^*-1$.
\ei

If $p_2(\bar\h)=l_2^*-2$, the circumscribed rectangle of $\bar\h$ is $\cR(2l_2^*+k-5,l_2^*-2)$, for any $k\geq1$, and if $p_2(\bar\h)=l_2^*-1$ it is $\cR(2l_2^*+k-6,l_2^*-1)$, for any $k\geq1$. By Lemma \ref{lemma3.2} we know that $n(\bar\h)\geq4$. If $p_2(\bar\h)=l_2^*-2$, since $k\geq1$, $\d<1$ and $\e\ll U_2$, we get $H(\bar\h)>\G$. Indeed
\be{}
\ba{lll}
H(\bar\h)\geq H(\cR(2l_2^*+k-5,l_2^*-2))+4\D=\\

U_1l_2^*+2U_2l_2^*+2U_1-U_2+kU_2-\e(2(l_2^*)^2+kl_2^*-2k-7l_2^*+14)>\G\Leftrightarrow\\

\Leftrightarrow \ U_1+4U_2>\e(12-4\d+k(\d-2)).
\ea
\ee

If $p_2(\bar\h)=l_2^*-1$, since $k\geq1$, $\d<1$ and $\e\ll U_2$, we get $H(\bar\h)>\G$. Indeed
\be{}
\ba{lll}
H(\bar\h)\geq H(\cR(2l_2^*+k-6,l_2^*-1))+4\D=\\

=U_1l_2^*+2U_2l_2^*+3U_1-2U_2+kU_2-\e(2(l_2^*)^2+kl_2^*-k-8l_2^*+10)>\G\Leftrightarrow\\

\Leftrightarrow \ 2U_1+4U_2>\e(8-5\d+k(\d-1)).
\ea
\ee

\medskip
\noindent
We give explicit argument to complete the proof of Proposition \ref{lemma4}, considering the cases that were left in subsection \ref{dim4}, because the proofs are analogue to the ones discussed in that subsection.

\medskip
\noindent
{\bf Additional material for Subsection \ref{dim4}}

\medskip
\noindent
{\bf Case $\D$s=4.} We report the explicit computations for the cases (b) and (c). In the case (b) we obtain
\be{}
\ba{lll}
H(\bar\h)\geq H(\cR(2l_2^*-k-x,l_2^*+k))+3\D=\\

=U_1l_2^*+2U_2l_2^*+k(U_1-U_2)+3U_1+3U_2-xU_2-2\e(l_2^*)^2-k\e l_2^*+\e k^2+x\e l_2^*+k\e x-3\e>\G\\

\Leftrightarrow \ \e k^2+k[U_1-U_2-\e l_2^*+x\e]+2U_1+4U_2-xU_2-\e+x\e l_2^*-3\e l_2^*>0.
\ea
\ee

\noindent
By (\ref{5.10}), since $k\geq0$, $0<\d<1$, $x\geq1$ and $\e\ll U_2$, we get $H(\bar{\h})>\G$. Indeed
\be{}
2U_1+4U_2-xU_2-\e+x\e l_2^*-3\e l_2^*=2U_1+U_2-\e-3\e\d+x\e\d\gg2U_1-2\e\d>0.
\ee

In the case (c), we get $H(\bar\h)>\G$. Indeed
\be{}
\ba{lll}
H(\cR(2l_2^*,l_2^*))+\e(2l_2^*-x-1)+3\D=\\

=U_1l_2^*+2U_2l_2^*-\e(2(l_2^*)^2-xl_2^*)+2\e l_2^*+3U_1+3U_2-xU_2-4\e-x\e>\G\Leftrightarrow\\

\Leftrightarrow \ 3U_2+2U_1>\e(2+\d+x(1-\d)), \quad \hbox{always since }  x\leq1, 1-\d>0 \hbox{ and } \e\ll U_2.
\ea
\ee

\medskip
\noindent
{\bf Case $\D s=5$.} We report the explicit computations for the cases (b) and (c). In the case (b) we obtain
\be{}
\ba{lll}
H(\bar\h)\geq H(\cR(2l_2^*-k-x,l_2^*+k))+4\D=\\

=U_1l_2^*+2U_2l_2^*+k(U_1-U_2)+4U_1+4U_2-xU_2-2\e(l_2^*)^2-k\e l_2^*+\e k^2+x\e l_2^*+k\e x-4\e>\G\\

\Leftrightarrow \ \e k^2+k[U_1-U_2-\e l_2^*+x\e]+3U_1+5U_2-xU_2-2\e+x\e l_2^*-3\e l_2^*>0.
\ea
\ee

\noindent
By (\ref{5.10}), since $k\geq0$, $0<\d<1$, $x\geq1$ and $\e\ll U_2$., we get $H(\bar{\h})>\G$. Indeed
\be{}
3U_1+5U_2-xU_2-2\e+x\e l_2^*-3\e l_2^*=3U_1+2U_2-2\e-3\e\d+x\e\d\gg3U_1-2U_2>0
\ee

In the case (c), since $x\leq1$, $1-\d>0$ and $\e\ll U_2$, we get $H(\bar\h)>\G$. Indeed
\be{}
\ba{lll}
H(\bar\h)\geq H(\cR(2l_2^*-x,l_2^*))+\e(2l_2^*-x-1)+4\D=\\

=U_1l_2^*+2U_2l_2^*+4U_1+4U_2-xU_2-\e(2(l_2^*)^2-xl_2^*)+2\e l_2^*-x\e-5\e>\G\Leftrightarrow\\

\Leftrightarrow \ 3U_1+4U_2>\e(3+\d+x(1-\d)).
\ea
\ee

\medskip
\noindent
We give explicit argument to complete the proof of Proposition \ref{lemma6}, considering the cases that were left in subsection \ref{dim6}, because the proofs are analogue to the ones discussed in that subsection.

\medskip
\noindent
{\bf Additional material for Subsection \ref{dim6}}

Let $\D s=0$. We analyze in detail the cases $(b-i)$ and $(b-ii)$. In the case $(b-i)$ we have $s(\bar{\h})\geq s^*$, $\D v\leq-1$ and $n(\bar{\h})\geq1$. By definition (\ref{defB}) for the case $s\geq s^*$, we get $p_2(\bar\h)=l_2^*$ and $v(\bar{\h})\geq p_{max}(\bar{\h})-1$. Thus the circumscribed rectangle of $\bar{\h}$ is $\cR(2l_2^*+k-1,l_2^*)$, for any $k\geq0$. Since $n(\bar{\h})\geq 1$, from (\ref{energeta}) we obtain
\be{}
\ba{ll}
H(\bar{\h})\geq H(\cR(2l_2^*+k-1,l_2^*))+\e v(\bar{\h})+\D\geq\\

\geq U_1l_2^*+2U_2l_2^*+kU_2-U_2-\e (2(l_2^*)^2+kl_2^*-l_2^*)+\e (2l_2^*+k-2)+U_1+U_2-\e.
\ea
\ee

We recall that $\G=U_1l_2^*+2U_2l_2^*+U_1-U_2-2\e (l_2^*)^2+3\e l_2^*-2\e$. Thus we get $H(\bar{\h})>\G$ if and only if $U_2>\e(1-k(1-\d))$, always since $\e\ll U_2$, $k\geq0$ and $\d<1$.

In the case $(b-ii)$ we have $s(\bar{\h})\geq s^*$, $\D v\leq-1$ and $n(\bar{\h})=0$, so from Remark \ref{remark1} it follows that $g_1'(\bar{\h})+g_2'(\bar{\h})\geq1$. We consider the following four cases:

\begin{itemize}
	\item [A.] $g_2'(\bar{\h})=1$;
	\item [B.] $g_1'(\bar{\h})=1$;
	\item [C.] either $g_1'(\bar{\h})=1$ and $g_2'(\bar{\h})=1$, or $g_2'(\bar\h)\geq2$;
	\item [D.] $g_1'(\bar\h)\geq2$ and $g_2'(\bar{\h})=0$.
\end{itemize}

\noindent
{\bf Case A.} Since $g_2'(\bar{\h})=1$, since $k\geq0$ and $\d<1$, we get $H(\bar\h)\geq\G$. Indeed
\be{}
\ba{lll}
H(\bar{\h})=H(\cR(2l_2^*+k-1,l_2^*))+\e v(\bar{\h})+U_1\geq\\

\geq U_1l_2^*+2U_2l_2^*+kU_2-U_2-\e (2(l_2^*)^2+kl_2^*-l_2^*)+\e (2l_2^*+k-2)+U_1\geq\G\Leftrightarrow\\

\Leftrightarrow \ k\e(1-\d)\geq0.
\ea
\ee

We note that $H(\bar{\h})=\G$ if $k=0$ and for those configurations $\bar{\h}$ such that $g_2'(\bar{\h})=1$ and $v(\bar{\h})=2l_2^*-1=p_{max}(\bar{\h})-1$, i.e., $\bar\h\in\cP_1$. Starting from such $\bar\h$, we note that in order to get $\D s=0$ the only admissible transitions are the movement of a single protuberance along the same side. In this way $\h\in\cP_1\subset\cB$. This contributes to Theorem \ref{3.1}{\it (ii)}. If $k\geq1$ or $v(\bar{\h})>p_{max}(\bar{\h})-1$, we get $H(\bar{\h})>\G$. 

\medskip
\noindent
{\bf Case B.} We have $g_1'(\bar{\h})=1$. By Remark \ref{remarkulteriore} we know that no line can become active, so $\D s=0$ is obtained by no line becoming active nor inactive. Referring to Figure \ref{fig:g1'}, if $\bar\h_{cl}$ is connected we note that the only admissible operations are moving protuberances or let a particle become free in $\h$ in such a way that $\D s=0$. In both cases we get $v(\h)\geq v(\bar\h)$. If $\bar\h\in\cB$, we get $\h\in\cB$, since $s(\h)=s(\bar\h)$ and due to the condition about the number of vacancies. If $\bar\h\notin\cB$ it is not a relevant case. If $\bar\h_{cl}$ is not connected we can argue similarly.

\medskip
\noindent
{\bf Case C.} We have either $g_1'(\bar\h)=)1$ and $g_2'(\bar\h)=1$, or $g_2'(\bar\h)\geq2$. Thus we can conclude in the same way as in the case II-C for $\D s=-1$.

\medskip
\noindent
{\bf Case D.} We have $g_1'(\bar{\h})\geq 2$ and $g_2'(\bar{\h})=0$, thus we can conclude in the same way as in the case II-B for $\D s=-1$.

\medskip
\noindent
We give explicit argument to complete the proof of Proposition \ref{lemma7}, considering the cases that were left in subsection \ref{dim7}, because the proofs are analogue to the ones discussed in that subsection.

\medskip
\noindent
{\bf Additional material for Subsection \ref{dim7}}

Let $\D s=1$. We analyze in detail the cases (b-i), (c-i) and (c-ii). In the case (b-i), the circumscribed rectangle of $\bar{\h}$ is $\cR(2l_2^*-k-2,l_2^*+k)$, for any $k\geq0$. Thus, since $\frac{U_1-2U_2}{\e}>-1$ and $\d-k-3<-2$, we get $H(\bar{\h})\geq\G$. Indeed
\be{}
\ba{lll}
H(\bar{\h})\geq H(\cR(2l_2^*-k-2,l_2^*+k))+\e(l_2^*+k-1)+\D=\\

=U_1l_2^*+2U_2l_2^*+U_1-U_2+k(U_1-U_2)-2\e(l_2^*)^2-\e kl_2^*+\e k^2+3\e l_2^*+3\e k-2\e>\G\\

\Leftrightarrow\e k^2+k[U_1-U_2+3\e-\e l_2^*]>0\Leftrightarrow \  \frac{U_1-2U_2}{\e}>\d-k-3.
\ea
\ee
\noindent
In particular, we obtain $H(\bar{\h})=\G$ if $k=0$ and $v(\bar{\h})=l_2^*-1$, i.e., $\bar\h\in\cP_2$, otherwise $H(\bar{\h})>\G.$ For such $\bar\h\in\cP_2$, we note that in order to get $\D s=1$, the only admissible operations are to attach the free particle to the protuberance, or to one of the three other sides. If we attach the particle on the vertical side, we get $v(\bar\h)+l_2^*-1=2l_2^*-2$, with $s(\h)=s^*$ and $p_2(\h)=l_2^*$, so $\h\in\cB$ and thus this is not a relevant case. If we attach the free particle on the horizontal side, we get $p_2(\h)=l_2^*+1$, with $s(\h)=s^*$, and thus $\h\notin\cB$. In this case we obtain a configuration in $\cP_2$. 

For the cases (c-i) and (c-ii) we refer again to the case $\D s=-1$. In the case (c-i), since $n(\bar{\h})\geq 1$, from (\ref{energeta}) we obtain

\be{}
\ba{lll}
H(\bar{\h})\geq H(\cR(2l_2^*-1,l_2^*))+\e v(\bar{\h})+\D\geq\\

\geq U_1l_2^*+2U_2l_2^*-U_2-\e (2(l_2^*)^2-l_2^*)+\e (2l_2^*-2)+U_1+U_2-\e=\\

=U_1l_2^*+2U_2l_2^*-2\e (l_2^*)^2+3\e l_2^*-3\e +U_1.
\ea
\ee

\noindent
Thus we get $H(\bar{\h})>\G$ if and only if $-3\e>-U_2-2\e\Leftrightarrow \ \e<U_2$, always since $\e\ll U_2$.

For the case (c-ii), again we distinguish the following subcases:

\begin{enumerate}
	\item [A.] $g_2'(\bar{\h})=1$;
	\item [B.] $g_1'(\bar{\h})=1$;
	\item [C.] either $g_1'(\bar{\h})=1$ and $g_2'(\bar{\h})=1$, or $g_2'(\bar\h)\geq2$;
	\item [D.] $g_1'(\bar{\h})\geq 2$ and $g_2'(\bar{\h})=0$.
\end{enumerate}

The reasonings are the same of the case $\D s=0$ shown in this Appendix.

\medskip
\noindent
We give explicit argument to complete the proof of Proposition \ref{lemma8}, considering the cases that were left in subsection \ref{dim8}, because the proofs are analogue to the ones discussed in that subsection.

\medskip
\noindent
{\bf Additional material for Subsection \ref{dim8}}

Let $\D s=2$. We analyze in detail the cases (c) and (d). For the case (c), since $\frac{U_1-2U_2}{\e}>-1$ and $\d-k-3<-2$, we get $H(\bar{\h})\geq\G$. Indeed
\be{}
\ba{llll}
H(\bar{\h})\geq H(\cR(2l_2^*-k-2,l_2^*+k))+\e(l_2^*+k-1)+\D=\\

=U_1l_2^*+2U_2l_2^*+U_1-U_2+k(U_1-U_2)-2\e(l_2^*)^2-\e kl_2^*+\e k^2+3\e l_2^*+3\e k-2\e>\G\\

\Leftrightarrow\e k^2+k[U_1-U_2+3\e-\e l_2^*]>0\Leftrightarrow \ \frac{U_1-2U_2}{\e}>\d-k-3.
\ea
\ee

\noindent
In particular, we obtain $H(\bar{\h})=\G$ if $k=0$ and $v(\bar{\h})=l_2^*-1$, i.e., $\bar\h\in\cP_2$, otherwise $H(\bar{\h})>\G.$ For such $\bar\h\in\cP_2$, there is no admissible exiting move from $\cB$ in order to get $\D s=2$, thus it is not a relevant case. 

In the case (d), since $n(\bar{\h})\geq1$, we get $H(\bar\h)>\G$. Indeed
\be{}
\ba{lll}
H(\bar{\h})\geq H(\cR(2l_2^*+k-1,l_2^*))+\e v(\bar{\h})+\D\geq\\

\geq U_1l_2^*+2U_2l_2^*+kU_2-U_2-\e (2(l_2^*)^2+kl_2^*-l_2^*)+\e (2l_2^*+k-2)+U_1+U_2-\e>\G\\

\Leftrightarrow \ U_2>\e(1-k(1-\d)), \quad \hbox{always since } \e\ll U_2, \ k\geq0 \hbox{ and } \d<1.
\ea
\ee

\end{document}